# Characters of locally analytic representations of $p$-adic reductive groups

Ralf Diepholz

**Abstract.** We propose a definition of characters in the context of Schneider-Teitelbaum's theory of locally analytic representations of $p$-adic reductive groups. This character will be a function on a compact subgroup $S_0$ in a maximal torus of the reductive group $G$. As an elementary tool we study evaluability of formal characters of commutative group representations. To an admissible $G$-representation $V$ we then associate a projective system of $S_0$-representations; the design of this approximating system constitutes the technical heart of this work. If the components of the projective system possess evaluable formal characters, and if their values converge to a function on $S_0$, then this limit function is the desired character of $V$. We show that our definition generalizes the one established in classical smooth representation theory. We determine the character in two examples.

**Contents**

*Ut! Olicrosse! Godemite!*
*(J. Neukirch, Algebraische Zahlentheorie, p. 46)*

## Introduction

In the smooth representation theory of $p$-adic groups one associates to every admissible representation $V$ of a group $G$ a trace $\Theta_V$ which is a linear form on the Hecke algebra $\mathcal{H}(G)$ of $G$. A theorem of Howe and Harish-Chandra states that if $G$ is reductive and $V$ is finitely generated as a $G$-module then the trace "coincides" with a smooth function $\theta$ (the so-called *character* of $V$) on the regular set $G^{\mathrm{reg}}$ of $G$; this means that $\Theta_V(f) = \int \theta(x)\, f(x)\, dx$ for every compactly supported smooth function $f$ on $G^{\mathrm{reg}}$.

The main obstacle to introducing the concept of character to the locally analytic representation theory of Schneider-Teitelbaum is the following fact: While on any admissible smooth $G$-representation the Hecke algebra $\mathcal{H}(G)$ acts by finite rank operators, on an admissible locally analytic $G$-representation the distribution algebra $D(G)$ of $G$ (which takes the role of the Hecke algebra in the locally analytic setting) acts by operators which in general are not even nuclear. Therefore there is no obvious way to assign a trace $\Theta(\lambda)$ to the operators $\lambda \in D(G)$ in order to obtain a linear form $\Theta$ on $D(G)$ which would - by reflexivity - point a way to define a character function on the group. In this work we propose a definition of what it means for an admissible locally analytic representations of a $p$-adic reductive group $G$ to possess a character on a subset $S'$ of a compact subgroup of a maximal torus of $G$. This definition will depend on the choice of a minimal parabolic subgroup.

In the following we explain our proposal, thereby giving an outline of the sections of this work. We first develop a theory of formal characters for representations of commutative groups which do not necessarily act by finite rank (or nuclear) operators but satisfy a certain finiteness condition on weight multiplicities, a condition we call *finite trigonalisability* (Sections 1-3). To such a finitely trigonalisable representation $V$ of a group $S$ over a field $K$ we associate a *formal character* $\Theta_V$ which is a "formal function", i.e. a formal sum of $K$-valued functions on $S$. Given a subset $S' \subset S$, we study the conditions under which a formal function $f$ on $S$ may be evaluated on $S'$, thus yielding an actual function $\mathrm{ev}_{S'}(f)$ on $S'$. If the formal character $\Theta_V$ is evaluable on $S'$ in this sense then the function

$$\mathrm{ev}_{S'}(\Theta_V) : S' \to K$$

is called the *character on* $S'$ of $V$.

Let $p$ be a prime number, and let $K\,|\,L\,|\,\mathbb{Q}_p$ be a tower of complete valued fields such that $L\,|\,\mathbb{Q}_p$ is finite and $K$ is discretely valued. Let $G$ be (the group of $L$-rational points of) a connected reductive group over $L$, $S \subset G$ a maximal split torus over $L$, and $N$ the unipotent radical of a minimal parabolic subgroup containing $S$, with Lie algebra $\mathfrak{n}$. Let $V$ be an admissible locally analytic representation of $G$ over $K$. Sections 4 to 8 prepare the definition of

an $\mathfrak{n}$-*character* of $V$, to be given in Section 9. (The reader who is familiar with the locally analytic representation theory, in particular with [32], is encouraged to start reading Section 9 right away and to use the earlier sections as a reference guide.) To this end we carry out a two-folded approximation process. Let $S_0$ be a compact subgroup of the torus $S$, and choose a compact open subgroup $G_0 \subset G$ containing $S_0$. We view $V$ only as a representation of the group $G_0$. By a result of Schneider-Teitelbaum $V$ is the compact inductive limit of a system $(V_i)_{i\in\mathbb{N}}$ of locally analytic $G_0$-representations on $K$-Banach spaces ([32], Proposition 6.5). In a second step, one can use the associated Lie algebra action to form the subspaces $V_i^{(\mathfrak{n})} \subset V_i$ of vectors which are $\mathfrak{n}$-*adically finite*, i.e. which are annihilated by some power of $\mathfrak{n}$ in the enveloping algebra $U(\mathfrak{n})$. The spaces $V_i^{(\mathfrak{n})}$ are stable under the compact commutative group $S_0$. Now we can formulate a first definition.

**Definition (preliminary version).** Let $S' \subset S_0$ be a subset. Suppose that each $S_0$- representation $V_i^{(\mathfrak{n})}$ possesses a character $\theta_i$ on $S'$ and that the sequence of functions $(\theta_i)_{i\in\mathbb{N}}$ converges pointwise. Then the limit function $\theta = \lim \theta_i$ is called the $\mathfrak{n}$-*character on* $S'$ of $V$.

The problem with this preliminary definition is that the limit of the above sequence of characters on $S'$ behaves not very well with respect to a change in the choice of the compact inductive system $(V_i)$ (for an illustration of this phenomenon compare the two examples in Section 3). As a consequence we have to specify a "canonical" compact inductive system giving $V$ in the limit. It turns out that this is easier to do in the dual setting.

Let $D(G_0, K)$ be the distribution algebra of $G_0$; this is a noetherian Fréchet algebra over $K$. Associated to any open subgroup $H \subset G_0$ which is a uniform pro- $p$-group there is a family $\|\ \|_{(H,r)}$ ($r \in p^{\mathbb{Q}}$, $p^{-1} \le r < 1$) of Banach algebra norms on $D(G_0, K)$ such that the completions $D_{(H,r)}(G_0, K)$ with respect to these norms realize the structure of a Fréchet-Stein algebra of $D(G_0, K)$; this implies in particular that $D(G_0, K) = \varprojlim_r D_{(H,r)}(G_0, K)$ as Fréchet spaces. Moreover, by admissibility the strong dual $M = V_b'$ of our $G_0$-representation is equal to the projective limit

$$M = \varprojlim_r M_{(H,r)}$$

where $M_{(H,r)} := D_{(H,r)}(G_0, K) \otimes_{D(G_0,K)} M$. In order to obtain a "canonical" projective system - not depending on the choice of $H$ - we consider all uniform open subgroups (satisfying certain technical assumptions) of $G_0$ at the same time. More exactly, we view $M$ as the projective limit of the system

$$(M_{(H,r)})_{(H,r)}$$

with a suitable directed ordering on the set of pairs $(H, r)$. A large part of Sections 5- 7 is devoted to the investigation of the transition maps of this projective system; here an elementary but important result is a kind of elementary divisor theorem for uniform pro- $p$-groups. The second approximation step, corresponding to the passage to $\mathfrak{n}$-adically finite vectors described above, is to consider the system of quotient representations $M_{(H,r)} / \mathfrak{n}^k M_{(H,r)}$ $(k \in \mathbb{N})$ of $S_0$.

**Definition (final version).** Set $S' \subset S_0$ be a subset. Suppose that

- each $S_0$-representation $M_{(H,r)} / \mathfrak{n}^k M_{(H,r)}$ possesses a formal character $\Theta_{(H,r),k}$;
- for each pair $(H, r)$ the sequence $(\Theta_{(H,r),k})_{k \in \mathbb{N}}$ of formal functions converges to a formal function $\Theta_{(H,r)}$ which is evaluable on $S'$;
- the net of evaluations $(\mathrm{ev}_{S'}(\Theta_{(H,r)}))_{(H,r)}$ converges pointwise to a function $\theta$ on $S'$.

Then the function $s \mapsto \theta(s^{-1})$ is called the $\mathfrak{n}$- *character on* $S'$ of $V$.

It turns out that this definition does not depend on the choice of the compact open subgroup $G_0 \subset G$ containing $S_0$. We emphasize that we do not claim that every admissible representation possesses an $\mathfrak{n}$-character.

The theory of admissible locally analytic representations includes as a special case the theory of admissible smooth representations. Section 10 is devoted to the verification that our $\mathfrak{n}$-characters generalize the classical concept of characters in the smooth theory.

We finally calculate the $\mathfrak{n}$-character in two examples of genuine locally analytic representations: The locally analytic principal series of $\mathrm{SL}_2(\mathbb{Q}_p)$ (Section 12) and of the Iwahori subgroup in a split reductive $p$-adic group (Section 13). These sections may be read independently of each other; they rely on results about explicit Banach space bases and filtrations of completed distribution algebras developped in Sections 7 and 8, and on general results about principal series representations (Section 11).

This article is based on the doctoral thesis of the author, who wishes to thank his advisor Peter Schneider for good advice, his colleagues Jan Kohlhaase, Enno Nagel, Tobias Schmidt and Matthias Strauch for a few helpful discussions, his friends and fellows Alexis Pangalos and Stefan Wiech who accompanied him from Hamburg to Münster, and the group *Sport* for ***Homerun.***

# I. Formal characters

This part is completely general. Simple results in infinite-dimensional linear algebra are given, mainly for lack of reference.

## 1. Evaluation of formal sums of functions

Let $S$ be a set, and let $R$ be an integral domain. The ring structure on $R$ induces a ring structure on the set $\mathrm{Map}(S, R)$ of all mappings $S \to R$ and a group structure on the set $\mathrm{Map}(S, R^\times)$ of all mappings $S \to R^\times$. Let $X \subset \mathrm{Map}(S, R^\times)$ be a subgroup. Let $\mathbb{Z}[\![X]\!]$ denote the $\mathbb{Z}$-module of all mappings $X \to \mathbb{Z}$, written as formal sums

$$\textstyle\sum_{\lambda \in X} n_\lambda \, e(\lambda) \quad (n_\lambda \in \mathbb{Z}).$$

We endow $\mathbb{Z}[\![X]\!]$ with the product topology of the discrete groups $\mathbb{Z}$. Thus $\mathbb{Z}[\![X]\!]$ is a Hausdorff and complete abelian group and contains the group ring $\mathbb{Z}[X]$ of $X$ as a dense subgroup (cf. [6]; Chap. 2, Sect. 3, Proposition 10; Chap. 3, Sect. 2, Proposition 25). This topology has the following features:

• Let $I$ be a directed set. A net $(f_i = \sum_\lambda n_{i,\lambda}\, e(\lambda))_{i \in I}$ in $\mathbb{Z}[\![X]\!]$ converges to an element $f = \sum_\lambda n_\lambda \, e(\lambda) \in \mathbb{Z}[\![X]\!]$ if and only if for every $\lambda \in X$ there exists an index $i \in I$ such that $n_{j,\lambda} = n_\lambda$ for all $j \geq i$.

• Given an arbitrary set $I$, a family $(f_i = \sum_\lambda n_{i,\lambda}\, e(\lambda))_{i \in I}$ in $\mathbb{Z}[\![X]\!]$ is summable if and only if for each $\lambda \in X$ the family $(n_{i,\lambda})_{i \in I}$ in $\mathbb{Z}$ has finite support; in this case the sum of the family $(f_i)_{i \in I}$ is given by

$$\textstyle\sum_i f_i = \sum_\lambda \left(\sum_i n_{i,\lambda}\right) e(\lambda) \in \mathbb{Z}[\![X]\!]$$

([5], Chap. IV, §4, No. 2, Lemma 1).

We will call two elements $f = \sum_\lambda n_\lambda \, e(\lambda)$, $g = \sum_\lambda m_\lambda \, e(\lambda) \in \mathbb{Z}[\![X]\!]$ *multipliable* if the family $\left(\sum_{\lambda \in X} n_\mu \, m_{\lambda \mu^{-1}} \, e(\lambda)\right)_{\mu \in X}$ is summable, i.e. if for each $\lambda \in X$ there are but finitely many pairs $(\mu_1, \mu_2) \in X \times X$ such that $\mu_1 \mu_2 = \lambda$ and $n_{\mu_1} m_{\mu_2} \neq 0$. If this is the case then we call the element

$$f\, g := \textstyle\sum_\lambda \left(\sum_\mu n_\mu \, m_{\lambda \mu^{-1}}\right) e(\lambda) \in \mathbb{Z}[\![X]\!]$$

the *product* of $f$ and $g$. This definition restricts to the usual multiplication on the group ring $\mathbb{Z}[X] \subset \mathbb{Z}[\![X]\!]$; moreover, it endows $\mathbb{Z}[\![X]\!]$ with the structure of a $\mathbb{Z}[X]$- module.

Let $S' \subset S$ be a subset. The ring homomorphism

$$\begin{array}{rccc} \mathrm{ev}_{S'}: & \mathbb{Z}[X] & \longrightarrow & \mathrm{Map}(S', R) \\ & \sum_\lambda n_\lambda\, e(\lambda) & \longmapsto & \sum_\lambda n_\lambda\, \lambda\,|_{S'} \end{array}$$

may be extended "by taking quotients" in the following way: Since $R$ is an integral domain the subset $\mathcal{S}_{S'} \subset \mathbb{Z}[X]$ of all $h$ such that $\mathrm{ev}_{S'}(h)$ has no zeros is multiplicative. Moreover, $\mathrm{ev}_{S'}(h)$ is invertible in the ring of mappings $S' \to \mathrm{Quot}\, R$ for any such $h$. Therefore we can form the localization $\mathcal{S}_{S'}^{-1}\, \mathbb{Z}[X]$, and $\mathrm{ev}_{S'}$ extends uniquely to a ring homomorphism

$$\mathcal{S}_{S'}^{-1}\, \mathbb{Z}[X] \to \mathrm{Map}(S', \mathrm{Quot}\, R).$$

Let $\mathbb{Z}[\![X]\!]_{S'}$ denote the subset of all elements $f \in \mathbb{Z}[\![X]\!]$ for which there exist $g \in \mathbb{Z}[X]$ and $h \in \mathcal{S}_{S'}$ with $h\, f = g$. This is in fact a $\mathbb{Z}$-submodule of $\mathbb{Z}[\![X]\!]$: indeed, given two equations $h\, f = g$, $h'\, f' = g'$ we have $h\, h'(f + f') = h'\, g + h\, g'$, and $\mathcal{S}_{S'}$ is multiplicative. We obtain a $\mathbb{Z}$-linear map

$$\mathbb{Z}[\![X]\!]_{S'} \to \mathcal{S}_{S'}^{-1}\, \mathbb{Z}[X].$$

by sending such an element $f$ to the quotient $g/h$ (this is well-defined: if $h\, f = g$, $h'\, f = g'$ then $g\, h' = h\, h'\, f = g'\, h$), and by composition the $\mathbb{Z}$-linear map

$$\mathrm{ev}_{S'}: \mathbb{Z}[\![X]\!]_{S'} \to \mathrm{Map}(S', \mathrm{Quot}\, R).$$

**Definition 1.1.** Let $f \in \mathbb{Z}[\![X]\!]$, $S' \subset S$ a subset, $x \in S'$. If $f$ is contained in $\mathbb{Z}[\![X]\!]_{S'}$ then $f$ is called *evaluable on* $S'$ and the value $\mathrm{ev}_{S'}(f)(x) \in \mathrm{Quot}(R)$ is called the *value of* $f$ *in* $x$.

**Remark 1.2.** (i) The notation "$\sum_{\lambda \in X} n_\lambda\, \lambda$" common in group rings would be ambiguous here: if $\lambda, \mu \in X$ are such that $\lambda + \mu \in X$ then one has to distinguish between the elements $e(\lambda + \mu)$ and $e(\lambda) + e(\mu)$. However, if $S$ is a generating subset of a group $G$ and the group $X$ consists of restrictions to $S$ of characters $G \to R^\times$ then $\mathrm{ev}_S$ is injective, and we could write "$\sum_\lambda n_\lambda\, \lambda$" instead of "$\sum_\lambda n_\lambda\, e(\lambda)$".

(ii) On the other hand the above applies to $X = R^\times$ (take $S$ a singleton); in this case we have the ring homomorphism

$$\begin{array}{rccc} \mathrm{ev}: & \mathbb{Z}[R^\times] & \longrightarrow & R \\ & \sum_\lambda n_\lambda\, e(\lambda) & \longmapsto & \sum_\lambda n_\lambda\, \lambda \end{array}$$

which of course is far from being injective.

(iii) The $\mathbb{Z}$-linear map $\mathrm{ev}_{S'}$ is multiplicative in the following sense: if $f$, $g$ are multipliable and

if $f, g, fg \in \mathbb{Z}[\![X]\!]_{S'}$ then $\mathrm{ev}_{S'}(fg) = \mathrm{ev}_{S'}(f)\,\mathrm{ev}_{S'}(g)$.

(iv) The value of $f$ in $x$ in the situation of Definition 1.1 does not depend on the ambient subset $S'$: If $S'' \subset S'$ is a subset then $\mathbb{Z}[\![X]\!]_{S'} \subset \mathbb{Z}[\![X]\!]_{S''}$ and $\mathrm{ev}_{S''}(f)(x) = \mathrm{ev}_{S'}(f)(x)$ for all $f \in \mathbb{Z}[\![X]\!]_{S'}$ and $x \in S''$.

(v) Suppose that the set $\mathcal{S}_{S'}$ contains no zero divisors. Then the natural map $\mathbb{Z}[\![X]\!]_{S'} \to \mathcal{S}_{S'}^{-1}\,\mathbb{Z}[X]$ is injective. Example: $S$ is a topological space, $S' \subset S$ is dense, $R$ is a complete topological ring, $X$ consists of continuous mappings, and the evaluation homomorphism $\mathrm{ev}_S$ is injective (cf. (i) above). Then $\mathcal{S}_{S'}$ contains no zero divisors.

## 2. Formal characters. Characters on a subset

Let $C\,|\,K$ be an extension of fields, $V$ a $K$-vector space, $S$ a commutative group, and let $\rho : S \to \mathrm{GL}(V)$ be a fixed $K$-linear representation of $S$ on $V$. Let $\rho_{(C)}$ be the $C$-linear representation of $S$ on $V_{(C)} = C \otimes_K V$ obtained by extension of scalars. Let $X(S)$ be the group of characters $S \to C^{\times}$.

**Definition 2.1.**

(i) Let $\lambda \in X(S)$ be a character. For $x \in S$ let $\lambda(x).$ denote the homothety on $V_{(C)}$ defined by $\lambda(x)$. The subspace of $V_{(C)}$

$$V_{(C)}^{\lambda} := \bigcap_{x\in S} \bigcup_{k\in\mathbb{N}} \ker\,(\rho_{(C)}(x) - \lambda(x).)^k$$

is called the *generalized weight space* of weight $\lambda$ of $\rho_{(C)}$. Its dimension is called the *algebraic multiplicity* of $\lambda$ in $\rho_{(C)}$.

(ii) Assume that $V_{(C)}$ is the sum of the generalized weight spaces of $\rho_{(C)}$ and that the algebraic multiplicity of each $\lambda \in X(S)$ in $\rho_{(C)}$ is finite. Then we say that the representation $\rho$ is *finitely trigonalisable over C* or *possesses a formal character over C,* and the element

$$\mathrm{Ch}(\rho) := \sum_{\lambda\in X(S)} \dim_C\!\left(V_{(C)}^{\lambda}\right) e(\lambda) \in \mathbb{Z}[\![X(S)]\!]$$

is called the *formal character* of $\rho$.

(iii) In addition to (ii), assume that $S' \subset S$ is a subset such that $\mathrm{Ch}(\rho)$ is evaluable on $S'$. Then we say that $\rho$ *possesses a character on* $S'$, and the map

$$\mathrm{ev}_{S'}(\mathrm{Ch}(\rho)) : S' \to C$$

is called the *character on $S'$* of $\rho$.

**Remark 2.2.** (i) Suppose $\rho$ is finitely trigonalisable over $C$. Since $(\rho_{(C)})_{(C)} = \rho_{(C)}$ we have $V^{\lambda}_{(C)} = (V_{(C)})^{\lambda}_{(C)}$ for all $\lambda$, hence $\rho_{(C)}$ is finitely trigonalisable, with $\mathrm{Ch}(\rho) = \mathrm{Ch}(\rho_{(C)})$. On the other hand, the formal character $\mathrm{Ch}(\rho)$ does not depend on the particular choice of $C$. Indeed, if $C' \mid C$ is an extension then $X(S)$ embeds into $X'(S) := \mathrm{Hom}\big(S, C'^{\times}\big)$, $\mathbb{Z}[\![X(S)]\!]$ embeds into $\mathbb{Z}[\![X'(S)]\!]$ (via "extension by zero"), and by exactness of the scalar extension functor we have

$$\ker(\rho_{(C')}(s) - \lambda.)^k = \big(\ker(\rho_{(C)}(s) - \lambda.)^k\big)_{(C')}$$

for all $k \in \mathbb{N}$, $\lambda \in X(S)$, $s \in S$; hence $\big(V^{\lambda}_{(C)}\big)_{(C')} = V^{\lambda}_{(C')}$ for all $\lambda \in X(S)$, and hence the element $\mathrm{Ch}(\rho) \in \mathbb{Z}[\![X'(S)]\!]$ stays the same whether calculated in $\mathbb{Z}[\![X(S)]\!]$ or in $\mathbb{Z}[\![X'(S)]\!]$.

(ii) If $\dim_K(V) < \infty$ and each $\rho(x)$ ($x \in S$) is trigonalisable over $C$ in the classical sense then $\rho$ is finitely trigonalisable over $C$ and possesses a ($K$-valued) character on $S$ which coincides with the usual character of the finite-dimensional representation $\rho$ (cf. [5], §5, Prop. 19).

**Lemma 2.3 (direct sums).** Let $V = \bigoplus_{i\in I} V_i$ be a decomposition into a direct sum of $\rho$-stable subspaces. Then $V$ is finitely trigonalisable over $C$ if and only if each subrepresentation $\rho\,|_{V_i}$ is finitely trigonalisable over $C$ and the family of formal characters $\big(\mathrm{Ch}\big(\rho_{V_i}\big)\big)_{i\in I}$ is summable in $\mathbb{Z}[\![X(S)]\!]$; in that case

$$\mathrm{Ch}(\rho) = \textstyle\sum_{i\in I} \mathrm{Ch}\big(\rho\,|_{V_i}\big).$$

*Proof.* This follows from $V_{(C)} = \oplus_{i\in I}(V_i)_{(C)}$ and $V^{\lambda}_{(C)} = \bigoplus_{i\in I}(V_i)^{\lambda}_{(C)}$ for each $\lambda \in X(S)$, where $(V_i)^{\lambda}_{(C)}$ denotes the generalized weight space of weight $\lambda$ of the subrepresentation $\rho_{(C)}\,|_{(V_i)_{(C)}}$. □

**Lemma 2.4 (subrepresentations).** Let $W \subset V$ be a $\rho$-stable subspace.

(i) The generalized weight space of weight $\lambda$ of $(\rho\,|_W)_{(C)}$ is equal to $W_{(C)} \cap V^{\lambda}_{(C)}$.

(ii) $W_{(C)} \cap \sum_{\lambda\in X(S)} V^{\lambda}_{(C)} = \sum_{\lambda\in X(S)} \big(W_{(C)} \cap V^{\lambda}_{(C)}\big)$.

(iii) $\rho\,|_W$ is finitely trigonalisable over $C$ if and only if $W_{(C)} \subset \sum_{\lambda\in X(S)} V^{\lambda}_{(C)}$ and $\dim_C\big(W_{(C)} \cap V^{\lambda}_{(C)}\big) < \infty$ for all $\lambda \in X(S)$; if this is the case then

$$\mathrm{Ch}(\rho\,|_W) = \sum_{\lambda\in X(S)} \dim_C\big(W_{(C)} \cap V^{\lambda}_{(C)}\big)\, e(\lambda).$$

*Proof.* (i) follows from the obviously equality $(\rho\,|_W)_{(C)} = \rho_{(C)}\,|_{W_{(C)}}$.

(ii) It is clear that $\sum_{\lambda\in X(S)}\left(W_{(C)} \cap V_{(C)}^{\lambda}\right)$ is contained in $W_{(C)} \cap \sum_{\lambda\in X(S)} V_{(C)}^{\lambda}$ . Vice versa, assume $w \in W_{(C)} \cap \sum_{\lambda\in X(S)} V_{(C)}^{\lambda}$. There is a finite subset $\Lambda \subset X(S)$ such that $w = \sum_{\lambda\in\Lambda} v_\lambda$ with $v_\lambda \in V_{(C)}^{\lambda}$. Put $U = \sum_{\lambda\in\Lambda} V_{(C)}^{\lambda}$, and for $x \in S$ and $a \in C$ define the generalized eigen space $U^{x,a} = \bigcup_{k\in\mathbb{N}} (\ker \rho_{(C)}(x)\,|_U - a)^k \subset U$. Thus $V_{(C)}^{\lambda} = \bigcap_{x\in S} U^{x,\lambda(x)}$ for $\lambda \in \Lambda$. We claim that there even is a finite subset $T \subset S$ such that

$$V_{(C)}^{\lambda} = \bigcap_{x\in T} U^{x,\lambda(x)} \text{ for all } \lambda \in \Lambda.$$

Indeed, there exists a finite subset $T \subset S$ which, for any two distinct $\lambda, \lambda' \in \Lambda$, contains an element $x$ such that $\lambda(x) \neq \lambda'(x)$. Since for any $x \in T$ the direct sum decomposition $U = \sum_{\lambda\in\Lambda} V_{(C)}^{\lambda}$ refines the direct sum decomposition $U = \sum_a U^{x,a}$ we have, for any $\lambda \in \Lambda$, $\bigcap_{x\in T} U^{x,\lambda(x)} = V_{(C)}^{\lambda_1} \oplus \ldots \oplus V_{(C)}^{\lambda_n}$ with $\lambda_1, \ldots, \lambda_n \in \Lambda$. But $\lambda_i \neq \lambda$ would imply $U^{x,\lambda(x)} \cap U^{x,\lambda_i(x)} = 0$ for some $x \in T$, and thereby $\bigcap_{x\in T} U^{x,\lambda(x)} \cap V_{(C)}^{\lambda_i} = 0$. This proves $V_{(C)}^{\lambda} = \bigcap_{x\in T} U^{x,\lambda(x)}$.

Now $w$ is contained in a finite-dimensional subspace $U' \subset U$ which is stable under each $\rho_{(C)}(x)$ $(x \in T)$. (Namely, write $T = \{x_1, \ldots, x_n\}$, for each $\lambda \in \Lambda$ and each $x_i \in T$ choose $k_{\lambda,i} \in \mathbb{N}$ such that $(\rho_{(C)}(x_i) - \lambda(x_i))^{k_{\lambda,i}}(v_\lambda)$ is zero, and let $U'$ be the span of the elements $\rho_{(C)}(x_1)^{k_1} \ldots \rho_{(C)}(x_n)^{k_n} v_\lambda$ $(\lambda \in \Lambda,\ 0 \leq k_i < k_{\lambda,i})$.) By the finite-dimensional theory (e.g. [5], §5, Prop. 19) the space $U'$ as well as the stable subspace $U' \cap W_{(C)}$ are trigonalisable simultaneously for the $\rho_{(C)}(x)$ $(x \in T)$, so that the decomposition $w = \sum_\lambda v_\lambda$ necessarily already takes place inside $U' \cap W_{(C)}$.

(iii) follows immediately from (i) and (ii). □

**Lemma 2.5 (increasing filtrations).** Let $V_0 \subset V_1 \subset \ldots \subset V$ be an increasing sequence of $\rho$-stable subspaces, and put $V_\infty := \bigcup_{i\in\mathbb{N}} V_i$. Then $\rho\,|_{V_\infty}$ is finitely trigonalisable over $C$ if and only if $\rho\,|_{V_i}$ is finitely trigonalisable over $C$ for each $i \in \mathbb{N}$ and the sequence $\left(\mathrm{Ch}\left(\rho\,|_{V_i}\right)\right)_{i\in\mathbb{N}}$ converges in $\mathbb{Z}[\![X(S)]\!]$; in that case

$$\lim_{i\in\mathbb{N}} \mathrm{Ch}\left(\rho\,|_{V_i}\right) = \mathrm{Ch}(\rho\,|_{V_\infty}).$$

*Proof.* Suppose that each $\rho\,|_{V_i}$ is finitely trigonalisable over $C$ and $\lim_{i\in\mathbb{N}} \mathrm{Ch}\left(\rho\,|_{V_i}\right)$ exists. According to Lemma 2.4, $(V_i)_{(C)} \subset \sum_{\lambda\in X(S)} V_{(C)}^{\lambda}$ for all $i$, hence $(V_\infty)_{(C)} \subset \sum_{\lambda\in X(S)} V_{(C)}^{\lambda}$, and for each $\lambda$, $(V_\infty)_{(C)} \cap V^\lambda_{(C)} = \bigcup_{i\in\mathbb{N}} \left((V_i)_{(C)} \cap V^\lambda_{(C)}\right)$. The existence of the above limit means that for each $\lambda$ there exists an $i_0$ such that $(V_i)_{(C)} \cap V^\lambda_{(C)} = \left(V_{i_0}\right)_{(C)} \cap V^\lambda_{(C)}$ for all $i \geq i_0$. Hence $(V_\infty)_{(C)} \cap V^\lambda_{(C)} = \left(V_{i_0}\right)_{(C)} \cap V^\lambda_{(C)}$, hence $\mathrm{Ch}_{V_\infty}(\rho) = \lim_{i\in\mathbb{N}} \mathrm{Ch}_{V_i}(\rho)$. - The reverse direction is clear. □

**Corollary 2.6 (block matrices).** (i) Suppose $\rho$ has a matrix representation of the form

$$s \mapsto \begin{pmatrix} A_0(s) & & * \\ & A_1(s) & \\ 0 & & \ddots \end{pmatrix}$$

with finite dimensional matrix representations $s \mapsto A_i(s)$ of $S$ over $K$ which are trigonalisable over $C$ ($i \in \mathbb{N}$). If the family $(\mathrm{Ch}(A_i))_{i\in\mathbb{N}}$ is summable in $\mathbb{Z}[\![X(S)]\!]$ then the representation $\rho$ is finitely trigonalisable over $C$ with formal character $\mathrm{Ch}\,\rho = \sum_{i\in\mathbb{N}} \mathrm{Ch}(A_i)$. □

**Lemma 2.7 (decreasing filtrations).** Let $V \supset W_0 \supset W_1 \supset \ldots$, $V \supset W'_0 \supset W'_1 \supset \ldots$ be two decreasing sequences of $\rho$-stable subspaces of finite codimension in $V$ such that for every $i \in \mathbb{N}$ there exists $j \in \mathbb{N}$ with $W'_j \subset W_i$ and $W_j \subset W'_i$. If each quotient representation $\rho_i$ of $\rho$ on $V/W_i$ is trigonalisable over $C$ and the sequence $(\mathrm{Ch}(\rho_i))_{i\in\mathbb{N}}$ converges in $\mathbb{Z}[\![X(S)]\!]$ then the same is true for the quotient representations $\rho'_i$ of $\rho$ on $V\big/W'_i$, and

$$\lim_{i\in\mathbb{N}} \mathrm{Ch}(\rho_i) = \lim_{i\in\mathbb{N}} \mathrm{Ch}\left(\rho'_i\right).$$

*Proof.* Let $i \in \mathbb{N}$. Choose $j \in \mathbb{N}$ such that $W_j \subset W'_i$. Then $V\big/W'_i$ is a quotient of the finite-dimensional vector space $V/W_j$. Hence with $\rho_j$ also $\rho'_i$ is trigonalisable over $C$, and the set of weights of $\rho'_i$ (with multiplicities) is a subset of the set of weights of $\rho_j$. This means that $\dim_C\left(\left(V\big/W'_i\right)^\lambda_{(C)}\right) \leq \dim_C\left((V/W_j)^\lambda_{(C)}\right)$ for all $\lambda$. Since $i$ was arbitrary and the sequence $(\mathrm{Ch}(V/W_j))_{j\in\mathbb{N}}$ converges it follows that the sequence $\left(\mathrm{Ch}\left(V\big/W'_i\right)\right)_{i\in\mathbb{N}}$ converges; namely, by symmetry, to the same limit. □

## 3. Convergent nets of characters on a subset

We keep the notations of the previous section. Suppose that $V$ is the union of an increasing sequence $V_0 \subset V_1 \subset \ldots$ of $\rho$-stable subspaces. We have seen (Lemma 2.5) that if each $V_i$ is trigonalisable over $C$ and the sequence of formal characters $\mathrm{Ch}(V_i)$ converges in $\mathbb{Z}[\![X(S)]\!]$ then $V$ is finitely trigonalisable over $C$, its formal character being the limit of the $\mathrm{Ch}(V_i)$. However, if the sequence $(\mathrm{Ch}(V_i))_i$ does not converge but each representation $V_i$ possesses a character on a subset $S' \subset S$ we may still ask if the sequence of characters on $S'$ converges in some sense.

**Convention.** By *weak convergence* of a net of $C$-valued functions we shall always mean pointwise convergence with respect to the discrete topology on $C$.

**Example 3.1.** Let $p$ be a prime number, and let $v_p$ denote the valuation of $\mathbb{Q}_p$ normalized by $v_p(p) = 1$. Let $S = \mathbb{Z}_p^\times$. Let $V = C^\infty(\mathbb{Z}_p, K)$ be the $K$-vector space of locally constant $K$- valued functions on $\mathbb{Z}_p$, endowed with the $S$-operation $\rho(s)(f) := (z \mapsto f(sz))$ $(s \in S,\ f \in V,\ z \in \mathbb{Z}_p)$. For any $h \in \mathbb{N}$ consider the disjoint open covering $U_{h,i} := i + p^h\,\mathbb{Z}_p$ $(0 \le i < p^h)$ of $\mathbb{Z}_p$. Then $V$ is the union of the $\rho$-stable subspaces

$$V_h := \left\{ f : \mathbb{Z}_p \to K;\ f \text{ constant on } U_{h,i} \text{ for all } 0 \le i < p^h \right\} \quad (h \in \mathbb{N}).$$

Since each $S$-representation $V_h$ is finite-dimensional it possesses a character $\theta_h = \mathrm{ev}_S(\mathrm{Ch}\,V_h)$ on $S$. The value of $\theta_h$ in $s \in S$ is the trace of the operator $s$ on $V_h$. The characteristic functions $\mathbb{1}_{U_{h,i}}$ $(0 \le i < p^h)$ form a basis of $V_h$. Letting $\left(\mathbb{1}^*_{U_{h,i}}\right)$ denote the dual basis we obtain

$$\begin{aligned}
\theta_h(s) &= \textstyle\sum_{0 \le i < p^h} \mathbb{1}^*_{U_{h,i}} \left(s\,\mathbb{1}_{U_{h,i}}\right) \\
&= \#\left\{0 \le i < p^h;\ s^{-1} i + p^h\,\mathbb{Z}_p = i + p^h\,\mathbb{Z}_p\right\} \\
&= \#\left\{0 \le i < p^h;\ v_p(i) \ge h - v_p\left(s^{-1} - 1\right)\right\} \\
&= p^{\min\left(h, v_p\left(s^{-1} - 1\right)\right)}
\end{aligned}$$

according to Lemma 3.2 below. Hence $\theta_h(1) = p^h$, and if $s \ne 1$ and $h \ge v_p\left(s^{-1} - 1\right)$ then $\theta_h(s) = p^{v_p\left(s^{-1}-1\right)} = \left|s^{-1} - 1\right|^{-1}$. It follows that the sequence $(\mathrm{ev}_{S-\{1\}}(\mathrm{Ch}\,V_h))_{h \in \mathbb{N}}$ converges weakly to the function $\theta : S - \{1\} \to K,\ s \mapsto \left|s^{-1} - 1\right|^{-1}$.

**Lemma 3.2.** For any two integers $0 \le \alpha \le \beta$ there are $p^\alpha$ elements $i$ contained in

$\{0, 1, \ldots, p^\beta - 1\}$ satisfying $v_p(i) \geq \beta - \alpha$.

*Proof*. Obvious. □

The next example shows that the limit function $\theta$ obtained in Example 3.1 is not an invariant of $\rho$ but depends on the particular system $(V_h)$ exhausting $V$.

**Example 3.3.** Define $S$, $V$, $\rho$ as in the previous example, but this time consider the subsets

$$\tilde{U}_{h,i} = \begin{cases} \{z^{-1};\, z \in i\,p^{1-h} + p^h\,\mathbb{Z}_p\} & ;\, i \in \mathbb{N} - \{0\} \\ p^h\,\mathbb{Z}_p & ;\, i = 0 \end{cases}$$

of $\mathbb{Q}_p$ $(h \in \mathbb{N})$. For $z \in \mathbb{Q}$ we have $z^{-1} \in \mathbb{Z}_p - p^h\,\mathbb{Z}_p \Leftrightarrow z \in p^{1-h}\,\mathbb{Z}_p - p\,\mathbb{Z}_p$, and

$$p^{1-h}\,\mathbb{Z}_p - p\,\mathbb{Z}_p = \dot{\bigcup_{i\in I}}\, i\,p^{1-h} + p^h\,\mathbb{Z}$$

where $I := \{1 \leq i < p^{2h-1},\, v_p(i) < h\}$. Hence the family $(\tilde{U}_{h,i})_{i\in I\cup\{0\}}$ is a disjoint open covering of $\mathbb{Z}_p$. Again $V$ is the union of the $\rho$-stable subspaces

$$\tilde{V}_h := \{f : \mathbb{Z}_p \to K;\ f \text{ constant on } \tilde{U}_{h,i} \text{ for all } i \in I \cup \{0\}\} \quad (h \in \mathbb{N}).$$

We calculate

$$\begin{aligned} \tilde{\theta}_h(s) := \mathrm{ev}_S(\mathrm{Ch}\,\tilde{V}_h) &= \textstyle\sum_{i\in I\cup\{0\}} \mathbb{1}^*_{\tilde{U}_{h,i}}\left(s\,\mathbb{1}_{\tilde{U}_{h,i}}\right) \\ &= 1 + \#\left\{i \in I;\, s\,i\,p^{1-h} + p^h\,\mathbb{Z}_p = i\,p^{1-h} + p^h\,\mathbb{Z}_p\right\} \\ &= 1 + \#\,\{i \in I;\, v_p(i) + v_p\,(s-1) \geq 2\,h - 1\}; \end{aligned}$$

here we used $s\,\mathbb{1}_{\tilde{U}_{h,0}} = \mathbb{1}_{\tilde{U}_{h,0}}$ for all $s$. For any fixed $s \neq 1$ and $h >> 0$ we have $v_p(i) + v_p(s-1) < 2\,h - 1$ for all $i \in I$, i.e. $\tilde{\theta}_h(s) = 1$. This means that the sequence $\left(\mathrm{ev}_{S-\{1\}}(\mathrm{Ch}\,\tilde{V}_h)\right)_{h\in\mathbb{N}}$ converges weakly to the function $S - \{1\} \to K$, $s \mapsto 1$.

## II. Distribution algebras of $p$-adic groups

This second part provides certain results concerning completed distribution algebras. The most important three are the following: the question when there are natural maps between completions with respect to different norms of the same distribution algebra (Section 5- 7), explicit bases of completed distribution algebras (viewed as $K$-Banach spaces; Section 7), and the description of the $\mathfrak{g}$-adic filtration by means of these bases (Section 8).

**Notations.** Let $p$ be a prime number, and let $C \mid K \mid L \mid \mathbb{Q}_p$ be a tower of complete valued fields such that $L \mid \mathbb{Q}_p$ is finite of degree $n$, $K$ is discretely valued, and $C$ is algebraically closed. Let $o_L$ be the ring of integers of $L$, $\pi_L \in o_L$ a prime element, and $|\ |_L$ the absolute value on $L$ normalized by $|\pi_L|_L = (\# o_L/\pi_L)^{-1}$. For the basic ideas of the theory of locally analytic representations of $p$-adic groups we refer to the series of papers [29], [30], [31], [32], [33]. By an $L$-analytic group we mean a finite-dimensional Lie group over $L$ in the sense of [7]. For such an $L$-analytic group $G$ with Lie algebra $\mathfrak{g}$ we denote by

- $C^{\mathrm{an}}(G, K)$ the locally convex $K$-vector space of locally analytic functions $G \to K$,
- $D(G, K) = C^{\mathrm{an}}(G, K)'_b$ (strong dual) the $K$-algebra of distributions on $G$, endowed with the convolution product,
- $U(\mathfrak{g})$ the enveloping $L$-algebra of $\mathfrak{g}$, viewed as a subalgebra of $D(G, K)$,
- $U(\mathfrak{g}, K)$ the closure of $U(\mathfrak{g}) \otimes_L K$ in $D(G, K)$

(cf. [29], Section 2). Finally, we let $R^L_{\mathbb{Q}_p}(G)$ denote the $\mathbb{Q}_p$-analytic group obtained from $G$ by scalar restriction ([9], 5.14); the underlying topological groups of $G$ and of $R^L_{\mathbb{Q}_p}(G)$ coincide.

### 4. Exp and Log

Let $G$ be an $L$-analytic group with $L$-Lie algebra $\mathfrak{g}$. On a suitable small additive open subgroup $V \subset \mathfrak{g}$ there exists an *exponential mapping* $\phi : V \to G$; this is a locally $L$-analytic map satisfying $\phi(m\,\mathfrak{x}) = \phi(\mathfrak{x})^m$ ($\mathfrak{x} \in V$, $m \in \mathbb{Z}$) and having the identity map $\mathfrak{g} \to \mathfrak{g}$ as the tangent map in 0 (cf. [7], III.4.3 and III.7.2). Exponential mappings are not unique. However, using the *logarithmic mapping* of $G$ one can make a natural choice as follows: Put

$$G_{\mathrm{f}} := \left\{ g \in G;\ \exists\, m \in \mathbb{N} : \lim_{n\to\infty} g^{m\,p^n} = 1 \right\}.$$

Then $G_{\mathrm{f}}$ is an open set of $G$, equal to the union of all compact subgroups of $G$, and closed under taking integral powers ([7], III.7.6, Lemma 1, Proposition 10 (i), and Corollary of Proposition 13). There is a *unique* map $\mathrm{Log}_G : G_{\mathrm{f}} \longrightarrow \mathfrak{g}$ with the following two properties:

- $\mathrm{Log}_G(g^m) = m\,\mathrm{Log}_G(g)$ for all $g \in G_{\mathrm{f}}$, $m \in \mathbb{Z}$,
- there are an open neighbourhood $U$ of 1 in $G$, an open neighbourhood $V$ of 0 in $\mathfrak{g}$, and a bijective exponential mapping $\phi : V \to U$ of $G$ such that $\mathrm{Log}_G|_U = \phi^{-1}$;

and this map is locally $L$-analytic ([7], III.7.6, Proposition 10).

**Proposition 4.1.** Let $G$ be an $L$-analytic group and $H \subset G$ an $L$-analytic subgroup. Then $H_{\mathrm{f}} = G_{\mathrm{f}} \bigcap H$ and $\mathrm{Log}_G|_{H_{\mathrm{f}}} = \mathrm{Log}_H$.

*Proof.* The first assertion is obvious. Let $V$ be an open neighbourhood of 0 in the Lie algebra $\mathfrak{g}$ of $G$ such that $\mathrm{Log}_G^{-1}|_V$ is an exponential mapping of $G$. According to [7], III.4.4, Proposition 8, there is an open neighbourhood $V'$ of 0 in $\mathrm{Lie}(H)$ such that $\mathrm{Log}_G^{-1}|_V$ restricts to an exponential mapping $\mathrm{Log}_G^{-1}|_{V'}$ of $H$. Since also $\mathrm{Log}_G(g^m) = m\,\mathrm{Log}(g)$ for all $g \in H_{\mathrm{f}}$ and $m \in \mathbb{Z}$, the restriction $\mathrm{Log}_G^{-1}|_{H_{\mathrm{f}}}$ is a logarithmic mapping of $H$, necessarily equal to $\mathrm{Log}_H$. □

Thus we shall often simply write $\mathrm{Log}_H = \mathrm{Log}$ when dealing with various $L$-analytic subgroups $H$ of a fixed $L$-analytic group $G$. We shall reserve the notation $\mathrm{Exp}_H = \mathrm{Exp}$ for exponential mappings which are inverse maps of logarithmic mappings.

Another notational remark: we capitalize the symbols Exp, Log in order to distinguish them from the exponential and logarithmic series in (completed) distribution algebras; these will always be denoted exp, log respectively. In a small neighbourhood of 1, 0 respectively both concepts coincide; for a more precise statement cf. Remark 7.3.

**Example 4.2.** Let $\mathsf{G}$ be a *unipotent* $L$-group, i.e. an algebraic $L$-group isomorphic to a closed subgroup of the upper strictly triangular subgroup of some $\mathrm{GL}_m/L$. Let $\mathsf{G}(L)$ be the corresponding group of $L$-rational points and $\mathfrak{g}$ its Lie algebra. Then $\mathsf{G}(L)_{\mathrm{f}} = \mathsf{G}(L)$, $\mathfrak{g}$ is nilpotent ([11], IV.2.2, Corollaire 2.13), and we have the global bijections

$$\mathsf{G}(L) \underset{\mathrm{Exp}}{\overset{\mathrm{Log}}{\rightleftarrows}} \mathfrak{g}$$

([11], IV.2.4 Proposition 4.1) which are given by the restrictions of the matrix logarithm resp. the matrix exponential of $\mathrm{GL}_m(L)$.

If $\mathfrak{x} \in \mathfrak{g}$ is an element in the domain of Exp and $H \subset G$ is an $L$-analytic subgroup with Lie algebra $\mathfrak{h} \subset \mathfrak{g}$ then the fact "$\mathfrak{x} \in \mathfrak{h}$" does not generally imply "$\mathrm{Exp}(\mathfrak{x}) \in H$": this is clear if one considers proper open subgroups $H \subset G$. However, in the case of algebraic groups we have the following:

**Proposition 4.3.** Let $\mathsf{H} \subset \mathsf{G}$ be an inclusion of affine algebraic $L$-groups. Let $G = \mathsf{G}(L)$, $H = \mathsf{H}(L)$ be the respective groups of $L$-rational points, viewed as $L$-analytic groups, with respective Lie algebras $\mathfrak{g}$ and $\mathfrak{h}$. There is an open neighbourhood $\mathfrak{g}^{\sim} \subset \mathfrak{g}$, depending only on $\mathsf{G}$, such that $\mathrm{Exp}(\mathfrak{g}^{\sim} \bigcap \mathfrak{h}) \subset H$.

*Proof.* There exist a number $m \in \mathbb{N}$ and an inclusion of affine algebraic $L$-groups $\mathsf{G} \to \mathrm{GL}_m/_L$ ([11], Corollaire II.5.5.2); hence by Proposition 4.1 we may assume $\mathsf{G} = \mathrm{GL}_m/_L$. Then $\mathfrak{g}$ is equal to the algebra of $m \times m$-matrices over $L$. Let $\mathrm{EXP} : \mathfrak{g} \to \mathsf{G}(L[\![T]\!])$ denote the formal exponential mapping defined in Section II.6.3.1 of [10]. In our case it is given by

$$\mathrm{EXP}(\mathfrak{x}) = \sum_{k \geq 0} \frac{1}{k!} \mathfrak{x}^k T^k = \left( \sum_{k \geq 0} \frac{1}{k!} x_{ij}^{(k)} T^k \right)_{ij}$$

where $\mathfrak{x}^k = \left( x_{ij}^{(k)} \right)_{ij}$ (cf. [11], Exemple II.6.3.3). Let $L\langle T \rangle \subset L[\![T]\!]$ be the subalgebra of power series converging on the closed unit ball. There is an open neighbourhood $\mathfrak{g}^{\sim} \subset \mathfrak{g}$ on which Exp is defined and given by the convergent series $\sum_{k \geq 0} \frac{1}{k!} \mathfrak{x}^k$. This means that for $\mathfrak{x} \in \mathfrak{g}^{\sim}$ the coefficients of the matrix $\mathrm{EXP}(\mathfrak{x})$ converge in $L$ if we substitute $T = 1$, hence they are contained in $L\langle T \rangle$. The map $\mathrm{Exp} : \mathfrak{g}^{\sim} \to \mathsf{G}(L)$ therefore factors as follows:

$$\mathfrak{g}^{\sim} \xrightarrow{\mathrm{EXP}} \mathsf{G}(L\langle T \rangle) \xrightarrow{T \mapsto 1} \mathsf{G}(L)$$

According to [11], Corollaire II.6.3.4 (c), the restriction of $\mathrm{EXP} : \mathfrak{g} \to \mathsf{G}(L[\![T]\!])$ to $\mathfrak{h}$ is equal to the intrinsically defined map $\mathrm{EXP} : \mathfrak{h} \to \mathsf{H}(L[\![T]\!])$; hence by restriction we obtain the map

$$\mathfrak{g}^{\sim} \bigcap \mathfrak{h} \xrightarrow{\mathrm{EXP}} \mathsf{G}(L\langle T \rangle) \bigcap \mathsf{H}(L[\![T]\!]) = \mathsf{H}(L\langle T \rangle) \xrightarrow{T \mapsto 1} \mathsf{H}(L)$$

which proves our assertion. □

## 5. Uniform pro- *p* groups

In this section we clarify how two arbitrary uniform pro-$p$-groups in the sense of [12] can intersect. Instead of recalling the formal definition of such groups (cf. [12], Section 4.1) we will collect their basic properties.

Let $H$ be a uniform pro-$p$-group of dimension $d \in \mathbb{N}$. Any $d$-tupel $\mathbf{h} = (h_1, \ldots, h_d)$ of topological generators of $H$ will be called an *ordered basis* of $H$ and gives rise to a homeomorphism

$$\begin{array}{rccc} \psi_{\mathbf{h}}: & \mathbb{Z}_p^d & \longrightarrow & H \\ & \alpha & \longmapsto & h_1^{\alpha_1} \ldots h_d^{\alpha_d} \end{array}$$

and thereby to an isomorphism of locally convex $K$-vector spaces

$$\psi_{\mathbf{h}}^* : C^{\mathrm{an}}(H, K) \to C^{\mathrm{an}}\left(\mathbb{Z}_p^d, K\right).$$

The uniform pro-$p$-group $H$ is a $\mathbb{Q}_p$-analytic group with global chart $\psi_{\mathbf{h}}^{-1}$. Conversely, every $\mathbb{Q}_p$-analytic group $G$ contains a uniform pro-$p$-group $H$ as an open subgroup. A basis of neighbourhoods of 1 in $H$ (as well as in $G$) is given by the *lower p-series* of $H$, i.e. the uniform pro-$p$-groups $H^{p^m} = \left\{h^{p^m}; h \in H\right\}$ $(m \in \mathbb{N})$.

In the book [12] a map $\log : H \to \Lambda$ into a certain $\mathbb{Z}_p$-Lie algebra $\Lambda$ (constructed inside a completion of the group ring of $H$ over $\mathbb{Q}_p$) is defined and shown to be an isomorphism when the underlying set of $H$ is endowed with a $\mathbb{Z}_p$- Lie algebra structure as follows:

$$z \cdot g := \lim_{n\to\infty} g^{z_n}$$

$$g + h := \lim_{n\to\infty}\left(g^{p^n} h^{p^n}\right)^{p^{-n}},$$

$$[g, h] := \lim_{n\to\infty}\left(g^{-p^n} h^{-p^n} g^{p^n} h^{p^n}\right)^{p^{-2n}}$$

($g, h \in H$, $z = \lim_{n\to\infty} z_n \in \mathbb{Z}_p$ with $z_n \in \mathbb{N}$; cf. [12], Definitions 1.25, 4.12, 4.29 and Corollary 7.14).

**Lemma 5.1.** Let $H$ be a uniform pro-$p$-group with Lie algebra $\mathfrak{h}$. Then $H_{\mathrm{f}} = H$, and the map $\mathrm{Log} : (H, +, [\,,\,]) \to \mathfrak{h}$ is an embedding of $\mathbb{Z}_p$-Lie algebras.

*Proof.* Since $H$ is compact we have $H = H_{\mathrm{f}}$.

Log is $\mathbb{Z}_p$-linear: Let $g \in H$ and $z = \lim_{n\to\infty} z_n \in \mathbb{Z}_p$ with $z_n \in \mathbb{N}$. Since addition in $(H, +, [\ ,\ ])$ satisfies $g^{z_n} = g + \ldots + g$ ($z_n$ times) and by continuity we have $\mathrm{Log}(z \cdot g) = \lim_{n\to\infty} \mathrm{Log}(g^{z_n}) = \lim_{n\to\infty} z_n \,\mathrm{Log}(g) = z\,\mathrm{Log}(g)$.

Log is injective: Proposition 12 of [7], III.7.6, assures that the kernel of Log consists of torsion elements, hence is trivial by Theorem 4.5 of [12].

Log is additive and respects the bracket: Let $g, h \in H$, $\mathfrak{x} = \mathrm{Log}(g)$, $\mathfrak{y} = \mathrm{Log}(h)$. Then

$$\begin{aligned}\mathrm{Log}\,(g+h) &= \lim_{n\to\infty} \mathrm{Log}\left(\left(g^{p^n} h^{p^n}\right)^{p^{-n}}\right)\\ &= \lim_{n\to\infty} p^{-n}\,\mathrm{Log}\,(\mathrm{Exp}\,(p^n\,\mathfrak{x})\,\mathrm{Exp}\,(p^n\,\mathfrak{y}))\\ &= \mathfrak{x} + \mathfrak{y}\,;\end{aligned}$$

the last equality is an instance of [7], III.7.2, Proposition 4 (1). Finally, using [7], III.7.2, Proposition 4 (2), we calculate

$$\begin{aligned}\mathrm{Log}\,([g, h]) &= \lim_{n\to\infty} \mathrm{Log}\left(\left[g^{p^n}, h^{p^n}\right]^{p^{-2n}}\right)\\ &= \lim_{n\to\infty} p^{-2n}\,\mathrm{Log}(\mathrm{Exp}(-p^n\,\mathfrak{x})\,\mathrm{Exp}(-p^n\,\mathfrak{y})\,\mathrm{Exp}(p^n\,\mathfrak{x})\,\mathrm{Exp}(p^n\,\mathfrak{y}))\\ &= [\mathfrak{x}, \mathfrak{y}]\,.\ \square\end{aligned}$$

Thus if $H$ is a uniform pro-$p$-group which at the same time is an $L$-analytic subgroup of some $L$-analytic group $G$ we will identify the $\mathbb{Z}_p$-Lie algebra $\Lambda$ mentioned above with the $\mathbb{Z}_p$-Lie subalgebra $\mathrm{Log}(H)$ of the $L$-Lie algebra of $G$, thereby identifying the map $\log : H \to \Lambda$ defined in [12] with our logarithmic mapping $\mathrm{Log} : H \to \mathrm{Log}(H)$.

We know that $\mathrm{Log}(H)$ is free as a $\mathbb{Z}_p$-module and that the ordered bases of $H$ correspond under Log to the $\mathbb{Z}_p$-bases of $\mathrm{Log}(H)$ ([12], Theorem 4.17 and Ex. 3 (i) of Section 8). Although we will not make any use of it we mention that the inverse map of $\mathrm{Log} : H \to \mathrm{Log}(H)$ is an exponential mapping (cf. [12], Ex. 3 (iv) of Section 8), and we will always denote it by Exp.

**Proposition 5.2 (elementary divisor theorem for uniform pro-$p$-groups).** Let $H$ be a uniform pro-$p$-group of dimension $d$. Let $H' \subset H$ be a closed subgroup which is a uniform pro-$p$-group. There are an ordered basis $(h_1, \ldots, h_d)$ of $H$ and integers $0 \le \alpha(1) \le \ldots \le \alpha(d')$ ($d' \le d$) such that $\left(h_1^{p^{\alpha(1)}}, \ldots, h_{d'}^{p^{\alpha(d')}}\right)$ is an ordered basis of $H'$. The number $d'$ and the tuple $(\alpha(1), \ldots, \alpha(d'))$ are uniquely determined by $H$ and $H'$. Moreover, $H'$ is open in $H$ if and only if

$d' = d$.

*Proof.* Let $\Lambda = \mathrm{Log}(H)$, $\Lambda' = \mathrm{Log}(H')$ denote the respective $\mathbb{Z}_p$-Lie algebras of $H$, $H'$. Then $\Lambda'$ is a submodule of the free rank-$d$-$\mathbb{Z}_p$-module $\Lambda$ (cf. [12], Theorem 4.17, Proposition 4.31 and Corollary 7.14). According to the elementary divisor theorem for principal ideal domains (cf. [5], VII, §4, Proposition 9) there are a $\mathbb{Z}_p$-basis $(\mathfrak{x}_1, \ldots, \mathfrak{x}_d)$ of $\Lambda$ and elements $c_1, \ldots, c_{d'} \in \mathbb{Z}_p$ $(d' \leq d)$, uniquely determined by $\Lambda$ and $\Lambda'$ up to a unit in $\mathbb{Z}_p$, such that $(c_1 \mathfrak{x}_1, \ldots, c_{d'} \mathfrak{x}_{d'})$ is a $\mathbb{Z}_p$-basis of $\Lambda'$ and $c_1 \mid c_2 \mid \ldots \mid c_{d'}$. After scaling by a unit in $\mathbb{Z}_p$ we may in fact attain $c_i = p^{\alpha(i)}$ with natural numbers $\alpha(1) \leq \ldots \leq \alpha(d')$. But then $(\mathrm{Exp}\,\mathfrak{x}_1, \ldots, \mathrm{Exp}\,\mathfrak{x}_d)$ is an ordered basis of $H$ and $(\mathrm{Exp}(c_1 \mathfrak{x}_1), \ldots, \mathrm{Exp}(c_{d'} \mathfrak{x}_{d'})) = ((\mathrm{Exp}\,\mathfrak{x}_1)^{p^{\alpha(1)}}, \ldots, (\mathrm{Exp}\,\mathfrak{x}_{d'})^{p^{\alpha(d')}})$ is an ordered basis of $H'$.

Let $(\tilde{h}_1, \ldots, \tilde{h}_d)$ be another ordered basis of $H$ such that $\tilde{h}^{p^{\beta(1)}}, \ldots, \tilde{h}^{p^{\beta(d'')}}$ is an ordered basis of $H'$, $0 \leq \beta(1) \leq \ldots \leq \beta(d'')$, $d'' \leq d$. Again by [12], Theorem 4.17 and Corollary 7.14, we have that $(\mathrm{Log}(\tilde{h}), \ldots, \mathrm{Log}(\tilde{h}))$ is a $\mathbb{Z}_p$-basis of $\Lambda$ and $(\mathrm{Log}(\tilde{h}^{p^{\beta(1)}}), \ldots, \mathrm{Log}(\tilde{h}^{p^{\beta(d'')}})) = (p^{\beta(1)} \mathrm{Log}(\tilde{h}), \ldots, p^{\beta(d'')} \mathrm{Log}(\tilde{h}))$ is a $\mathbb{Z}_p$-basis of the free submodule $\Lambda'$. By the uniqueness of the above elementary divisors $c_1, \ldots, c_d$ we deduce $d' = d''$, $\beta(1) = \alpha(1), \ldots, \beta(d') = \alpha(d')$.

If $d'$ equals $d$ then $H'$ contains a member of the lower $p$-series of $H$ and hence is open in $H$. The converse follows from [12], Proposition 4.4. □

**Remark 5.3.** The above proposition has no converse, in the following sense: Let $H$ be a uniform pro-$p$-group with ordered basis $(h_1, \ldots, h_d)$, and let $\alpha(1), \ldots, \alpha(d') \geq 0$. Then the closed subgroup $H'$ of $H$ generated by $\left(h_1^{p^{\alpha(1)}}, \ldots, h_{d'}^{p^{\alpha(d')}}\right)$ in general is not uniform. Example: the pro-$p$-groups $\begin{pmatrix} 1 & p\mathbb{Z}_p & \mathbb{Z}_p \\ & 1 & p\mathbb{Z}_p \\ & & 1 \end{pmatrix}$, $\begin{pmatrix} 1 & p\mathbb{Z}_p & p\mathbb{Z}_p \\ & 1 & p\mathbb{Z}_p \\ & & 1 \end{pmatrix}$ are uniform while $\begin{pmatrix} 1 & p\mathbb{Z}_p & p^2\mathbb{Z}_p \\ & 1 & p\mathbb{Z}_p \\ & & 1 \end{pmatrix}$ is not.

**Corollary 5.4.** Let $H$ be a uniform pro-$p$-group, $g \in H$. Then there is an ordered basis $(h_1, \ldots, h_d)$ of $H$ and a number $\alpha \in \mathbb{N}$ such that $g = h_1^{p^{\alpha}}$.

*Proof.* The group $H$ has no torsion ([12], Theorem 4.5), hence the closed subgroup generated by $g$ is a uniform pro-$p$-group (isomorphic to $\mathbb{Z}_p$, cf. [12], Proposition 1.26 (iii)). Now apply Proposition 5.2. □

**Definition 5.5.** (i) Given a pair of uniform pro-$p$ groups $H, H'$ such that $H' \subset H$ is a closed subgroup we call the sequence $\alpha = (\alpha(1), \ldots, \alpha(d'))$ determined by the above proposition the *sequence of p-elementary divisors* of $H'$ in $H$, and we call $\alpha(d')$ the *highest p-elementary divisor* of $H'$ in $H$.

(ii) Let $H, H'$ be uniform pro-$p$-groups which are open subgroups of some topological group $G$. Choose $k \in \mathbb{N}$ such that $H^{p^k} \subset H'$, and let $\alpha$ be the sequence of $p$-elementary divisors of $H^{p^k}$ in $H'$. Then the number

$$\mathrm{ud}(H, H') := \alpha(d) - \alpha(1)$$

is called the *uniform defect* of $H$ and $H'$.

Part (ii) of the definition is justified by the subsequent lemma; note that the uniform defect of the two groups is zero if and only if one is a member of the lower $p$-series of the other. Note also that, in the situation of part (i) of the definition, the highest $p$-elementary divisor of $H'$ in $H$ is zero if and only if $H'$ is *compatible* with $H$ in the sense of [20], Section 1.3.

**Lemma 5.6.** Let $G$ be a topological group. Let $H, H'$ be open subgroups of $G$ which are uniform pro-$p$-groups. Let $k, k', l, l' \in \mathbb{N}$ such that $H^{p^k} \subset H'^{p^{k'}}$, $H'^{p^{l'}} \subset H^{p^l}$. Let $\alpha$ (resp. $\beta$) be the sequence of $p$-elementary divisors of $H^{p^k}$ in $H'^{p^{k'}}$ (resp. of $H'^{p^{l'}}$ in $H^{p^l}$). Then the differences $\alpha(d) - \alpha(1)$ and $\beta(d) - \beta(1)$ are equal and do not depend on the choices of $k, k', l, l'$.

*Proof.* If we replace $k, k'$ by $\tilde{k}, \tilde{k}'$ with $\tilde{k} \geq k$, say, then the numbers $\alpha(i)$ are replaced by $\tilde{\alpha}(i) = \alpha(i) + \tilde{k} - k - \tilde{k}' + k'$ (note that regarding Lemma 4.10 of [12] this makes sense even in the case $\tilde{\alpha}(i) < \alpha(i)$). Hence the quantity $\alpha(d) - \alpha(1)$ does not depend on the choice of $k, k'$. Similarly for $\beta(d) - \beta(1)$.

Let $(h_1, \ldots, h_d)$ be an ordered basis of $H'^{p^{k'}}$ such that $(h_1^{p^{\alpha(1)}}, \ldots, h_d^{p^{\alpha(d)}})$ is an ordered base of $H^{p^k}$. Fix $N \geq \alpha(d)$. Then $H'^{p^{k'+N}}$ has the ordered basis $(h_1^{p^N}, \ldots, h_d^{p^N})$ and is contained in $H^{p^k}$. We see that the unique $p$-elementary divisors of $H'^{p^{k'+N}}$ in $H^{p^k}$ are given by $N - \alpha(d) \leq \ldots \leq N - \alpha(1)$. Consequently, by the first part of the proof, $\beta(d) - \beta(1) = (N - \alpha(1)) - (N - \alpha(1)) = \alpha(d) - \alpha(1)$. □

**Lemma 5.7.** Let $H$ be a uniform pro-$p$-group and $H' \subset H$ an open normal subgroup which is a uniform pro-$p$-group. Let $\alpha$ be the sequence of $p$-elementary divisors of $H'$ in $H$ and choose an ordered basis $(h_1, \dots, h_d)$ of $H$ such that $\left(h_1^{p^{\alpha(1)}}, \dots, h_d^{p^{\alpha(d)}}\right)$ is an ordered basis of $H'$. Then the elements $h_1^{l_1}, \dots, h_d^{l_d}$ $(0 \le l_\nu < p^{\alpha(\nu)})$ constitute a system of representatives for $H/H'$.

*Proof.* A general element of $H$ may be written as $h_1^{k_1 p^{\alpha(1)}+l_1} \dots h_d^{k_d p^{\alpha(d)}+l_d}$ with $k_\nu \in \mathbb{Z}_p$, $0 \le l_\nu < p^{\alpha(\nu)}$; since $H'$ is normal in $H$ this element is contained in $h_1^{l_1} \dots h_d^{l_d} H'$. Hence the family $\left(h_1^{l_1}, \dots, h_d^{l_d}\right)_{0 \le l_\nu < p^{\alpha(\nu)}}$ *contains* a system of representatives. The same argument applies to the normal subgroup group $H^{p^{\alpha(d)}}$ of $H'$. But we know that the group $H \big/ H^{p^{\alpha(d)}}$ has cardinality $p^{d\alpha(d)}$ (cf. Section 4.2 of [12]), so the mentioned family already has to be a system of representatives. □

## 6. Distribution algebras of uniform pro- $p$ groups

We study the restriction of norms of distribution algebras of a uniform group to distribution algebras of uniform subgroups, thereby generalizing a result of T. Schmidt (cf. [26], Proposition 4.5 and Lemma 4.6).

Let $H$ be a uniform pro- $p$-group of dimension $d$, viewed as a $\mathbb{Q}_p$-analytic group. The distribution algebra $D(H, K)$ is a $K$- Fréchet algebra and contains the group ring $\mathbb{Z}_p[H]$. The image of a group element $h \in H$ in $\mathbb{Z}_p[H] \subset D(H, K)$ will be denoted by $\delta_h$ (the "Dirac distribution" of $h$).

**Convention.** In the following we will view $\mathbb{R}^d$ and $\mathbb{N}^d$ as partially ordered sets via the product ordering

$$(a_1, \dots, a_d) \le (b_1, \dots, b_d) \Leftrightarrow a_\nu \le b_\nu \text{ for all } 0 \le \nu \le d.$$

Then

$$I(d) := \left\{(r_1, \dots, r_d);\, r_1, \dots, r_d \in p^{\mathbb{Q}} \cap \left[p^{-1}, 1\right)\right\} \subset \mathbb{R}^d$$

is a directed set, and the subset $\{(r, \dots, r);\, r \in p^{\mathbb{Q}} \cap \left[p^{-1}, 1\right)\}$ of diagonal tuples is cofinal in $I(d)$. As a matter of notation we write

$$|\alpha| = \sum_{1 \le \nu \le d} \alpha_\nu,$$

$$r^\alpha = \prod_{1\le \nu\le d} r_\nu^{\alpha_\nu},$$
$$\boldsymbol{\lambda}^\alpha = \lambda_1^{\alpha_1} \dots \lambda_d^{\alpha_d}$$

for all $r \in \mathbb{R}^d$, $\boldsymbol{\lambda} = (\lambda_1, \dots, \lambda_d) \in D(H, K)^d$, $\alpha \in \mathbb{N}^d$.

Choose an ordered basis $\mathbf{h} = (h_1, \dots, h_d)$ of $H$, and put $\mathbf{b} = (b_1, \dots, b_d) := (\delta_{h_1} - 1, \dots, \delta_{h_d} - 1) \in \mathbb{Z}_p[H]^d$. Then for each $f \in C^{\mathrm{an}}(H, K)$ the Mahler expansion ([21], III, 1.2.4) of the function $\psi_{\mathbf{h}}^*(f)$ is given by

$$\psi_{\mathbf{h}}^*(f) = \sum_{\alpha\in\mathbb{N}^d} \mathbf{b}^\alpha(f) \binom{\cdot}{\alpha}.$$

From the characterization of locally analytic functions on $\mathbb{Z}_p^d$ by means of their Mahler coefficients ([21], III, 1.3.9) it follows that every distribution $\lambda \in D(H, K)$ has a unique convergent expansion

$$\lambda = \sum_{\alpha\in\mathbb{N}^d} d_\alpha\, \mathbf{b}^\alpha$$

with a family $(d_\alpha)_{\alpha\in\mathbb{N}^d}$ in $K$ satisfying $\lim_{|\alpha|\to\infty}|d_\alpha|\, r^\alpha = 0$ for any $r \in I(d)$. Conversely, any such family $(d_\alpha)$ defines an element $\lambda \in D(H, K)$ by the above formula. For every $r \in I(d)$ we define a norm $\|\ \|_{\mathbf{h},r}$ on $D(H, K)$ by

$$\|\lambda\|_{\mathbf{h},r} := \sup_{\alpha\in\mathbb{N}^d} |d_\alpha|\, r^\alpha;$$

the *dominant index* of $\lambda$, i.e. the maximal (cf. the above convention) index $\alpha \in \mathbb{N}^d$ satisfying $\|\lambda\|_{\mathbf{h},r} = |d_\alpha|\, r^\alpha$ will be denoted by

$$\alpha = \mathrm{dom}_{\mathbf{h},r}(\lambda).$$

Then the original Fréchet topology on $D(H, K)$ is defined by the family of norms $\|\ \|_{\mathbf{h},r}$ $(r \in I(d))$. Let $D_{\mathbf{h},r}(H, K)$ denote the completion of $D(H, K)$ with respect to the norm $\|\ \|_{\mathbf{h},r}$; this is a $K$-Banach space allowing the family $\left(\mathbf{b}^\alpha\right)_{\alpha\in\mathbb{N}^d}$ as an orthogonal basis.

Put $\kappa = 1$ if $p \neq 2$, $\kappa = 2$ if $p = 2$. In case $(r, \dots, r) \in I(d)$ is a *diagonal* tuple we know by the work of Schneider-Teitelbaum that the norm

$$\|\ \|_r := \|\ \|_{\mathbf{h},(r^\kappa,\dots,r^\kappa)}$$

is multiplicative and does not depend on the choice of the ordered basis $\mathbf{h}$, and that the completion

$$D_r(H,K) := D_{\mathbf{h},(r^\kappa,\ldots,r^\kappa)}(H,K)$$

is a noetherian $K$-Banach algebra (cf. [32], Theorem 4.5, Remark 4.6 and the remarks following Theorem 4.10; take into account that the map $h \mapsto \max(m+\kappa;\ h \in H^{p^m})$ is a $p$-valuation of $H$). Before stating the main result of this section we provide two lemmata.

**Lemma 6.1.** Let $H$ be a uniform pro- $p$-group with ordered basis $\mathbf{h} = (h_1, \ldots, h_d)$ and put $\mathbf{b} = (\delta_{h_1} - 1, \ldots, \delta_{h_d} - 1)$. Let $r \in I(d)$. Let $\left(\lambda_i = \sum_{\alpha\in\mathbb{N}^d} a_{i,\alpha}\, \mathbf{b}^\alpha\right)_{i\in I}$ be a family in $D_{\mathbf{h},r}(H,K)$ such that, for each $i \in I$, $S_i = \mathrm{dom}_{\mathbf{h},r}(\lambda_i)$ is the only element of $\mathbb{N}^d$ satisfying

$$\|\lambda_i\|_{\mathbf{h},r} = \left\|a_{i,S_i}\, \mathbf{b}^{S_i}\right\|_{\mathbf{h},r}.$$

If the mapping $i \mapsto \mathrm{dom}_{\mathbf{h},r}(\lambda_i) : I \to \mathbb{N}^d$ is injective (resp. bijective) then $(\lambda_i)_{i\in I}$ is an orthogonal family (resp. an orthogonal basis) in $D_{\mathbf{h},r}(H,K)$.

*Proof.* For diagonal tuples $r = (r_0, \ldots, r_0)$ this is [26], Corollary 4.2, based on Lemma 2 of Section 1.4 in [16]. One easily checks that the proofs of these results literally carry over to general tuples $r \in I(d)$. □

**Lemma 6.2.** (i) Let $a, b, c, d \in \mathbb{N}$. Then

$$\sum_{k\geq 0} (-1)^{a-k} \binom{a}{k}\binom{c\,k+d}{b} = \begin{cases} c^a & \text{if } b = a\,, \\ 0 & \text{if } b < a\,. \end{cases}$$

(ii) Let $a, b, h \in \mathbb{N}$. Then

$$\sum_{k\geq 0} (-1)^{a-k} \binom{a}{k}\binom{p^h k}{b} \begin{cases} = 1 & \text{if } b = a\,p^h, \\ = 0 & \text{if } b < a \text{ or } b > a\,p^h \\ \equiv 0 \pmod p & \text{if } b \neq a\,p^h. \end{cases}$$

*Proof.* (i) We use induction on $c$, the case $c = 0$ being well-known. The addition theorem gives $\binom{a}{k}\binom{c\,k+d}{b} = \sum_{l=0}^{b} \binom{a}{k}\binom{k}{l}\binom{(c-1)\,k+d}{b-l} = \sum_{l=0}^{b} \binom{a}{l}\binom{a-l}{k-l}\binom{(c-1)\,k+d}{b-l}$, hence

$$\sum_{k\geq 0} (-1)^{a-k} \binom{a}{k}\binom{c\,k+d}{b} = \sum_{l=0}^{b} \binom{a}{l} \sum_{k\geq 0} (-1)^{a-k+l} \binom{a-l}{k}\binom{(c-1)\,(k+l)+d}{b-l}.$$

If $b < a$ then $b - l < a - l$, and the assertion follows by induction. If $b = a$ then by induction $\sum_{k\geq 0} (-1)^{a-k+l} \binom{a-l}{k}\binom{(c-1)\,(k+l)+d}{b-l} = (c-1)^{a-l}$, and the assertion follows from the binomial formula.

(ii) (owed to E. Nagel) We have $\binom{a}{k}\binom{p^h k}{b} \neq 0$ only for $a \geq k$ and $p^h k \geq p^h b$, whence the cases $b = a\, p^h$ and $b > a\, p^h$. The case $b < a$ is a special case of part (i). - We claim the following:

$$\binom{p^h k}{b} \equiv \begin{cases} 0 & (\text{mod } p) \text{ if } b \neq p^h\, i \text{ for all } 0 \leq i \leq k, \\ \binom{k}{i} & (\text{mod } p) \text{ if } b = p^h\, i \text{ for some } 0 \leq i \leq k. \end{cases}$$

A summand in the sum $\binom{p^h k}{b} = \sum_{b_1 + \ldots + b_k = b} \binom{p^h}{b_1} \ldots \binom{p^h}{b_k}$ is congruent to 1 mod $p$ if and only if

$(*)$ each $b_\nu$ $(0 \leq \nu \leq k)$ is equal to either 0 or $p^h$

and is zero mod $p$ in all other cases. But the number of tuples $(b_1, \ldots, b_k)$ satisfying $b_1 + \ldots + b_k = b$ and $(*)$ is equal to

$$\begin{cases} 0 & \text{if } b \neq p^h\, i \text{ for all } 0 \leq i \leq a, \\ \binom{k}{i} & \text{if } b = p^h\, i \text{ for some } 0 \leq i \leq a \end{cases}$$

as in the latter case $i$ tokens of $p^h$ have to be distributed onto $k$ possible indices (without consideration of order). This proves our claim.

Since also $\binom{p^h k}{b} = \binom{k}{i}$ for $b > p^h k$ and $i > k$ we conclude

$$\sum_{k \geq 0} (-1)^{a-k} \binom{a}{k}\binom{p^h k}{b} \equiv \begin{cases} 0 & \text{if } b \neq p^h\, i \text{ for all } 0 \leq i \leq a, \\ \sum_{k \geq 0} (-1)^{a-k} \binom{a}{k}\binom{k}{i} & \text{if } b = p^h\, i \text{ for some } 0 \leq i \leq a, \end{cases}$$

and the lemma follows from the part (i). □

**Proposition 6.3.** Let $H$ be a uniform pro-$p$-group of dimension $d$. Let $H' \subset H$ be an open subgroup which is a uniform pro-$p$-group. Let $\gamma \in \mathbb{N}^d$ be the sequence of $p$-elementary divisors of $H'$ in $H$. Choose $r \in I(d)$ and an ordered basis $\mathbf{h} = (h_1, \ldots, h_d)$ of $H$ such that $r' := \left(r_1^{p^{\gamma(1)}}, \ldots, r_d^{p^{\gamma(d)}}\right)$ is contained in $I(d)$, the norm $\| \ \|_{\mathbf{h},r}$ is multiplicative, and $\mathbf{h}' := \left(h_1^{p^{\gamma(1)}}, \ldots, h_d^{p^{\gamma(d)}}\right)$ is an ordered basis of $H'$. Then the norm $\| \ \|_{\mathbf{h}',r'}$ on $D(H', K)$ is equal to the restriction of the norm $\| \ \|_{\mathbf{h},r}$ on $D(H, K)$ to $D(H', K)$ and in particular is multiplicative.

*Proof.* Write $\mathbf{b} = \left(\delta_{h_1} - 1, \ldots, \delta_{h_d} - 1\right)$, $\mathbf{b}' = \left(\delta_{h'_1} - 1, \ldots, \delta_{h'_d} - 1\right)$. We have to show that

$$\left\|\sum_\alpha d'_\alpha\, \mathbf{b}'^\alpha\right\|_{\mathbf{h},r} = \left\|\sum_\alpha d'_\alpha\, \mathbf{b}'^{\,\alpha}\right\|_{\mathbf{h}',r'}.$$

for every convergent series

$$(*)\qquad \textstyle\sum_{\alpha\in\mathbb{N}^d} d'_\alpha\, \mathbf{b}'^{\alpha} \quad (d'_\alpha \in K)$$

in $D(H', K)$.

We first determine, for each $\beta \in \mathbb{N}^d$, the expansion

$$\mathbf{b}'^{\beta} = \sum_{\alpha\in\mathbb{N}^d} t_{\beta,\alpha}\, \mathbf{b}^{\alpha} \quad (t_{\beta,\alpha} \in K).$$

We calculate $b_\nu'^{\beta_\nu} = \left((b_\nu+1)^{p^{\gamma(\nu)}} - 1\right)^{\beta_\nu} = \sum_{l\geq 0} (-1)^{\beta_\nu - l} \binom{\beta_\nu}{l} (b_\nu+1)^{l\,p^{\gamma(\nu)}} = \sum_{k\geq 0}\sum_{l\geq 0} (-1)^{\beta_\nu - l} \binom{\beta_\nu}{l}\binom{l\,p^{\gamma(\nu)}}{k} b_\nu^k$; therefore

$$t_{\beta,\alpha} = \prod_{1\leq\nu\leq d}\sum_{l\geq 0} (-1)^{\beta_\nu - l}\binom{\beta_\nu}{l}\binom{l\,p^{\gamma(\nu)}}{\alpha_\nu}.$$

In case $\beta = (0, \ldots, 1, \ldots, 0)$ is the $\nu$-th unit vector this reduces to

$$t_{(0,\ldots,1,\ldots,0),\alpha} = \begin{cases} \binom{p^{\gamma(\nu)}}{\alpha_\nu} & ; 1 \leq \alpha_\nu \leq p^{\gamma(\nu)} \text{ and } \alpha_{\nu'} = 0 \text{ for all } \nu' \neq \nu, \\ 0 & ; \text{otherwise.} \end{cases}$$

Our next claim is that, for every $\beta \in \mathbb{N}^d$, there is a unique element $S_\beta$ of $\mathbb{N}^d$ such that $\left\|\mathbf{b}'^{\beta}\right\|_{\mathbf{h},r} = \left\|t_{\beta,S_\beta}\, \mathbf{b}^{S_\beta}\right\|_{\mathbf{h},r}$. First observe that the assumption on $r'$ implies that, for every $1 \leq \nu \leq d$, $r_\nu^{p^{\gamma(\nu)}-1} > |p|$. Using the fact that $p$ divides $\binom{p^{\gamma(\nu)}}{k}$ for all $1 \leq k < p^{\gamma(\nu)}$ we obtain

$$\left|\binom{p^{\gamma(\nu)}}{p^{\gamma(\nu)}}\right| r_\nu^{p^{\gamma(\nu)}} = r_\nu^{p^{\gamma(\nu)}} > \left|\binom{p^{\gamma(\nu)}}{k}\right| r_\nu \geq \left|\binom{p^{\gamma(\nu)}}{k}\right| r_\nu^k$$

for all $1 \leq k < p^{\gamma(\nu)}$. This proves our claim in case $\beta = (0, \ldots, 1, \ldots, 0)$ is the $\nu$-th unit vector: letting $S_\beta := (0, \ldots, p^{\gamma(\nu)}, \ldots, 0) = p^{\gamma(\nu)}\,\beta$, we have

$$\left\|t_{\beta,S_\beta}\, \mathbf{b}^{\beta}\right\|_{\mathbf{h},r} > \left\|t_{\beta,\alpha}\, \mathbf{b}^{\alpha}\right\|_{\mathbf{h},r} \quad \text{for all } \alpha \neq S_\beta,$$

$$\left\|\mathbf{b}'^{\,\beta}\right\|_{\mathbf{h},r} = \max_\alpha \left\|t_{\beta,\alpha}\, \mathbf{b}^{\alpha}\right\|_{\mathbf{h},r} = \left\|t_{\beta,S_\beta}\, \mathbf{b}^{\beta}\right\|_{\mathbf{h},r} = r_\nu^{p^{\gamma(\nu)}}.$$

For general $\beta$ we have, by multiplicativity of the norm $\|\ \|_{\mathbf{h},r}$,

$$(**)\qquad \left\|\mathbf{b}'^{\beta}\right\|_{\mathbf{h},r} = \prod_{1\leq\nu\leq d} \|b'_\nu\|_{\mathbf{h},r}^{\beta_i} = \prod_{1\leq\nu\leq d} r_\nu{}^{p^{\gamma(\nu)}\,\beta_\nu} = r'^{\beta}.$$

On the other hand,

$$|t_{\beta,\alpha}| = \prod_{1\le\nu\le d} \left| \sum_{l\ge 0} (-1)^{\beta_\nu - l} \binom{\beta_\nu}{l} \binom{l\, p^{\gamma(\nu)}}{\alpha_\nu} \right| \begin{cases} = 1 & ;\text{ if } \alpha = \left(\beta_1\, p^{\gamma(1)}, \dots, \beta_d\, p^{\gamma(d)}\right) \\ < 1 & ;\text{ otherwise} \end{cases}.$$

by Lemma 6.2. It follows that $S_\beta := \left(\beta_1\, p^{\gamma(1)}, \dots, \beta_d\, p^{\gamma(d)}\right) \in \mathbb{N}^d$ is the unique element such that

$$\left\| t_{\beta,S_\beta}\, \mathbf{b}^{S_\beta} \right\|_{\mathbf{h},r} = r^{S_\beta} = r'^{\beta} = \left\| \mathbf{b}'^{\beta} \right\|_{\mathbf{h},r},$$

and our claim is established.

Since the elements $S_\beta$ ($\beta \in \mathbb{N}^d$) are pairwise different it follows from Lemma 6.1 that the family $\left(\mathbf{b}'^{\beta}\right)_{\beta\in\mathbb{N}}$ is orthogonal in the $K$-Banach space $D_{\mathbf{h},r}(H, K)$. This means that

$$\left\| \sum_\alpha d_\alpha\, \mathbf{b}'^{\alpha} \right\|_{\mathbf{h},r} = \max_\alpha |d_\alpha| \left\| \mathbf{b}'^{\alpha} \right\|_{\mathbf{h},r}$$

for every convergent series $\sum_\alpha d_\alpha\, \mathbf{b}'^{\alpha}$ in $D(H, K)$ and a fortiori for our convergent series $(*)$ in $D(H', K) \subset D(H, K)$. Together with $(**)$ this implies

$$\left\| \sum_\alpha d'_\alpha\, \mathbf{b}'^{\alpha} \right\|_{\mathbf{h},r} = \max_\alpha |d_\alpha|\, r'^{\alpha} = \left\| \sum_\alpha d_\alpha\, \mathbf{b}'^{\,\alpha} \right\|_{\mathbf{h}',r'}$$

as was to be shown. □

## 7. Distribution algebras of $L$-analytic groups

Let $H$ be an $L$-analytic group which is a uniform pro- $p$ group, and let $r \in p^{\mathbb{Q}} \cap \left[p^{-1}, 1\right)$. The map $D\left(R^L_{\mathbb{Q}_p}(H), K\right) \to D(H, K)$ dual to the canonical injection $C^{\mathrm{an}}(H, K) \to C^{\mathrm{an}}\left(R^L_{\mathbb{Q}_p}(H), K\right)$ is an epimorphism of $K$- Fréchet algebras ([31], Section 1). We let $\|\ \|_{\bar r}$ denote the residue norm on $D(H, K)$ induced from the norm $\|\ \|_r$ on $D\left(R^L_{\mathbb{Q}_p}(H), K\right)$ via this epimorphism, and we let $D_r(H, K)$ (resp. $U_r(\mathfrak{h}, K)$) denote the completion of $D(H, K)$ (resp. of $U(\mathfrak{h}, K)$) with respect to the norm $\|\ \|_{\bar r}$.

The property of an $L$-analytic group $H$ of being uniform pro- $p$ only depends on its underlying topological group, not on its being viewed as an $L$-analytic or $\mathbb{Q}_p$-analytic group. There is a finer concept of uniformity adapted to the circumstances of $L$-analytic groups (cf. [27], Section 2):

**Definition 7.1.** An $L$-analytic group $H$ with Lie algebra $\mathfrak{h}$ is called $L$- *uniform* if it is a uniform pro- $p$-group and the $\mathbb{Z}_p$-lattice $\mathrm{Log}(H) \subset \mathfrak{h}$ is an $o_L$- lattice.

If $H$ is $L$-uniform, $(\mathfrak{x}_1, \ldots, \mathfrak{x}_d)$ is an $o_L$-basis of $\mathrm{Log}(H)$ and $(v_1, \ldots, v_n)$ is a $\mathbb{Z}_p$-basis of $o_L$ then the family

$$(h_{ij} := \mathrm{Exp}\, v_i\, \mathfrak{x}_j)_{1\le i\le n, 1\le j\le d}$$

is an ordered basis of $H$; conversely, any $L$-analytic group $H$ which is a uniform pro- $p$- group with an ordered basis of the above form is $L$-uniform. According to [20], Proposition 1.3.5, every $L$-analytic group possesses a basis of neighbourhoods of 1 consisting of open $L$- uniform subgroups. We remark that if $H$ is $L$-uniform with ordered basis $(h_{ij})$ as above then every uniform pro- $p$-group with an ordered basis of the form $\left({h_{ij}}^{p^{\alpha(j)}}\right)_{i,j}$ with $\alpha(1), \ldots, \alpha(d) \in \mathbb{N}$, and in particular every lower $p$-series member of $H$, is $L$-uniform as well.

Recall that $\kappa = 1$ if $p \neq 2$, $\kappa = 2$ if $p = 2$. For $r \in (0, 1)$ define

$$\epsilon(r) := \sup\left\{m_0 \in \mathbb{N}_{\ge 1};\ \left|m_0^{-1}\right| r^{\kappa\, m_0} \ge \left|m^{-1}\right| r^{\kappa\, m} \text{ for all } m \in \mathbb{N}_{\ge 1}\right\}.$$

**Theorem 7.2 (Frommer-Kohlhaase).** Let $H$ be an $L$-uniform group with Lie algebra $\mathfrak{h}$. Let $r \in p^{\mathbb{Q}} \bigcap \left[p^{-1}, 1\right)$.

(i) Let $\mathfrak{X} = (\mathfrak{x}_1, \ldots, \mathfrak{x}_d)$ be an $o_L$-basis of $\mathrm{Log}(H)$. Then $(\mathfrak{X}^{\alpha})_{\alpha\in\mathbb{N}^d}$ is an orthogonal basis of $U_r(\mathfrak{h}, K)$.

(ii) $D_r(H, K)$ is a free right $U_r(\mathfrak{h}, K)$-module of rank $n \cdot d \cdot \epsilon(r)$. More exactly: Let $(h_1, \ldots, h_{n\,d})$ be an ordered basis of $H$ and put $\mathbf{b} := (\delta_{h_1} - 1, \ldots, \delta_{h_{n\,d}} - 1)$. Then $\left(\mathbf{b}^{\alpha}\right)_{\alpha\in\mathbb{N}^{n\,d}, \alpha_1<\epsilon(r), \ldots, \alpha_{n\,d}<\epsilon(r)}$ is a $U_r(\mathfrak{h}, K)$-basis of $D_r(H, K)$.

*Proof.* This is essentially Theorem 1.4.2 of [20]. The orthogonal basis property of the family $(\mathfrak{X}^{\alpha})$ for the norm $\|\ \|_{\bar{r}}$ is proved in [27], Proposition 2.4. For the description of $\epsilon(r)$ cf. [16], Section 1.4. □

**Remark 7.3.** Let $H$ be an $L$-analytic group which is a uniform pro- $p$-group, and let $r \in p^{\mathbb{Q}} \bigcap \left[p^{-1}, 1\right)$.

(1) Let $g \in H$. Then the series $\log \delta_g$ converges in $D_r(H, K)$, and

$$\log \delta_g = \mathrm{Log}\, g.$$

(2) What is the significance of the number $\epsilon(r)$? We have:

a) if $(\mathfrak{x}_1, \ldots, \mathfrak{x}_d)$ is a $\mathbb{Z}_p$-basis of $\mathrm{Log}(H)$ then $\epsilon(r) = \mathrm{dom}_r(\mathfrak{x}_i)$ for each $1 \leq i \leq d$; in particular $\|\mathfrak{x}_i\|_{\bar r} = \epsilon(r)\, r^{\kappa\, \epsilon(r)}$;

b) $\epsilon(r)$ is a $p$-power, and $\epsilon\left(r^{p^{-k}}\right) = p^k\, \epsilon(r)$ for $k \in \mathbb{N}$;

c) $|p|^{p/(p-1)} \leq r^{\kappa\, \epsilon(r)} < |p|^{1/(p-1)}$, and $\epsilon(r) \to \infty$ monotonely for $r \to 1$;

d) Let $g \in H^{\epsilon(r)}$, $\mathfrak{x} := \mathrm{Log}\, g$. Then the series $\exp \mathfrak{x}$ converges in $D_r(H, K)$. In particular, $\delta_g = \exp \mathfrak{x} \in U_r(\mathfrak{h}, K)$ and

$$\exp \mathfrak{x} = \delta_{\mathrm{Exp}\, \mathfrak{x}}.$$

*Proof.* (1) Let $(h_1, \ldots, h_d)$ be an ordered basis of $H$, and write as usual $\mathbf{b} = (b_1, \ldots, b_d)$, $b_i = \delta_{h_i} - 1$. According to Corollary 5.4 we may assume $g = h_1{}^{p^m}$ $(m \in \mathbb{N})$. Hence $\|\delta_g - 1\|_{\bar r} = \left\|\sum_{k\geq 1} \binom{p^m}{k} b_1^k\right\|_{\bar r} \leq r < 1$, so the series $\log \delta_g$ converges ([7], Section II.8.4).

We identify the Lie algebra $\mathfrak{h}$ of $H$ with the Lie algebra of $R^L_{\mathbb{Q}_p}(H)$. Then $\mathrm{Log}_H$ coincides with $\mathrm{Log}_{R^L_{\mathbb{Q}_p}(H)}$, and the epimorphism $D\left(R^L_{\mathbb{Q}_p} H, K\right) \to D(H, K)$ is compatible with the embeddings $H \to D\left(R^L_{\mathbb{Q}_p} H, K\right)$, $H \to D(H, K)$ and with the embeddings $\mathfrak{h} \to D\left(R^L_{\mathbb{Q}_p} H, K\right)$, $\mathfrak{h} \to D(H, K)$. Hence by continuity of this epimorphism, in order to show $\log \delta_g = \mathrm{Log}\, g$ we may restrict to the case $L = \mathbb{Q}_p$. But then

$$\mathrm{Log}(g^t) = t\, \mathrm{Log}(g), \qquad \log\left(\delta_g^t\right) = t \log \delta_g$$

for all $t \in \mathbb{Z}$ and hence for all $t \in \mathbb{Z}_p$. For sufficiently small $t$ we deduce

$$\delta_{\mathrm{Exp}(t\, \mathrm{Log}\, g)} = \exp(t \log \delta_g) = \sum_{n\geq 0} \binom{t}{n} (\delta_g - 1)^n.$$

Using the formula $\lim_{t\to 0} \frac{1}{t} \binom{t}{n} = \frac{(-1)^{n-1}}{n}$ we compute

$$\log \delta_g = \sum_{n\geq 1} \frac{(-1)^{n-1}}{n} (\delta_g - 1)^n = \lim_{t\to 0} \frac{1}{t} \sum_{n\geq 1} \binom{t}{n} (\delta_g - 1)^n = \lim_{t\to 0} \frac{1}{t} (\delta_{\mathrm{Exp}(t\, \mathrm{Log}\, g)} - 1);$$

the last expression is precisely $\mathrm{Log}\, g$ (cf. the description of the embedding $\mathfrak{h} \subset D(H, K)$ at the end of Section 2 in [30]).

(2) a) This is clear by definition of $\epsilon(r)$, since $\mathfrak{x}_i = \sum_{n>1} (-1)^{n-1} \frac{1}{n} b_i^n$, $b_i := \delta_{\mathrm{Exp}\, \mathfrak{x}_i} - 1$, by part (1).

b) First observe that the expression $\left|m^{-1}\right| r^{\kappa m}$ $(m \in \mathbb{N}_{\geq 1})$ attains its maximal values at $p$-powers $m$. Thus we may describe $\epsilon(r)$ as the maximum of the set

$$\{p^{k_0};\, k_0 \in \mathbb{N},\, p^{k_0}\, r^{\kappa\, p^{k_0}} \geq p^k\, r^{\kappa\, p^k} \text{ for all } k \in \mathbb{N}\}.$$

Let $k_1 \in \mathbb{N}$. The condition "$p^{k_0}\, r^{\kappa\, p^{k_0}} \geq p^k\, r^{\kappa\, p^k}$ for all $k \in \mathbb{N}$" is equivalent to "$p^{k_0}\, r^{\kappa\, p^{k_0}} \geq p^{k-k_1}\, r^{\kappa\, p^{k-k_1}}$ for all $k \in \mathbb{N}$" (since $p^k\, r^{\kappa\, p^k} \leq p^0\, r^{\kappa\, p^0}$ for negative $k$), and further to "$p^{k_0+k_1}\left(r^{p^{-k_1}}\right)^{\kappa\, p^{k_0+k_1}} \geq p^k\left(r^{p^{-k_1}}\right)^{\kappa\, p^k}$ for all $k \in \mathbb{N}$". This means that $\epsilon\left(r^{p^{-k_1}}\right) = \epsilon(r)\, p^{k_1}$.

c) Write $\kappa\, \epsilon(r) = p^k$. By definition, $\left|p^{-k}\right| r^{p^k} > \left|p^{-(k+1)}\right| r^{p^{k+1}}$, hence $p^{-1} = |p| > r^{(p-1)\, p^k}$, hence $p^{-1/(p-1)} > r^{p^k}$. Also by definition, $\left|p^{-k}\right| r^{p^k} \geq \left|p^{-(k-1)}\right| r^{p^{k-1}}$, hence $r^{\left(1-p^{-1}\right) p^k} \geq |p| = p^{-1}$, hence $r^{p^k} \geq p^{-1/\left(1-p^{-1}\right)} = p^{p/(1-p)}$. This also implies $\epsilon(r) \to \infty$ for $r \to 1$. In order to show monotony, let $r' \leq r$, and write $\epsilon(r) = p^k$, $\epsilon(r') = p^{k'}$, $r' = r\, s$. By definition, $p^{k'}\, r^{\kappa\, p^{k'}}\, s^{\kappa\, p^{k'}} \geq p^k\, r^{\kappa\, p^k}\, s^{\kappa\, p^k} \geq p^{k'}\, r^{\kappa\, p^{k'}}\, s^{\kappa\, p^k}$, which implies $k' \leq k$.

d) According to [7], Section II.8.4, we have to show that $\|\mathfrak{x}\|_{\bar{r}} < |p|^{1/(p-1)}$. There is a basis $(\mathfrak{x}_i)$ of $\mathrm{Log}(H)$ satisfying $\mathfrak{x} = p^m\, \mathfrak{x}_1$ (Corollary 5.4), and by assumption $p^m \geq \epsilon(r)$. By part a), $\|\mathfrak{x}\|_{\bar{r}} = p^{-m}\, \|\mathfrak{x}_1\|_{\bar{r}} = p^{-m}\, \epsilon(r)\, r^{\kappa\, \epsilon(r)} \leq r^{\kappa\, \epsilon(r)}$, and our claim follows from part c). □

Now fix a *compact* $L$-analytic group $G$. Let $H \subset G$ be an open normal subgroup which is a uniform pro-$p$-group, and let $r \in p^{\mathbb{Q}} \cap \left[p^{-1}, 1\right)$. In Section 5 of [32] the algebra $D(G, K)$ is endowed with a norm depending on the pair $(H, r)$ in the following way: Let $\mathcal{R}$ be a system of representatives of $G/H$. Let $\lambda \in D\left(R^L_{\mathbb{Q}_p}(G), K\right)$, and write

$$\lambda = \sum_{x \in \mathcal{R}} \lambda_x\, \delta_x$$

according to the decomposition $D\left(R^L_{\mathbb{Q}_p}(G), K\right) = \oplus_{x \in \mathcal{R}} D\left(R^L_{\mathbb{Q}_p}(H), K\right) \delta_x$. Then

$$\lambda \mapsto \max_{x \in \mathcal{R}} \|\lambda_x\|_r$$

defines a norm on $D\left(R^L_{\mathbb{Q}_p}(G), K\right)$ which is independent of the choice of $\mathcal{R}$ (cf. [26], Section 4.2). We let $\|\;\|_{(H,r)}$ denote the residue norm on $D(G, K)$ induced from this norm via the epimorphism $D\left(R^L_{\mathbb{Q}_p}(G), K\right) \to D(G, K)$, and we denote the completion of $D(G, K)$ with respect to this norm by $D_{(H,r)}(G, K)$. Then

$$D_{(H,r)}(G, K) = \bigoplus_{x\in\mathcal{R}} D_r(H, K)\,\delta_x$$

as $K$-Banach spaces (cf. [23], Remark 2.2.5). By [32], Theorem 5.1 (and proof), $D_{(H,r)}(G, K)$ is a noetherian $K$-Banach algebra.

Recall that a *basis* of a $K$-Banach space $(M, \|\ \|)$ is a family $(T_i)_{i\in I}$ in $M$ such that every $\mu \in M$ has a representation as a convergent sum $\mu = \sum_{i\in I} a_i\, T_i$ with a uniquely determined family $(a_i)_{i\in I}$ in $K$ (cf. [24], p. 53).

**Corollary 7.4.** Let $G$ be a compact $L$-analytic group, $H \subset G$ an $L$-uniform open normal subgroup, and $r \in p^{\mathbb{Q}} \cap [p^{-1}, 1)$. Let $\mathfrak{X} = (\mathfrak{x}_1, \ldots, \mathfrak{x}_d)$ be an $o_L$-basis of $\mathrm{Log}(H)$. Let $\mathcal{R}$ be a system of representatives of $G/H^{\epsilon(r)}$. Then the family $(\mathfrak{X}^\alpha\, \delta_x)_{\alpha\in\mathbb{N}^d, x\in\mathcal{R}}$ is a basis of the $K$-Banach space $D_{(H,r)}(G, K)$.

*Proof.* Let $(h_i)$ be an ordered basis of $H$, $\mathbf{h} := \left(\delta_{h_i}\right)_i$, $\mathbf{b} := \left(\delta_{h_i} - 1\right)_i$. In Theorem 7.2 (ii) we may replace the $U_r(\mathfrak{g}, K)$-basis $\left(\mathbf{b}^\alpha\right)_{\alpha<\epsilon(r)}$ of $D_r(H, K)$ by the basis $\left(\mathbf{h}^\alpha\right)_{\alpha<\epsilon(r)}$, since

$$\left(\mathbf{h}^\alpha\right)_{\alpha<\epsilon(r)} = \left(\mathbf{b}^\alpha\right)_{\alpha<\epsilon(r)} \cdot T, \qquad \left(\mathbf{b}^\alpha\right)_{\alpha<\epsilon(r)} = \left(\mathbf{h}^\alpha\right)_{\alpha<\epsilon(r)} \cdot T^{-1}$$

with base change matrices $T = (\binom{\alpha}{\beta})_{\alpha,\beta<\epsilon(r)}$, $T^{-1} = ((-1)^{|\beta|}\binom{\alpha}{\beta})_{\alpha,\beta<\epsilon(r)}$. Let $\mathcal{R}'$ be any system of representatives of $G/H$. From the decomposition $D_{(H,r)}(G, K) = \oplus_{y\in\mathcal{R}'} D_r(H, K)\,\delta_y$ we obtain that $\left(\mathbf{h}^\alpha\, \delta_y\right)_{\alpha<\epsilon(r), y\in\mathcal{R}'}$ is a $U_r(\mathfrak{g}, K)$-basis of $D_{(H,r)}(G, K)$. Now if $\mathbf{h}^\alpha\, \delta_y$ ($\alpha < \epsilon(r)$, $y \in \mathcal{R}'$) and $\delta_x$ ($x \in \mathcal{R}$) are two representatives contained in the same coset of $H^{\epsilon(r)}$ then $\mathbf{h}^\alpha\, \delta_y\, \delta_x^{-1}$ is a unit in $U_r(\mathfrak{g}, K)$ (Remark 7.3 (2) d)). Hence $(\delta_x)_{x\in\mathcal{R}}$ is also a $U_r(\mathfrak{g}, K)$-basis of $D_{(H,r)}(G, K)$. The corollary follows then from Theorem 7.2 (i). □

**Proposition 7.5.** Let $G$ be a compact $L$-analytic group. Let $H, H' \subset G$ be open normal subgroups which are uniform pro-$p$-groups. Let $r, r' \in p^{\mathbb{Q}} \cap [p^{-1}, 1)$. Suppose that one of the following conditions holds:

(i) $H' = H$ and $r' \geq r$;

(ii) $H' \subset H$, $r' = r$, $r^{\kappa p^m} \geq p^{-1}$ where $m$ is the highest $p$-elementary divisor of $H'$ in $H$;

(iii) there is an $m \in \mathbb{N}$ such that $H' = H^{p^m}$, $r' = r^{p^m}$.

Then there is a continuous $K$-algebra homomorphism $D_{(H',r')}(G, K) \to D_{(H,r)}(G, K)$ extending the identity map $D(G, K) \to D(G, K)$. In case (iii) this homomorphism is a topological isomorphism.

*Proof.* First note that a continuous $K$-linear map $D_{(H',r')}(G, K) \to D_{(H,r)}(G, K)$ fixing $G$ (that is, fixing the Dirac distributions $\delta_g$ for $g \in G$) is automatically a $K$-algebra homomorphism extending the identity on $D(G, K)$.

(i) We have $\| \ \|_r \le \| \ \|_{r'}$ on $D\big(R^L_{\mathbb{Q}_p}(H), K\big)$, hence $\| \ \|_{\bar{r}} \le \| \ \|_{\bar{r}'}$ on $D(H, K)$, hence the identity on $D(H, K)$ extends to a norm-decreasing map $D_{r'}(H, K) \to D_r(H, K)$. This induces a norm-decreasing $K$-linear map $D_{(H,r')}(G, K) = \oplus_{x\in G_0/H} D_{r'}(H, K)\,\delta_x \to \oplus_{x\in G_0/H} D_r(H, K)\,\delta_x = D_{(H,r)}(G, K)$ which fixes $G$.

(ii) Let $(\alpha(1), \ldots, \alpha(d) = m)$ be the sequence of $p$-elementary divisors of $H'$ in $H$. By assumption we have $r'' := (r^{\kappa\, p^{\alpha(1)}}, \ldots, r^{\kappa\, p^{\alpha(d)}}) \in I(d)$, so we may apply Proposition 6.3 to obtain a an ordered basis $\mathbf{h}'$ of $H'$ such that the norm $\| \ \|_{\mathbf{h}',r''}$ on $D\big(R^L_{\mathbb{Q}_p}(H'), K\big)$ is equal to the restriction of the norm $\| \ \|_r$ on $D\big(R^L_{\mathbb{Q}_p}(H), K\big)$ to $D\big(R^L_{\mathbb{Q}_p}(H'), K\big)$. It follows from [26], Corollary 4.9, that the residue norm on $D(H', K)$ induced from $\| \ \|_{\mathbf{h}',r''}$ is equivalent to the restriction of the norm $\| \ \|_{\bar{r}}$ on $D(H, K)$ to $D(H', K)$, and for the respective completions we have

$$D_r(H, K) = \bigoplus_{x\in H/H'} D_{\mathbf{h}',r''}(H', K)\,\delta_x$$

as $K$-Banach spaces. On the other hand there is, similarly as in part (i) of this proof, the natural continuous $K$-linear map $D_r(H', K) \to D_{\mathbf{h}',r''}(H', K)$, and thereby the continuous $K$-linear map

$$D_{(H',r)}(G, K) = \left(\bigoplus_{y\in G_0/H}\ \bigoplus_{x\in H/H'} D_r(H', K)\,\delta_y\,\delta_x\right) \longrightarrow \left(\bigoplus_{y\in G_0/H} D_r(H, K)\,\delta_y\right) = D_{(H,r)}(G, K),$$

fixing $G$.

(iii) The restriction of the norm $\| \ \|_r$ on $D\big(R^L_{\mathbb{Q}_p}(H), K\big)$ to $D\big(R^L_{\mathbb{Q}_p}\big(H^{p^m}\big), K\big)$ is equal to the norm $\| \ \|_{r^{p^m}}$ on $D\big(R^L_{\mathbb{Q}_p}\big(H^{p^m}\big), K\big)$ (Proposition 6.3), hence ([26], Corollary 4.9) $D_r(H, K) = \oplus_{x\in H/H^{p^m}} D_{r^{p^m}}\big(H^{p^m}, K\big)$. Consequently, $D_{\left(H^{p^m}, r^{p^m}\right)}(G, K) = \oplus_{y\in G_0/H} \oplus_{x\in H/H^{p^m}} D_{r^{p^m}}\big(H^{p^m}, K\big)\,\delta_y\,\delta_x = \oplus_{y\in G_0/H} D_r(H, K)\,\delta_y = D_{(H,r)}(G, K)$. □

## 8. Decreasing filtrations of the enveloping algebra

Let $G$ be an $L$-analytic group with Lie algebra $\mathfrak{g}$ and exponential mapping Exp. Let $V$ be a locally analytic $G$-representation. The Lie algebra action on $V$ is given by

$$\mathfrak{x}\, v = \tfrac{d}{dt} \operatorname{Exp}(t\, \mathfrak{x})\, v\,|_{t=0}$$

($\mathfrak{x} \in \mathfrak{g}$, $v \in V$); it naturally extends to an action of the enveloping algebra $U(\mathfrak{g})$ of $\mathfrak{g}$. For $k \in \mathbb{N}$ let $\mathfrak{g}^k \subset U(\mathfrak{g})$ denote the subset of elements of the form $\mathfrak{x}_1 \dots \mathfrak{x}_k$ with $\mathfrak{x}_1, \dots, \mathfrak{x}_k \in \mathfrak{g}$.

**Definition 8.1.** A vector $v \in V$ is called $\mathfrak{g}$-*adically finite* if $\mathfrak{g}^k v = 0$ for some $k \in \mathbb{N}$.

From the formula $\mathfrak{x}\, g\, v = g \operatorname{Ad}(g)^{-1}(\mathfrak{x})\, v$ ($\mathfrak{x} \in \mathfrak{g}$, $g \in G$, $v \in V$, cf. [15] Satz 3.1.3) it follows that the $\mathfrak{g}$-adically finite vectors form an invariant subspace of $V$.

Let $T(\mathfrak{g}) = \oplus_{k\in\mathbb{N}} T^k(\mathfrak{g})$ be the graded tensor $L$-algebra of $\mathfrak{g}$, with $T^k(\mathfrak{g})$ denoting the subspace of homogeneous tensors of order $k$ ($k \in \mathbb{N}$). We have the canonical epimorphism

$$\pi : T(\mathfrak{g}) \to U(\mathfrak{g}),\ \mathfrak{z}_1 \otimes \dots \otimes \mathfrak{z}_k \mapsto \mathfrak{z}_1 \dots \mathfrak{z}_k,$$

and $\pi\big(T^k(\mathfrak{g})\big)$ is equal to the $L$-vector subspace of $U(\mathfrak{g})$ generated by $\mathfrak{g}^k$. The sum of the $\pi\big(T^k(\mathfrak{g})\big)$ is of course not direct anymore, we rather have $\pi\big(T^k(\mathfrak{g})\big) \subset \sum_{1\le k'\le k} L\, \mathfrak{g}^{k'}$ for $k \ge 1$. The $L$- vector spaces

$$U_m(\mathfrak{g}) := \sum_{k\le m} \pi\big(T^k(\mathfrak{g})\big) \quad (m \in \mathbb{N})$$

constitute the (increasing) *canonical filtration* of the $L$-algebra $U(\mathfrak{g})$ considered in [7], I.2.6, and [13], 2.3.1. On the other hand, for each $m \in \mathbb{N}$ we obviously have

$$\mathfrak{g}^m\, U(\mathfrak{g}) := \{\textstyle\sum_i \mathfrak{x}_i\, \mathfrak{y}_i;\ \mathfrak{x}_i \in \mathfrak{g}^m,\ \mathfrak{y}_i \in U(\mathfrak{g})\} = \sum_{k\ge m} \pi\big(T^k(\mathfrak{g})\big),$$

and this $L$-vector space is equal to the left, as well as to the two-sided ideal in $U(\mathfrak{g})$ generated by $\mathfrak{g}^m$. As the spaces $\pi\big(T^k(\mathfrak{g})\big)$ are finite-dimensional the ideals $\mathfrak{g}^m\, U(\mathfrak{g})$ are of finite codimension in $U(\mathfrak{g})$. Moreover, $\mathfrak{g}^k\, U(\mathfrak{g}) \cdot \mathfrak{g}^m\, U(\mathfrak{g}) \subset \mathfrak{g}^{k+m}\, U(\mathfrak{g})$ for all $k, m \in \mathbb{N}$.

**Definition 8.2.** The decreasing filtration $\big(\mathfrak{g}^k\, U(\mathfrak{g})\big)_{k\in\mathbb{N}}$ is called the $\mathfrak{g}$-*adic filtration* of the algebra $U(\mathfrak{g})$.

Clearly, a vector is $\mathfrak{g}$-adically finite if and only if it is annihilated by the ideal $\mathfrak{g}^k\, U(\mathfrak{g})$ for some $k \in \mathbb{N}$. Recall that a vector $v$ is called $U(\mathfrak{g})$- *finite* (cf. [14], Definition 4.1.10) if it is contained in a finite-dimensional $U(\mathfrak{g})$-stable subspace of $V$, or equivalently: if it is annihilated by a left ideal in $U(\mathfrak{g})$ of finite codimension. By a lemma of Emerton ([14] 4.1.11) this is equivalent to the condition that $v$ is *locally finite*, i.e. contained in a finite-dimensional subspace of $V$ stable under an open subgroup of $G$. Obviously, every $\mathfrak{g}$-adically finite vector is $U(\mathfrak{g})$-finite. The converse does not hold as Example 8.3 below indicates.

**Example 8.3.** Let $G = o_L$, the additive group of integers of $L$. Then Exp is the identity map on $o_L$. Let $V$ be the locally analytic regular representation of $G$, that is, $V$ is the locally convex $K$-vector space $C^{\mathrm{an}}(G, K)$ endowed with the translation action of $G$. The Lie algebra action is given by $\mathfrak{x}_1\, f = -f'$ ($\mathfrak{x}_1 = 1 \in L$, $f \in V$; here $f'$ denotes the dervative of $f$).

(a) A vector $f \in V$ is $\mathfrak{g}$-adically finite if and only if it is *locally polynomial*, i.e. if and only if for each $x \in G$ there is a polynomial $P$ with coefficients in $K$ such that $f(z) = P(y - x)$ for all $y$ in a neighbourhood of $x$.

(b) Choose an open subgroup $U \subset G$ such that the series $\sum_{m\in\mathbb{N}} x^m/m!$ converges for $x \in U$. The vector

$$e: G \longrightarrow K$$
$$x \longmapsto \begin{cases} \sum_{m\in\mathbb{N}} x^m/m! & ; x \in U \\ 0 & ; x \notin U \end{cases}$$

is not $\mathfrak{g}$-adically finite but is $U(\mathfrak{g})$-finite: We have $\mathfrak{x}_1\, e = -e$, hence $\mathfrak{g}^k\, e = L\, e \neq 0$ for all $k \in \mathbb{N}$, while the ideal generated by $\mathfrak{x}_1 + 1$ in $U(\mathfrak{g})$ is of finite codimension and annihilates $e$.

More generally one can show the following: Let $G$ be the group of $L$-rational points of a unipotent algebraic $L$-group and let $V$ be the left regular representation on $C^{\mathrm{an}}(G, K)$. A vector $f \in V$ is $\mathfrak{g}$-adically finite if and only if it is locally polynomial with respect to the canonical chart of the second species associated to any ordered basis of $\mathfrak{g}$.

We denote by $\mathfrak{g} = C^0(\mathfrak{g}) \supset C^1(\mathfrak{g}) \supset \ldots$ the lower central series of $\mathfrak{g}$.

**Lemma 8.4.** Let $G$ be the group of $L$-rational points of a unipotent $L$-group. For any $g \in G$ and $m \in \mathbb{N}$,

$$(\mathrm{Ad}(g) - \mathrm{id}_{\mathfrak{g}})\,(\mathbf{C}^m(\mathfrak{g})) \subset \mathbf{C}^{m+1}(\mathfrak{g}).$$

*Proof.* Let $\mathfrak{x} = \mathrm{Log}(g)$. We may view $G$ as a group of strictly upper triangular matrices in some $\mathrm{GL}_m$ (Example 4.2). Hence Exp is defined on the whole of $\mathfrak{g}$, and we have the equality $\mathrm{Ad}(\mathrm{Exp}(\mathfrak{x})) = \exp(\mathrm{ad}(\mathfrak{x}))$ in $\mathrm{GL}(\mathfrak{g})$. For all $\mathfrak{y} \in \mathbf{C}^m(\mathfrak{g})$,

$$(\mathrm{Ad}(g) - \mathrm{id}_{\mathfrak{g}})\,(\mathfrak{y}) = \exp(\mathrm{ad}(\mathfrak{x}))\,(\mathfrak{y}) - \mathfrak{y} = [\mathfrak{x}, \mathfrak{y}] + \tfrac{1}{2}[\mathfrak{x}, [\mathfrak{x}, \mathfrak{y}]] + \ldots\,.$$

But $\mathbf{C}^m(\mathfrak{g})$ is an ideal in $\mathfrak{g}$, so $(\mathrm{Ad}(g) - \mathrm{id}_{\mathfrak{g}})\,(\mathfrak{y})$ is contained in $[\mathfrak{x}, \mathbf{C}^m(\mathfrak{g})] \subset \mathbf{C}^{m+1}(\mathfrak{g})$. □

For the rest of this section, we fix an open normal $L$-uniform subgroup $H \subset G$ and a number $r \in p^{\mathbb{Q}} \cap [p^{-1}, 1)$. Let $(\mathfrak{x}_1, \ldots, \mathfrak{x}_d)$ be an $o_L$-basis of $\mathrm{Log}(H)$ and $\mathcal{R}$ a system of representatives for $H^{\epsilon(r)} \backslash G$. Then the family

$$\left(\mathfrak{x}_1^{\alpha_1} \ldots \mathfrak{x}_d^{\alpha_d}\, \delta_g\right)_{\alpha \in \mathbb{N}^d, g \in \mathcal{R}}$$

is a basis of the $K$-Banach space $D_{(H,r)}(G, K)$ (Corollary 7.4). We further assume $\mathcal{R}$ to be of the following form: There exists a family $\mathbf{g} = (g_{i,\beta})_{1 \le i \le d, \beta \in I}$ in $G$ such that $\mathcal{R} = \{g_{1,\beta} \cdot \ldots \cdot g_{d,\beta};\ \beta \in I\}$. For $k \in \mathbb{N}$, let $M(\mathbf{g}, k) \subset D_{(H,r)}(G, K)$ denote the closed $K$-vector subspace generated by the elements

$$\mathfrak{x}_1^{\alpha_1}\, \delta_{g_{1,\beta}} \ldots \mathfrak{x}_d^{\alpha_d}\, \delta_{g_{d,\beta}} \qquad (\beta \in I,\ \alpha \in \mathbb{N}^d, \text{ not all } \alpha_i < k).$$

The next two lemmata compare this space with the closed subspace

$$\mathfrak{g}^k D_{(H,r)}(G, K) = \left\{\textstyle\sum_i \mathfrak{x}_i\, \lambda_i;\ \mathfrak{x}_i \in \mathfrak{g}^k,\ \lambda_i \in D_{(H,r)}(G, K)\right\}$$

of $D_{(H,r)}(G, K)$.

**Lemma 8.5.** (i) Let $g \in G$, $\mathfrak{x} \in \mathfrak{g}$. In $D_{(H,r)}(G, K)$ we have

$$\delta_g\, \mathfrak{x}\, \delta_g^{-1} = \mathrm{Ad}(g)\,(\mathfrak{x}).$$

(ii) For any $k \in \mathbb{N}$,

$$M(\mathbf{g}, k) \subset \mathfrak{g}^k D_{(H,r)}(G, K).$$

*Proof.* There exists $c \in L^\times$ such that $c\,\mathfrak{x} \in \mathrm{Log}(H^{\epsilon(r)})$. Put $h := \mathrm{Exp}(c\,\mathfrak{x})$. By Remark 7.3, $\delta_g\, \mathfrak{x}\, \delta_g^{-1} = c^{-1}\, \delta_g \log(\delta_h)\, \delta_g^{-1} = c^{-1} \sum_{n \ge 1} (-1)^{n-1}\, n^{-1}\, \delta_g\, \delta_h^n\, \delta_g^{-1} = c^{-1} \sum_{n \ge 1} (-1)^{n-1}\, n^{-1} \left(\delta_{g\,h\,g^{-1}}\right)^n = c^{-1} \log\left(\delta_{g\,h\,g^{-1}}\right) = c^{-1}\, \mathrm{Log}\left(g\,h\,g^{-1}\right)$, since $g\,h\,g^{-1}$ is contained in $H^{\epsilon(r)}$. After eventually making $c$ smaller we may apply Corollary 3 (ii) of [7], III.4.4, to obtain $c^{-1}\, \mathrm{Log}\left(g\,h\,g^{-1}\right) = c^{-1}\, \mathrm{Ad}(g)\,(\mathrm{Exp}(h)) = \mathrm{Ad}(g)\,(\mathfrak{x})$, whence (i). In order to prove (ii), we have to show that $\mathfrak{x}_1^{\alpha_1}\, \delta_{g_{1,\beta}} \ldots \mathfrak{x}_d^{\alpha_d}\, \delta_{g_{d,\beta}}$ is contained in $\mathfrak{g}^k D_{(H,r)}(G, K)$ for $\beta \in I$, $\alpha \in \mathbb{N}^d$, not all $\alpha_i < k$. But from part (i) it follows that there exist $k' \ge k$ and $\lambda \in \mathfrak{g}^{k'}$ such that $\mathfrak{x}_1^{\alpha_1}\, \delta_{g_{1,\beta}} \ldots \mathfrak{x}_d^{\alpha_d}\, \delta_{g_{d,\beta}} = \lambda\, \delta_{g_{1,\beta}} \ldots \delta_{g_{d,\beta}}$. □

**Lemma 8.6.** Let $G$ be the group of $L$-rational points of a unipotent $L$-group, and choose $m \ge 2$ such that $C^m(\mathfrak{g}) = 0$. Suppose that each $C^i(\mathfrak{g})$ is spanned by a subfamily of the basis $(\mathfrak{x}_1, \ldots, \mathfrak{x}_d)$.

Let $k_1, k_2 \in \mathbb{N}$ such that $d\,m\,k_2 < k_1$. Then

$$\mathfrak{g}^{k_1}\, D_{(H,r)}(G, K) \subset M(\mathbf{g}, k_2).$$

*Proof.* By our assumption on $\mathcal{R}$, $\mathfrak{g}^{k_1}\, D_{(H,r)}(G, K)$ is the closed $K$-vector subspace of $D_{(H,r)}(G, K)$ generated by the elements

$$\lambda\, \mathfrak{x}_1^{\alpha_1} \ldots \mathfrak{x}_d^{\alpha_d}\, \delta_{g_{1,\beta}} \ldots \delta_{g_{d,\beta}} \quad (\lambda \in \mathfrak{g}^{k_1},\ \alpha \in \mathbb{N}^d,\ \beta \in I),$$

and hence also by the elements

$$\lambda\, \delta_{g_{1,\beta}} \ldots \delta_{g_{d,\beta}} \quad (\lambda \in \mathfrak{g}^{k},\ k \geq k_1,\ \beta \in I).$$

We fix such an element and introduce some notations: For $\mathfrak{z} \in \mathfrak{g}$ define $\mathcal{N}(\mathfrak{z}) := \max\left\{i \in \mathbb{N};\ \mathfrak{z} \in C^i(\mathfrak{g})\right\}$ if $\mathfrak{z} \neq 0$, $\mathcal{N}(0) := \infty$. Put

$$\mathcal{A} := \left\{\left((\mathfrak{z}_{i,j})_{1\leq j\leq k(i)}\right)_{1\leq i\leq d};\ \mathfrak{z}_{i,j} \in \mathfrak{g},\ k(i) \in \mathbb{N}\right\}.$$

Every family $\mathfrak{z} = \left((\mathfrak{z}_{i,j})_{1\leq j\leq k(i)}\right)_{1\leq i\leq d} \in \mathcal{A}$ gives rise to an element

$$\lambda(\mathfrak{z}) := \mathfrak{z}_{1,1} \ldots \mathfrak{z}_{1,k(1)}\, \delta_{g_{1,\beta}}\ \ldots\ \mathfrak{z}_{d,1} \ldots \mathfrak{z}_{d,k(d)}\, \delta_{g_{d,\beta}} \in D_{(H,r)}(G, K),$$

and we define $\mathcal{N}(\mathfrak{z}) := \sum_{i,j} \mathcal{N}(\mathfrak{z}_{i,j})$, $\mathcal{L}(\mathfrak{z}) := \sum_i k(i)$. Finally, let $\mathcal{B}$ denote the subset of $D_{(H,r)}(G, K)$ of elements of the form

$$\mu = c\, \mathfrak{x}_1^{\alpha_1}\, \delta_{g_{1,\beta}} \ldots \mathfrak{x}_d^{\alpha_d}\, \delta_{g_{d,\beta}} \quad (c \in K,\ \alpha \in \mathbb{N}^d),$$

and for such an element $\mu \in \mathcal{B}$ put $\mathcal{L}(\mu) = \sum_i \alpha_i$. Clearly, if $\mathcal{L}(\mu) \geq d\,k_2$ then $\mu \in M(\mathbf{g}, k_2)$. Thus it will suffice to show that our fixed generator $\lambda\, \delta_{g_{1,\beta}} \ldots \delta_{g_{d,\beta}}$ is a finite sum of elements $\mu_\nu \in \mathcal{B}$ with $\mathcal{L}(\mu_\nu) \geq d\,k_2$.

We claim the following:

(+) Let $\mathfrak{z} \in \mathcal{A}$ with $\mathcal{L}(\mathfrak{z}) = k$. There exists a finite family $(\mu_\nu)$ in $\mathcal{B}$ and a finite family $(\mathfrak{z}^{(\nu)})$ in $\mathcal{A}$ such that

(a) $\lambda(\mathfrak{z}) = \sum_\nu \mu_\nu + \sum_\nu \lambda(\mathfrak{z}^{(\nu)})$,

(b) $\mathcal{L}(\mu_\nu) = k$ for all $\nu$,

(c) $k - 1 \leq \mathcal{L}(\mathfrak{z}^{(\nu)}) \leq k$ and $\mathcal{N}(\mathfrak{z}^{(\nu)}) \geq \mathcal{N}(\mathfrak{z}) + 1$ for all $\nu$.

In order to prove (+), we first observe that $\lambda(\mathfrak{z})$ is a linear combination of elements of the form $\lambda(\mathfrak{y})$ with $\mathfrak{y} = (\mathfrak{y}_{i,j}) \in \mathcal{A}$, all $\mathfrak{y}_{i,j} \in \{\mathfrak{x}_1, \ldots, \mathfrak{x}_d\}$, $\mathcal{L}(\mathfrak{y}) = \mathcal{L}(\mathfrak{z})$ and (using our assumption on the basis) $\mathcal{N}(\mathfrak{y}) \geq \mathcal{N}(\mathfrak{z})$. Hence we may assume that $\lambda(\mathfrak{z})$ itself has this form. Next, for all $1 \leq i,\ j \leq d$ we have

(I) $\mathfrak{x}_j\, \mathfrak{x}_i = \mathfrak{x}_i\, \mathfrak{x}_j + [\mathfrak{x}_j, \mathfrak{x}_i]$ with $\mathcal{N}([\mathfrak{x}_j, \mathfrak{x}_i]) \geq \mathcal{N}(\mathfrak{x}_i) + \mathcal{N}(\mathfrak{x}_j) + 1$

and, according to Lemma 8.5 (i) and Lemma 8.4,

(II) $\mathfrak{x}_j\, \delta_{g_{i,\beta}} = \delta_{g_{i,\beta}}\, \mathfrak{x}_j - \mathfrak{y}\, \delta_{g_{i,\beta}}$ with $\mathcal{N}(\mathfrak{y}) \geq \mathcal{N}(\mathfrak{x}_j) + 1$,

where $\mathfrak{y} := \delta_{g_{i,\beta}}\, \mathfrak{x}_j\, \delta_{g_{i,\beta}}^{-1} - \mathfrak{x}_j = \mathrm{Ad}(g_{i,\beta})\,(\mathfrak{x}_j) - \mathfrak{x}_j$. Applying a finite number of operations of type (I) and (II) to $\lambda(\mathfrak{z})$ we indeed obtain

$$\lambda(\mathfrak{z}) = \mathfrak{x}_1^{\alpha_1}\, \delta_{g_{1,\beta}} \ldots \mathfrak{x}_d^{\alpha_d}\, \delta_{g_{d,\beta}} + \textstyle\sum_\nu \lambda(\mathfrak{z}^{(\nu)}) \quad (\alpha \in \mathbb{N}^d,\ \sum \alpha_i = k,\ \mathfrak{z}^{(\nu)} \in \mathcal{A})$$

such that the elements $\mathfrak{z}^{(\nu)}$ satisfy the condition (c) of (+).

To conclude the proof, we observe that our chosen generator $\lambda\, \delta_{g_{1,\beta}} \ldots \delta_{g_{d,\beta}}$ has the form $\lambda(\mathfrak{z})$ for some $\mathfrak{z} \in \mathcal{A}$ with $\mathcal{L}(\mathfrak{z}) = k$. Since $k \geq k_1 \geq d\, k_2$, and since the procedure (+) "increases $\mathcal{N}$" and "decreases $\mathcal{L}$" according to (c), we may apply it to $\lambda(\mathfrak{z})$ and then recursively to the summands $\lambda(\mathfrak{z}^{(\nu)})$ in (a) until we arrive at an equation

$$\lambda(\mathfrak{z}) = \textstyle\sum_\nu \mu_\nu + \sum_\nu \lambda(\mathfrak{z}^{(\nu)}) \quad (\mu_\nu \in \mathcal{B},\ \mathfrak{z}^{(\nu)} \in \mathcal{A})$$

with either $\mathcal{L}(\mathfrak{z}^{(\nu)}) = d\, k_2$, or $\mathcal{L}(\mathfrak{z}^{(\nu)}) > d\, k_2$ and $\mathcal{N}(\mathfrak{z}^{(\nu)}) > (m-1)\, \mathcal{L}(\mathfrak{z}^{(\nu)})$ for each $\nu$. In view of (b) this also implies $\mathcal{L}(\mu_\nu) \geq d\, k_2$ for each $\nu$. Moreover, if $\mathcal{L}(\mathfrak{z}^{(\nu)}) = d\, k_2$ then in virtue of (c) we have $\mathcal{N}(\mathfrak{z}^{(\nu)}) \geq k - d\, k_2 > (m-1)\, d\, k_2 = (m-1)\, \mathcal{L}(\mathfrak{z}^{(\nu)})$ as well in this case. But for each $\mathfrak{z}^{(\nu)} = \left(\mathfrak{z}_{i,j}^{(\nu)}\right)$ the relation $\mathcal{N}(\mathfrak{z}^{(\nu)}) > (m-1)\, \mathcal{L}(\mathfrak{z}^{(\nu)})$ implies that at least one $\mathfrak{z}_{i,j}^{(\nu)}$ is contained in $C^m(\mathfrak{g})$, and thereby $\lambda(\mathfrak{z}^{(\nu)}) = 0$ by the choice of $m$. Hence

$$\lambda(\mathfrak{z}) = \textstyle\sum_\nu \mu_\nu \quad \text{with } \mu_\nu \in \mathcal{B},\ \mathcal{L}(\mu_\nu) \geq d\, k_2. \ \square$$

**Remark 8.7.** The $\mathfrak{g}$-adic filtration $\left(\mathfrak{g}^k\, U(\mathfrak{g})\right)$ is in general not separated. In fact, in case $[\mathfrak{g}, \mathfrak{g}] = \mathfrak{g}$ it is even stationary with $\mathfrak{g}^k\, U(\mathfrak{g}) = \mathfrak{g}\, U(\mathfrak{g})$ for all $k \geq 1$. If on the other hand $\mathfrak{g}$ is nilpotent then a simplified version of Lemma 8.6 shows that the $\mathfrak{g}$-adic filtration is cofinal with the filtration given by the subspaces generated by the elements

$$\mathfrak{x}_1^{\alpha_1} \ldots \mathfrak{x}_d^{\alpha_d} \quad (\alpha \in \mathbb{N}^d, \text{ not all } \alpha_i < k)$$

($k \in \mathbb{N}$, $(\mathfrak{x}_1, \ldots, \mathfrak{x}_d)$ a fixed basis of $\mathfrak{g}$); and this filtration clearly is separated. Bourbaki considers yet another decreasing filtration which is independent of the choice of a particular basis: letting $S'^k(\mathfrak{g}) \subset T^k(\mathfrak{g})$ denote the subspace of *symmetric* tensors, the canonical map $\pi : T(\mathfrak{g}) \to U(\mathfrak{g})$ restricts to isomorphisms of $L$-vector spaces

$$S'^k(\mathfrak{g}) \xrightarrow{\sim} U^k(\mathfrak{g}) := \pi\big(S'^k(\mathfrak{g})\big),$$

($k \in \mathbb{N}$), and each $U_m(\mathfrak{g})$ is the direct sum of the $U^k(\mathfrak{g})$ ($k \leq m$) (compare [4], III.6.3, Remark; [7], I.2.4, Corollary 4 and subsequent discussion). Again, a proof similar to that of Lemma 8.6 shows that the $\mathfrak{g}$-adic filtration is cofinal with the filtration $\left(\sum_{k' \geq k} U^k(\mathfrak{g})\right)_{k \in \mathbb{N}}$ provided $\mathfrak{g}$ is nilpotent.

## III. $\mathfrak{n}$-characters of admissible representations of reductive $p$-adic groups

**Notations.** For the remaining sections we fix the following situation: Let $\mathfrak{G}$ be a connected reductive group over $L$. Let $\mathfrak{S} \subset \mathfrak{G}$ be a maximal $L$-split $L$-torus. Let $W$ be the Weyl group and $\Phi_{\mathrm{nd}}$ the system of non-divisible roots of $(\mathfrak{G}, \mathfrak{S})$. Let $\mathfrak{P} \subset \mathfrak{G}$ be a minimal parabolic subgroup containing $\mathfrak{S}$, and let $\mathfrak{N}$ denote the unipotent radical of the parabolic subgroup containing $\mathfrak{S}$ and opposite to $\mathfrak{P}$. Let $G = \mathfrak{G}(L)$, $S = \mathfrak{S}(L)$, $P = \mathfrak{P}(L)$, $N = \mathfrak{N}(L)$ denote the corresponding groups of $L$-rational points, considered as $L$-analytic groups. Finally, let

$$\mathfrak{g} = \mathfrak{g}_0 \oplus \bigoplus_{a \in \Phi_{\mathrm{nd}}} \mathfrak{g}_a$$

be the $L$-linear root space decomposition of the Lie algebra $\mathfrak{g}$ of $G$; here $\mathfrak{g}_0$ denotes the Lie algebra of $S$. By an abuse of language we shall call *root space vectors* the nonzero elements of $\mathfrak{g}$ which are contained in some $\mathfrak{g}_a$ ($a \in \Phi_{\mathrm{nd}}$ or $a = 0$).

We will consider admissible representations of open subgroups $G_1 \subset G$ over $K$.

### 9. Definition of $\mathfrak{n}$- characters

**Definition 9.1.** Let $\mathfrak{G}' \subset \mathfrak{G}$ be an $L$-subgroup such that the Lie algebra $\mathfrak{g}'$ of $\mathfrak{G}'$ is spanned by root space vectors. Let $G_0 \subset \mathfrak{G}'(L)$ be an open compact subgroup. We define $\Xi(G_0)$ to be the set of pairs $(H, r)$ where

(a) $H \subset G_0$ is an open normal $L$-uniform subgroup such that the $o_L$-lattice $\mathrm{Log}(H)$ is spanned by root space vectors and is contained in the neighbourhood $(\mathfrak{g}')^{\sim} \subset \mathfrak{g}'$ specified in Proposition 4.3, and

(b) $r \in p^{\mathbb{Q}} \cap \left[p^{-1}, 1\right)$.

Let us keep the notations of the definition. The groups $H \subset G_0$ with property (a) form a basis of neighbourhoods of 1 in $G_0$. Indeed, let $H' \subset G_0$ be an arbitrary open normal $L$-uniform subgroup. Then the $o_L$-lattice $\Lambda' = \mathrm{Log}(H')$ contains an $o_L$-lattice $\Lambda$ spanned by root space vectors and contained in $(\mathfrak{g}')^{\sim}$; moreover, $\Lambda \supset p^m \Lambda'$ for some $m \in \mathbb{N}$. The property of $H'$ of being uniform pro- $p$ translates to the fact that the lattice $\Lambda'$ is *powerful*, i.e. $[\Lambda', \Lambda'] \subset p^{\kappa} \Lambda'$ ([12], Section 9.4). But $[p^m \Lambda, p^m \Lambda] \subset [p^m \Lambda', p^m \Lambda'] \subset p^{\kappa+2m} \Lambda' \subset p^{\kappa} p^m \Lambda$, hence $p^m \Lambda$ is powerful and the uniform pro- $p$-group $H := \mathrm{Exp}(p^m \Lambda)$ is an open subgroup of $H'$ satisfying the above properties. This proves our claim.

For a pair $((H, r), (H', r'))$ in $\Xi(G_0)$ consider the following condition:

$$(+) \qquad \text{there exist } k, k' \in \mathbb{N} \text{ such that } \begin{aligned} & {H'}^{p^{k'}} \subset H^{p^k}, \\ & {r'}^{p^{k'}} \geq r^{p^k} \geq p^{-1}, \\ & {r'}^{\kappa p^{k'+m}} \geq p^{-1}, \end{aligned}$$

where $m$ denotes the highest $p$-elementary divisor of ${H'}^{p^{k'}}$ in $H^{p^k}$. We define a binary relation on $\Xi(G_0)$ by declaring $(H, r) \leq (H', r')$ if and only if there is a finite sequence

$$(H, r) = (H_0, r_0), (H_1, r_1), \ldots, (H_j, r_j) = (H', r') \quad (j \geq 0)$$

in $\Xi(G_0)$ such that $((H_{i-1}, r_{i-1}), (H_i, r_i))$ satisfies (+) for all $1 \leq i \leq j$. This relation is reflexive and transitive by construction, and in fact is a directed preordering: indeed, given two elements $(H, r)$, $(H', r')$ of $\Xi(G_0)$ we find a $k'$ with ${H'}^{p^{k'}} \subset H$; letting $m$ denote the highest $p$- elementary divisor of ${H'}^{p^{k'}}$ in $H$ and choosing $s < 1$ such that $s \geq \max\left\{r^{p^{-k'}}, p^{-\left(\kappa p^{k+m}\right)^{-1}}, r'\right\}$ we have $(H, r) \leq (H', s)$ and $(H', r') \leq (H', s)$.

The following result exhibits the crucial properties of the $\leq$- relation:

**Lemma 9.2.** Let $(H, r), (H', r') \in \Xi(G_0)$ such that $(H, r) \leq (H', r')$. Then:

(i) The identity map $D(G_0, K) \to D(G_0, K)$ extends to a continuous $K$-algebra homomorphism

$$D_{(H',r')}(G_0, K) \to D_{(H,r)}(G_0, K).$$

(ii) ${H'}^{\epsilon(r')} \subset H^{\epsilon(r)}$.

*Proof*. For both (i) and (ii) we may assume that (+) already holds for the pair $(H, r)$, $(H', r')$. Proposition 7.5. yields maps $D_{(H',r')}(G_0, K) \to D_{\left({H'}^{p^{k'}}, {r'}^{p^{k'}}\right)}(G_0, K) \to D_{\left(H^{p^k}, {r'}^{p^{k'}}\right)}(G_0, K) \to D_{\left(H^{p^k}, r^{p^k}\right)}(G_0, K) \to D_{(H,r)}(G_0, K)$. This proves (i). By Remark 7.3 b), $H^{\epsilon(r)} = \left(H^{p^k}\right)^{\epsilon\left(r^{p^k}\right)}$ and ${H'}^{\epsilon(r')} = \left({H'}^{p^{k'}}\right)^{\epsilon\left({r'}^{p^{k'}}\right)}$. Now (ii) follows from (+) and the monotony of the map $r \mapsto \epsilon(r)$. □

Let $G_1$ be an open subgroup of $G$, let $S_0$ be a compact subgroup of $G_1 \bigcap S$, and let $V$ be an admissible locally analytic $G_1$-representation over $K$, with strong dual $M = V_b'$. From now we take $G_0$ to be a compact open subgroup of $G_1$ containing $S_0$, and we simply write $\Xi(G_0) = \Xi$. (Such a group $G_0$ always exists, cf. [35], 3.2.) For each pair $(H, r) \in \Xi$ we define the left $D_{(H,r)}(G_0, K)$- module

$$M_{(H,r)} := D_{(H,r)}(G_0, K) \underset{D(G_0,K)}{\otimes} M;$$

here $D(G_0, K)$ acts on $M$ through the inclusion $G_0 \to G_1$. Note that since $V$ is admissible $M_{(H,r)}$ is a finitely generated module over the noetherian $K$-Banach algebra $D_{(H,r)}(G_0, K)$ and hence carries a canonical $K$-Banach space topology ([32], Proposition 2.1). Using Lemma 9.2 (i) we obtain a projective system

$$(M_{(H,r)})_{(H,r)\in\Xi}.$$

**Remark 9.3.** $M = \varprojlim_{(H,r)\in\Xi} M_{(H,r)}$.

*Proof.* Since $V$ is admissible $M$ is a coadmissible module over the Fréchet-Stein algebra $(D(H', K), (\| \|_{\bar{r}})_{r<1})$ for any fixed choice of a uniform open normal subgroup $H' \subset G_0$ (viewed as an $L$-analytic group). Regarding [32], Corollary 3.1., this means that $M = \varprojlim_{r<1} M_{(H',r)}$. Hence it suffices to show that the system $(M_{(H',r)})_r$ is cofinal in $(M_{(H,r)})_{(H,r)}$. Let $(H, r) \in \Xi$. There exists a $k' \in \mathbb{N}$ with $H'^{p^{k'}} \subset H$. Let $m$ denote the highest $p$-elementary divisor of $H'^{p^{k'}}$ in $H$. Choose $r' < 1$ such that $r' \geq \max\left\{r^{p^{-k'}}, p^{-\left(\kappa p^{k'+m}\right)^{-1}}\right\}$. Then $(H, r) \leq (H', r')$. □

Fix a pair $(H, r) \in \Xi$. Let $\mathfrak{n}$ denote the Lie algebra of $N$. Since $G_0 \subset G$ is open the universal enveloping algebra $U(\mathfrak{n})$ is contained in $D_{(H,r)}(G_0, K)$ and thereby acts on $M_{(H,r)}$. Moreover, since $\mathfrak{S}$ normalizes $\mathfrak{N}$, the subspace $\mathfrak{n}^k M_{(H,r)}$ is $S_0$-invariant for any $k \in \mathbb{N}$, and we will view the quotients $M_{(H,r)} / \mathfrak{n}^k M_{(H,r)}$ as representations of $S_0$.

Let $X(S_0)$ denote the group of characters $S_0 \to C$. Recall that by *weak convergence* of a net of $C$-valued functions we always mean pointwise convergence with respect to the discrete topology on $C$.

**Definition 9.4.** Let $S' \subset S_0$ be a subset. Suppose there is a pair $(H_0, r_0) \in \Xi$ and a number $k_0 \in \mathbb{N}$ such that the following holds:

(i) for every pair $(H, r) \geq (H_0, r_0)$ the representations $M_{(H,r)} / \mathfrak{n}^k M_{(H,r)}$ $(k \geq k_0)$ are finitely trigonalisable over $C$ and the sequence $\mathrm{Ch}\left(M_{(H,r)} / \mathfrak{n}^k M_{(H,r)}\right)_{k\geq k_0}$ converges to an element $\Theta_{(H,r)}$ in $\mathbb{Z}[\![X(S_0)]\!]$;

(ii) $\Theta_{(H,r)}$ is evaluable on $S'$ for every $(H, r) \geq (H_0, r_0)$;

(iii) the net of functions $(\mathrm{ev}_{S'}(\Theta_{(H,r)}))_{(H,r)\geq(H_0,r_0)}$ converges weakly to a function $\theta : S' \to C$.

Then we say that the $G_1$-representation $V$ possesses an $\mathfrak{n}$- *character on* $S'$, and the function $\theta_V : S' \to C$ defined by $\theta_V(s) = \theta\left(s^{-1}\right)$ is called the $\mathfrak{n}$- *character* (or simply the *character*) *on* $S'$ of $V$.

**Remark 9.4.** The character $\theta_V$ (existence and value) depends only on the $G_1$-representation $V$, the groups $N$ and $S_0$, and the subset $S'$, while $G_0$ should be considered as an auxiliary group which may be replaced by any other compact open subgroup of $G_1$ containing $S_0$. Indeed, let $H_0 \subset G_0$ an open subgroup containing $S_0$. Then $\Xi(H_0)$ is a preordered subset of $\Xi(G_0)$ which is cofinal: if $(H, r) \in \Xi(G_0)$ then $\left(H^{p^k}, r\right) \in \Xi(H_0)$ for sufficiently large $k \in \mathbb{N}$, and $\left(H^{p^k}, r\right) \geq (H, r)$. Let $(H, r) \in \Xi(H_0)$. Then

$$D_{(H,r)}(H_0, K) \underset{D(H_0,K)}{\otimes} D(G_0, K) = D_{(H,r)}(H_0, K) \underset{D(H_0,K)}{\otimes} \left( \bigoplus_{x \in H_0 \backslash G_0} D(H_0, K)\, \delta_x \right)$$

$$= \bigoplus_{x \in H_0 \backslash G_0} \left( D_{(H,r)}(H_0, K) \underset{D(H_0,K)}{\otimes} D(H_0, K)\, \delta_x \right) = D_{(H,r)}(G_0, K),$$

hence $D_{(H,r)}(H_0, K) \otimes_{D(H_0,K)} M = D_{(H,r)}(H_0, K) \otimes_{D(H_0,K)} D(G_0, K) \otimes_{D(G_0,K)} M = D_{(H,r)}(G_0, K) \otimes_{D(G_0,K)} M$ as $S_0$-representations.

**Proposition 9.5.** Let $S' \subset S_0$ be a subset. Suppose that $V = \bigoplus_{i \in I} V_i$ is the finite direct sum of $G_1$-stable subspaces $V_i$. If each $V_i$ possesses an $\mathfrak{n}$-character $\theta_{V_i}$ on $S'$ then $V$ possesses an $\mathfrak{n}$-character $\theta_V$ on $S'$, and

$$\theta_V = \sum_{i \in I} \theta_{V_i}.$$

*Proof.* For the strong duals we have $M = \oplus_{i \in I} M^{(i)}$ where $M^{(i)} := (V_i)'_b$. Furthermore, $M_{(H,r)} = D_{(H,r)}(G_0, K) \otimes_{D(G_0,K)} \left(\oplus_{i \in I} M^{(i)}\right) = \oplus_{i \in I}\left(D_{(H,r)}(G_0, K) \otimes_{D(G_0,K)} M^{(i)}\right) = \oplus_{i \in I} M^{(i)}_{(H,r)}$ for all $(H, r) \in \Xi$, and finally $M_{(H,r)} \big/ \mathfrak{n}^k M_{(H,r)} = \oplus_{i \in I}\left(M^{(i)}_{(H,r)} \big/ \mathfrak{n}^k M^{(i)}_{(H,r)}\right)$ for all $k$ because the direct sum is $\mathfrak{n}$-stable. By assumption we find a pair $(H_0, r_0)$ such that for all $(H, r) \geq (H_0, r_0)$ and *for every* $i$ the formal characters $\mathrm{Ch}\left(M^{(i)}_{(H,r)} \big/ \mathfrak{n}^k M^{(i)}_{(H,r)}\right)$ exist and converge to an element

$$\Theta^{(i)}_{(H,r)} := \lim_{k \to \infty} \mathrm{Ch}\left(M^{(i)}_{(H,r)} \big/ \mathfrak{n}^k M^{(i)}_{(H,r)}\right).$$

Hence (Lemma 2.3) each $M_{(H,r)} \big/ \mathfrak{n}^k M_{(H,r)}$ possesses a formal character satisfying

$$\mathrm{Ch}\Big(M_{(H,r)} \Big/ \mathfrak{n}^k\, M_{(H,r)}\Big) = \sum_{i\in I} \mathrm{Ch}\Big(M^{(i)}_{(H,r)} \Big/ \mathfrak{n}^k\, M^{(i)}_{(H,r)}\Big).$$

Therefore

$$\begin{aligned} \sum_{i\in I} \Theta^{(i)}_{(H,r)} &= \sum_{i\in I} \lim_{k\to\infty} \mathrm{Ch}\left(M^{(i)}_{(H,r)} \Big/ \mathfrak{n}^k\, M^{(i)}_{(H,r)}\right) \\ &= \lim_{k\to\infty} \sum_{i\in I} \mathrm{Ch}\left(M^{(i)}_{(H,r)} \Big/ \mathfrak{n}^k\, M^{(i)}_{(H,r)}\right) \\ &= \lim_{k\to\infty} \mathrm{Ch}\left(M_{(H,r)} \Big/ \mathfrak{n}^k\, M_{(H,r)}\right) \quad = \quad : \Theta_{(H,r)}\,. \end{aligned}$$

For each $i$, by assumption, the elements $\Theta^{(i)}_{(H,r)}$ are evaluable on $S'$ and the evaluations converge weakly to a function $\theta^{(i)} : S' \to C$. Therefore, since $\mathbb{Z}[\![X(S_0)]\!]_{S'}$ is a $\mathbb{Z}$-module and $\mathrm{ev}_{S'}$ is $\mathbb{Z}$-linear, the elements $\Theta_{(H,r)}$ are evaluable on $S'$ and

$$\mathrm{ev}_{S'}(\Theta_{(H,r)})\,(s) = \sum_{i\in I} \mathrm{ev}_{S'}\Big(\Theta^{(i)}_{(H,r)}\Big)\,(s).$$

Hence the functions $\mathrm{ev}_{S'}(\Theta_{(H,r)})$ converge weakly to $\sum_{i\in I} \theta^{(i)}$. □

## 10. The smooth case

In this section we intend to show that, modulo restriction to $S_0$, Definition 9.4 generalizes the concept of character of the smooth representation theory.

Let $V$ be a smooth representation of an open subgroup $G_1 \subset G$ over $K$. Let $\mathcal{H}(G_1)$ be the Hecke algebra of $G_1$, i.e. the $K$-vector space of compactly supported locally constant $K$- valued functions on $G_1$, endowed with the convolution product. Choose a Haar measure $\mu$ on $G_1$. Then $\mathcal{H}(G_1)$ acts on $V$ by the rule

$$(f, v) \mapsto \int_{G_1} f \cdot v^{\bullet}\, d\mu$$

($f \in \mathcal{H}(G_1)$, $v \in V$; here $v^{\bullet} : G_1 \to V$ denotes the locally constant orbit map $g \mapsto g\,v$). The smooth representation $V$ is called *admissible* if every $f \in \mathcal{H}(G_1)$ has finite rank as an operator on $V$, equivalently: if the subspace of $H$-invariants $V^H \subset V$ is finite-dimensional for every compact open subgroup $H \subset G_1$ (cf. [10], p. 119). If this is the case then the linear form

$$\Theta_V : \mathcal{H}(G_1) \to K,\ f \mapsto \operatorname{tr}\left( V \to V,\ v \mapsto \int_{G_1} f \cdot v^{\bullet}\, d\mu \right)$$

is defined.

**Definition 10.1.** Let $G' \subset G_1$ be an open subset. Let $V$ be an admissible smooth representation of $G_1$. Let $\theta : G' \to C$ be a locally constant function such that

$$\Theta_V(f) = \int_{G_1} f \cdot \theta\, d\mu$$

for all functions $f \in \mathcal{H}(G_1)$ whose support is contained in $G'$. Then we say that $V$ *possesses a character on* $G'$, and the function $\theta$ is called the *character on* $G'$ of $V$.

For finitely generated smooth representations there is the following classical result (cf. [19], Corollary of Theorem 2; [34], Corollary 4.8.2):

**Theorem 10.2 (R. Howe, Harish-Chandra).** Let $G^{\mathrm{reg}} \subset G$ be the subset of regular elements. Let $V$ be an admissible smooth representation of $G$ such that $V$ is finitely generated as a $G$-module. Then $V$ possesses a character on $G^{\mathrm{reg}}$.

Now let $G_0 \subset G_1$ be a compact open subgroup, and consider the quotient algebra

$$D^{\infty}(G_0, K) := D(G_0, K)/(\mathfrak{g})$$

where $(\mathfrak{g}) \subset D(G_0, K)$ is the closed two-sided ideal generated by the Lie algebra $\mathfrak{g}$ of $G$. This quotient algebra may be viewed as the strong dual of the subspace $\mathcal{H}(G_0) = C^{\infty}(G_0, K) \subset C^{\mathrm{an}}(G_0, K)$ of locally constant functions (cf. [30], Section 2). As a topological $K$-vector space $D^{\infty}(G_0, K)$ is generated by the distributions $\delta_x\, \mu_H$ where $\delta_x$ denotes the Dirac distribution of an element $x \in G_0$ and $\mu_H$ denotes the normalized Haar measure of an open subgroup $H \subset G_0$ ([2], Section 2.1; compare also [8], Ex. 22 of IV.2). Note that for any open normal subgroup $H \subset G_0$ the finite-dimensional group ring $K[G_0/H]$ may be viewed as the subring of $D^{\infty}(G_0, K)$ with $K$-basis $(\delta_x\, \mu_H)_{x \in G_0/H}$.

For a pair $(H, r) \in \Xi = \Xi(G_0)$ we put

$$D^{\infty}_{(H,r)}(G_0, K) := D_{(H,r)}(G_0, K)/(\mathfrak{g})$$

where this time $(\mathfrak{g})$ is the closed two-sided ideal in $D_{(H,r)}(G_0, K)$ generated by $\mathfrak{g}$. Recall (Remark 7.3) that $H^{\epsilon(r)}$ is an open normal subgroup of $G_0$ depending on $(H, r)$.

**Lemma 10.3.** Let $(H, r) \in \Xi$. We have a canonical isomorphism of $K$- algebras

$$D^{\infty}_{(H,r)}(G_0, K) \simeq K[G_0 / H^{\epsilon(r)}].$$

*Proof.* By Corollary 7.4 we have $D_{(H,r)}(G_0, K) = \oplus_{x \in \mathcal{R}} U_r(\mathfrak{g}, K)\,\delta_x$ where $\mathcal{R}$ is a system of representatives of $G_0 / H^{\epsilon(r)}$. We claim that the closed ideal in $D_{(H,r)}(G_0, K)$ generated by $\mathfrak{g}$ is equal to

$$I := \bigoplus_{x \in \mathcal{R}} \mathfrak{g}\, U_r(\mathfrak{g}, K)\,\delta_x.$$

It is clear that $I$ is contained in that ideal and contains $\mathfrak{g}$ because $\mathfrak{g}\, U_r(\mathfrak{g}, K)$ is the closed ideal in $U_r(\mathfrak{g}, K)$ generated by $\mathfrak{g}$. Now for any $\mathfrak{x}, \mathfrak{y} \in U(\mathfrak{g})$, $x, y \in G_0$ we have

$$\mathfrak{x}\,\delta_x\,\mathfrak{y}\,\delta_y = \mathfrak{z}\,\delta_{xy} = \mathfrak{z}\,\delta_h\,\delta_z$$

with $\mathfrak{z} \in U(\mathfrak{g})$, $z \in \mathcal{R}$ and $h \in H^{\epsilon(r)}$; moreover, if either $\mathfrak{x}$ or $\mathfrak{y}$ are non-constant then $\mathfrak{z}$ is non-constant (Lemma 8.5 (i)). Since $\delta_h \in U_r(\mathfrak{g}, K)$ (Remark 7.3) this shows that $I$ is an ideal, thereby proving our claim.

It follows that $D^{\infty}_{(H,r)}(G_0, K) = D_{(H,r)}(G_0, K)/I$ as a $K$-vector space is isomorphic to $\oplus_{x \in \mathcal{R}} K\,\delta_x = K[G_0 / H^{\epsilon(r)}]$ under the map $\delta_x + I \mapsto \delta_x\,\mu_{H^{\epsilon(r)}}$, and this map is obviously multiplicative. □

Let $V$ be an admissible smooth representation of $G_1$. Following [30], Section 2, we may view $V$ as a locally analytic representation by endowing the underlying $K$-vector space with the finest locally convex topology; among the locally analytic representations such a representation is characterized by the fact that $V$ is admissible (in the sense of locally analytic representation theory) with trivial Lie algebra action ([32], Theorem 6.6). The dual $M = V'_b$ is in this case the full linear dual $V^*$, and the action of $D(G_0, K)$ on $M$ factors through $D^{\infty}(G_0, K)$. Using the above lemma we can describe the coherent sheaf structure of the coadmissible $D(G_0, K)$-module $M$ in terms of the pairs $(H, r) \in \Xi$: We have

$$M_{(H,r)} = D_{(H,r)}(G_0, K) \underset{D(G_0,K)}{\otimes} V^* = D^{\infty}_{(H,r)}(G_0, K) \underset{D^{\infty}(G_0,K)}{\otimes} V^*$$

as $D^{\infty}(G_0, K)$-modules; hence

$$M_{(H,r)} = K[G_0 / H^{\epsilon(r)}] \underset{D^{\infty}(G_0,K)}{\otimes} V^* = \left(V^{H^{\epsilon(r)}}\right)^*$$

(compare also the proof of [32], Theorem 6.6). Note that the space of $H^{\epsilon(r)}$-invariants $V^{H^{\epsilon(r)}}$ is $G_0$-stable since $H^{\epsilon(r)} \subset G_0$ is normal.

The restriction of $V^{H^{\epsilon(r)}}$ to the commutative group $S_0$ is finitely trigonalisable over $C$ and possesses a character on $S_0$, namely the restriction to $S_0$ of the usual character

$$\theta_{(H,r)} : G_0 \to K,\ g \mapsto \mathrm{tr}\big(g_\bullet : V^{H^{\epsilon(r)}} \to V^{H^{\epsilon(r)}}\big)$$

of $V^{H^{\epsilon(r)}}$ (cf. Remark 2.2 (ii)). The function $\theta_{(H,r)}$ is constant on the cosets of $H^{\epsilon(r)}$, and its connection to the linear form $\Theta_V$ is given by

$$(*)\qquad \Theta_V\left(\mathbb{1}_{g\,H^{\epsilon(r)}}\right) = \int_{G_1} \mathbb{1}_{g\,H^{\epsilon(r)}}\, \theta_{(H,r)}\, d\mu \;=\; \mu(H^{\epsilon(r)})\, \theta_{(H,r)}(g)$$

for every $g \in G_0$ (compare [10], Section 1.5, equation (12)).

**Lemma 10.4.** Let $G' \subset G_1$ be an open subset, $S' := S_0 \cap G'$.

(i) If the smooth representation $V$ possesses a character $\theta$ on $G'$ then the net of functions $(\theta_{(H,r)}\,|_{S'})_{(H,r)\in\Xi}$ converges weakly to the restriction $\theta\,|_{S'}$.

(ii) If the net of functions $(\theta_{(H,r)}\,|_{S'})_{(H,r)\in\Xi}$ converges weakly to a function $\theta : S' \to K$ which extends to a locally constant $K$-valued function $G' \to K$ then there is an open subset $G'' \subset G'$ containing $S'$ such that the smooth representation $V$ possesses a character on $G''$ which coincides with $\theta$ on $S'$.

*Proof.* (i) Let $g \in S'$. Since $G' \subset G_1$ is open and regarding Lemma 9.2 we find a pair $(H_0, r_0)$ such that $g\,H^{\epsilon(r)} \subset G'$ and such that the function $\theta$ is constant on $g\,H^{\epsilon(r)}$ for all $(H, r) \geq (H_0\, r_0)$; these two conditions imply

$$\theta(g) = \mu(H^{\epsilon(r)})^{-1} \int_{G_1} \mathbb{1}_{g\,H^{\epsilon(r)}}\, \theta\, d\mu = \mu(H^{\epsilon(r)})^{-1}\, \Theta_V\big(\mathbb{1}_{g\,H^{\epsilon(r)}}\big)$$

and therefore, using the equality $(*)$, that the value $\theta_{(H,r)}(g)$ is constant for all $(H, r) \geq (H_0\, r_0)$.

(ii) We denote the locally constant extension of $\theta$ to $G'$ again by $\theta$, and for each $g \in S'$ we choose a pair $(H_g, r_g)$ such that $\theta_{(H,r)}(g) = \theta(g)$ for all $(H, r) \geq (H_g, r_g)$. By eventually increasing $(H_g, r_g)$ we attain that $\theta$ is constant on $g\,H_g^{\epsilon(r_g)} \subset G'$ (Lemma 9.2) and hence, using $(*)$ again,

$$(**)\qquad \int_{G_1} \mathbb{1}_{g\,H^{\epsilon(r)}}\, \theta\, d\mu = \mu(H^{\epsilon(r)})\, \theta(g) = \Theta_V\big(\mathbb{1}_{g\,H^{\epsilon(r)}}\big)$$

for all $(H, r) \geq (H_g, r_g)$. Let $G''$ be the union of all $g\, H_g{}^{\epsilon(r_g)}$ $(g \in S')$. Then $(**)$ shows that $\Theta_V(f) = \int_{G_1} f\, \theta\, d\mu$ for all $f$ with $\mathrm{supp}(f) \subset G''$, i.e. the restriction of $\theta$ to $G''$ is a character of $V$ on $G''$. □

**Theorem 10.5.** Let $G_1 \subset G$ be an open subgroup. Let $V$ be the admissible locally analytic representation associated to an admissible smooth representation of $G_1$ over $K$ according to the above method. Let $S' \subset S_0$ be an open subset. The following assertions are equivalent:

(i) $V$ possesses an $\mathfrak{n}$-character $\theta$ on $S'$, and there is an open subset $G' \subset G_1$ containing $S'$ such that $\theta$ extends to a locally constant function $\theta' : G' \to K$;

(ii) there is an open subset $G'' \subset G_1$ containing $S'$ such that the underlying admissible smooth representation of $V$ possesses a character $\theta''$ on $G''$.

If these assertions hold then $\theta'$ and $\theta''$ agree on $S'$.

*Proof.* Since $\mathfrak{n}$ acts trivially on $V$ we have $M_{(H,r)} \big/ \mathfrak{n}^k M_{(H,r)} = M_{(H,r)}$ for all $k \in \mathbb{N}$. The theorem then follows from the preceding lemma. □

**Remark 10.6.** Lemma 10.3 indicates in particular that we can characterize the subspace $C^\infty_{H^{\epsilon(r)}}(G_0, K) \subset C^\infty(G_0, K)$ of functions which are constant on the cosets of $H^{\epsilon(r)}$ by means of the $\|\ \|_{(H,r)}$-norm: We have

$$C^\infty_{H^{\epsilon(r)}}(G_0, K) \simeq K[G_0 / H^{\epsilon(r)}]^* = D^\infty_{(H,r)}(G_0, K)^*.$$

This phenomenon is very peculiar to the smooth situation as we are going to illustrate in the case $G_0 = \mathbb{Z}_p$. Fix $h \in \mathbb{N}$, and let $\mathcal{O}_h(\mathbb{Z}_p, \mathbb{Q}_p) \subset C^{\mathrm{an}}(\mathbb{Z}_p, \mathbb{Q}_p)$ denote the subspace of functions which are holomorphic on the cosets of $p^h\, \mathbb{Z}_p$. On the other hand, for $r \in (0, 1]$ (resp. $(0, 1)$) we consider the space $C_r(\mathbb{Z}_p, \mathbb{Q}_p)$ (resp. $C_{r+}(\mathbb{Z}_p, \mathbb{Q}_p)$) of all continuous functions $f : \mathbb{Z}_p \to \mathbb{Q}_p$ with Mahler expansion $f = \sum_{n \in \mathbb{N}} b_n \binom{\bullet}{n}$ satisfying $\lim_{n\to\infty} |b_n|\, r^{-n} = 0$ (resp. $(|b_n|\, r^{-n})_{n\in\mathbb{N}}$ is bounded). Put $\lambda_h := p^{-h}(p-1)^{-1}$. If $0 < r \leq p^{-\lambda_h} < s \leq 1$ then

$$C_{r+}(\mathbb{Z}_p, \mathbb{Q}_p) \underset{\neq}{\subset} \mathcal{O}_h(\mathbb{Z}_p, \mathbb{Q}_p) \underset{\neq}{\subset} C_s(\mathbb{Z}_p, \mathbb{Q}_p).$$

*Proof.* First observe that for any $\lambda \in \mathbb{Q}$ and $t := p^{-\lambda}$ we have

$$f \in C_{t+} \iff (v_p(b_n) - n\,\lambda)_{n\in\mathbb{N}} \text{ bounded below.}$$

Next, by a criterion of Y. Amice ([21], 1.3.8) $f$ is holomorphic on the cosets of $p^h\, \mathbb{Z}_p$ if and only if $\lim_{n\to\infty} v_p(b_n) - \frac{p^{-h}\, n - \mathrm{ch}(n)}{p-1} = +\infty$; here $\mathrm{ch}(n) = \sum a_i$ denotes the sum of $p$-adic digits of $n = \sum a_i\, p^i$ $(0 \le a_i < p)$. In other words:

$$f \in \mathcal{O}_h(\mathbb{Z}_p, \mathbb{Q}_p) \iff \lim v_p(b_n) - n\,\lambda_h + p^h\,\lambda_h\, \mathrm{ch}(n) = +\infty.$$

"$C_{r+}(\mathbb{Z}_p, \mathbb{Q}_p) \subset \mathcal{O}_h(\mathbb{Z}_p, \mathbb{Q}_p)$": follows from $\mathrm{ch}(n) \to +\infty$.

"$\mathcal{O}_h(\mathbb{Z}_p, \mathbb{Q}_p) \subset C_s(\mathbb{Z}_p, \mathbb{Q}_p)$": Choose $\lambda = \lambda_h - \delta < \lambda_h$ $(\delta > 0)$ such that $s \ge p^{-\lambda} > p^{-\lambda_h}$. Then $v_p(b_n) - n\,\lambda = v_p(b_n) - n\,\lambda_h + n\,\delta > v_p(b_n) - n\,\lambda_h + p^h\,\lambda_h\, \mathrm{ch}(n)$ for $n >> 0$, since $\mathrm{ch}(n)/n \to 0$.

"$C_{r+}(\mathbb{Z}_p, \mathbb{Q}_p) \ne \mathcal{O}_h(\mathbb{Z}_p, \mathbb{Q}_p)$": Choose $(b_n)$ such that $v_p(b_n) = n\,\lambda_n - p^h\,\lambda_h\, \mathrm{ch}(n)/2$. Then $v_p(b_n) - n\,\lambda_h + p^h\,\lambda_h\, \mathrm{ch}(n) = p^h\,\lambda_h\, \mathrm{ch}(n)/2 \to +\infty$, hence $f \in \mathcal{O}_h(\mathbb{Z}_p, \mathbb{Q}_p)$. On the other hand $v_p(b_n) - n\,\lambda_h = -p^h\,\lambda_h\, \mathrm{ch}(n)/2$ is not bounded below.

"$\mathcal{O}_h(\mathbb{Z}_p, \mathbb{Q}_p) \ne C_s(\mathbb{Z}_p, \mathbb{Q}_p)$": Choose $(b_n)$ such that $v_p(b_n) = n\,\lambda_n - p^h\,\lambda_h\, \mathrm{ch}(n)$. Then $v_p(b) - n\,\lambda_n - p^h\,\lambda_h\, \mathrm{ch}(n)$ is constant, hence $f \notin \mathcal{O}_h(\mathbb{Z}_p, \mathbb{Q}_p)$. It suffices to show that $v_p(b_n) - n\,\lambda \to +\infty$ for all $\lambda < \lambda_h$. Put $\varepsilon = \lambda - \lambda_h$. Choose $N \in \mathbb{N}$ such that $p^h\,\lambda_h\, \mathrm{ch}(n)/n \le \varepsilon/2$ and $n\,\varepsilon/2 \ge M$ for all $n \ge N$. Then $v_p(b_n) - n\,\lambda = n\left(\varepsilon - p^h\,\lambda_h\, \mathrm{ch}(n)/n\right) \ge M$ for all $n \ge N$. □

## 11. Preliminaries on principal series representations

Let $\mathfrak{M} \subset \mathfrak{P}$ be the Levi subgroup containing $\mathfrak{S}$. Let $\chi$ be a $K$-linear locally analytic $P$-representation which is the composite of a locally analytic character $\mathfrak{M}(L) \to K^\times$ and the canonical map $P \to \mathfrak{M}(L)$. The *locally analytic principal series representation of G induced from* $\chi$ is defined to be the locally convex $K$-vector space

$$V = \mathrm{Ind}_P^G(\chi) = \left\{F \in C^{\mathrm{an}}(G, K);\, F(g\,p) = \chi(p)^{-1}\, F(g)\right\}$$

together with the left translation action of the group $G$. Below we will see that $V$ is indeed admissible.

Using the Bruhat-Tits decomposition the locally convex vector space $V$ may be identified with a finite direct sum of function spaces on compact unipotent groups. For the purpose of explicit calculations we will give the identification in extenso.

Let $G_0$ be the Iwahori subgroup of $G$ of the same type as $P$. Then $G_0$ is open and is contained in a special maximal compact subgroup $G_1 \subset G$. The Iwasawa decomposition $G = G_1\,P$ ([10], Section 3.5) induces an isomorphism of $G_1$- representations

$$
\begin{array}{rcl}
V_0 := \mathrm{Ind}_{P\bigcap G_1}^{G_1}(\chi\,|_{G_1}) & \longrightarrow & V \\
F & \longmapsto & \big(g\,p \mapsto \chi^{-1}(p)\,F(g)\big) \qquad (g \in G_1,\ p \in P) \\
F\,|_{G_1} & \longleftarrow & F
\end{array}
$$

For $w \in W$ put $P_w := w\,P\,w^{-1} \bigcap G_0$, $N_w := w\,N\,w^{-1} \bigcap G_0$. Then $G_1 = \bigcup_{w\in W} G_0\,w\,P$ (disjoint union; cf. [10], section 3.5). On the other hand, for each $w \in W$ the map

$$
\begin{array}{rcl}
V_0 \bigcap C^{\mathrm{an}}(G_0\,w\,P,\,E) & \longrightarrow & \mathrm{Ind}_{P_w}^{G_0}(\chi_w) \quad =: \ V_w \\
F & \longmapsto & F \circ R_w
\end{array}
$$

is a well-defined isomorphism of $G_0$-representations; here $R_w : G_0 \to G_0\,w\,P$ denotes right multiplication with some fixed representative of $w$ in $N_G(S)$, and $\chi_w : P_w \to E$ is defined by composing $\chi$ with conjugation by $w^{-1}$. Furthermore, the product map $N_w \times P_w \to G_0$ is an isomorphism of $L$-analytic manifolds (cf. [23], Lemma 3.3.2). Let $\pi_w^- : G_0 \to N_w$, $\pi_w^+ : G_0 \to P_w$ denote the components of the inverse map. An element $F \in V_w$ is determined by its restriction to $N_w$ which in turn can be any function, whence the restriction isomorphism $V_w \to C^{\mathrm{an}}(N_w, K)$. The inverse map $C^{\mathrm{an}}(N_w, K) \to V_w$ is explicitely given by $f \mapsto \big(\chi_w^{-1} \circ \pi_w^+\big) \cdot (f \circ \pi_w^-)$. Identifying the $G_0$-representation $V_w$ with $C^{\mathrm{an}}(N_w, K)$ we finally obtain the direct sum decomposition of $G_0$- representations

$$
\begin{array}{rcl}
V & \longrightarrow & \bigoplus_{w\in W} V_w \\
F & \longmapsto & \big(F \circ R_w\,|_{N_w}\big)_{w\in W} \quad .
\end{array}
$$

The preimage $F \in V$ of an element $(f_w)_{w\in W} \in \bigoplus_{w\in W} V_w$ is given by

$$
F(b\,w\,p) = \chi\big(w^{-1}\,\pi_w^+(b)\,w\,p\big)^{-1}\ f_w(\pi_w^-(b)) \qquad (b \in G_0,\ w \in W,\ p \in P).
$$

Note that since $G_0 \subset G$ is open we also have the Lie algebra action of $\mathfrak{g}$ on each $V_w$. We let $\mathfrak{n}_w$ denote the Lie algebra of $N_w$, and we put $S_0 = S \bigcap G_0$. As an easy consequence of the explicit description of the above decomposition we obtain the following:

**Lemma 11.1.** Fix a Weyl group element $w \in W$.

(i) The action of $G_0$ on $V_w$ is given by

$$
(x\,f)\,(y) = \chi_w^{-1}\big(\pi_w^+\big(x^{-1}\,y\big)\big)\ f\big(\pi_w^-\big(x^{-1}\,y\big)\big)
$$

($x \in G_0$, $f \in V_w$, $y \in N_w$). Special cases:

$$
(x\,f)\,(y) = \chi_w(x)\ f\big(x^{-1}\,y\,x\big) \quad \text{if } x \in S_0,
$$

$$(x\,f)\,(y) = f\big(x^{-1}\,y\big) \quad \text{if } x \in N_w.$$

(ii) The Lie algebra action on $V_w$ is given by

$$\begin{aligned}(\mathfrak{x}\,f)\,(y) &= f(y)\cdot \tfrac{d}{d\,t}\,\chi_w^{-1}(\pi_w^+(\exp(-t\,\mathfrak{x})\,y))\,|_{t=0} \\ &+ \tfrac{d}{d\,t}\,f(\pi_w^-(\exp(-t\,\mathfrak{x})\,y))\,|_{t=0}\end{aligned}$$

($\mathfrak{x} \in \mathfrak{g}$, $f \in V_w$, $y \in N_w$). Special case:

$$(\mathfrak{x}\,f)\,(y) = \tfrac{d}{d\,t}\,f(\exp(-t\,\mathfrak{x})\,y)\,|_{t=0} \quad \text{if } x \in \mathfrak{n}_w. \ \square$$

Passing to the duals $M := V_b' = D(G,K) \otimes_{D(P,K)} \chi'$ and $M^w := (V_w)_b' = D(G_0,K) \otimes_{D(P_w,K)} \chi_w'$ ($w \in W$) we get the direct sum decomposition of $D(G_0,K)$- modules

$$M = \bigoplus_{w\in W} M^w.$$

The above decomposition $G_0 = N_w\,P_w$ induces an isomorphism of $K$- Fréchet spaces

$$\begin{aligned} D(G_0,K) &\longrightarrow D(N_w,K)\,\hat{\otimes}_K\,D(P_w,K) \\ \lambda &\longmapsto \lambda\circ(\pi_w^-)^* \otimes \lambda\circ(\pi_w^+)^* \end{aligned}$$

(cf. [23], Proposition 3.3.4). Hence the natural map

$$\begin{aligned} D(N_w,K) \otimes_K \chi_w' &\longrightarrow D(G_0,K) \otimes_{D(P_w,K)} \chi_w' = M^w \\ \lambda\otimes a &\longmapsto \lambda\otimes a \end{aligned}$$

is an isomorphism. If we identify $M^w = D(N_w,K)$ then the restriction of the $D(G_0,K)$-action on $M^w$ to $D(N_w,K)$ is simply left multiplication in $D(N_w,K)$. Hence $M^w$ is finitely generated as a $D(G_0,K)$ module; in other words, $V_w$ is strongly admissible in the sense of [29], Section 3. Consequently:

**Corollary 11.2.** $V$ is strongly admissible as a locally analytic $G_0$-representation. $\square$

This result is also proved in [33], Section 6 (for the case $L = \mathbb{Q}_p$; the argument given however works also for general $L$).

For the remainder of this section we fix a Weyl group element $w \in W$. Let $\Delta_w \subset \Phi_{\mathrm{nd}}$ be the root system basis corresponding to the minimal parabolic subgroup $w\,\mathbf{P}\,w^{-1}$ and let $\Phi_w \subset \Phi_{\mathrm{nd}}$ denote the subset of negative roots with respect to $\Delta_w$; thus $\mathfrak{n}_w = \oplus_{a\in\Phi_w}\mathfrak{g}_a$ and $N_w = \prod_{a\in\Phi_w} U_a$ (direct span in any fixed order), where $U_a$ denotes the intersection of the root subgroup $\mathfrak{U}_a(L)$ with $G_0$ ([3] 21.11). Since $w\,\mathfrak{N}\,w^{-1}$ is unipotent we have the global bijections

$$w\,N\,w^{-1} \underset{\mathrm{Exp}}{\overset{\mathrm{Log}}{\rightleftarrows}} \mathfrak{n}_w,$$

inducing bijections $N_w \rightleftarrows \mathrm{Log}(N_w) \subset \mathfrak{n}_w$, $U_a \rightleftarrows \mathrm{Log}(U_a) \subset \mathfrak{g}_a$ (cf. Example 4.2).

**Proposition 11.3.** Let $H \subset G_0$ be an open $L$-uniform subgroup such that the $o_L$-lattice $\mathrm{Log}(H) \subset \mathfrak{g}$ has a basis $\mathfrak{X}$ consisting of root space vectors and is contained in the neighbourhood $\mathfrak{g}^{\sim} \subset \mathfrak{g}$ specified in Proposition 4.3. Then $H \bigcap N_w$ and $H \bigcap P_w$ are closed $L$-uniform subgroups of $H$ and the $o_L$-lattices $\mathrm{Log}(H \bigcap N_w)$ and $\mathrm{Log}(H \bigcap P_w)$ are spanned by the elements of $\mathfrak{X}$ they contain.

*Proof.* Let $\mathfrak{p}_w$ denote the Lie algebra of $P_w$. The intersections $\mathrm{Log}(H) \bigcap \mathfrak{n}_w$ and $\mathrm{Log}(H) \bigcap \mathfrak{p}_w$ are $o_L$-Lie subalgebras of $\mathrm{Log}(H)$. Moreover, the $L$-vector spaces $\mathfrak{n}_w$ and $\mathfrak{p}_w$ are direct sums of the root spaces $\mathfrak{g}_a$ ($a \in \Phi \bigcup \{0\}$) they contain. Therefore the $o_L$-lattices $\mathrm{Log}(H) \bigcap \mathfrak{n}_w$ and $\mathrm{Log}(H) \bigcap \mathfrak{p}_w$ are spanned by the elements of $\mathfrak{X}$ they contain. This also implies that the quotients $\mathrm{Log}(H)/(\mathrm{Log}(H) \bigcap \mathfrak{n}_w)$ and $\mathrm{Log}(H)/(\mathrm{Log}(H) \bigcap \mathfrak{p}_w)$ are $\mathbb{Z}_p$-torsion free. Hence ([12], Proposition 7.15 (i)) those lattices correspond to closed $\mathbb{Q}_p$-uniform subgroups $H' := \mathrm{Exp}(\mathrm{Log}(H) \bigcap \mathfrak{n}_w)$, $H'' := \mathrm{Exp}(\mathrm{Log}(H) \bigcap \mathfrak{p}_w)$ of $H$. By construction $H, H'$ are in fact $L$-uniform.

It remains to be seen that $H' = H \bigcap N_w$, $H'' = H \bigcap P_w$. If $g \in H \bigcap N_w$ then $\mathrm{Log}(g) \in \mathrm{Log}(H) \bigcap \mathrm{Log}(N_w) \subset \mathrm{Log}(H) \bigcap \mathfrak{n}_w$, hence $g \in H'$. In the other direction, if $g \in H'$ then $\mathrm{Log}(g) \in \mathrm{Log}(H) \bigcap \mathfrak{n}_w$, and from Proposition 4.3 we obtain $g \in H' \bigcap w\,N\,w^{-1} = H' \bigcap N_w$. Thus we have proved $H' = H \bigcap N_w$. The equality $H'' = H \bigcap P_w$ is shown analogously. □

**Corollary 11.4.** Let $(H, r), (H', r') \in \Xi(G_0)$. Then $(H \bigcap N_w, r), (H' \bigcap N_w, r') \in \Xi(N_w)$ and $(H \bigcap P_w, r), (H' \bigcap P_w, r') \in \Xi(P_w)$. If $(H, r) \leq (H', r')$ then $(H \bigcap N_w, r) \leq (H' \bigcap N_w, r')$ and $(H \bigcap P_w, r) \leq (H' \bigcap P_w, r')$.

*Proof.* The first assertion follows from Proposition 11.3. For the second assertion we certainly may assume that there already exist natural numbers $k, k'$ satisfying the condition (+) from Section 9. By Proposition 11.3 there is an ordered basis $\mathbf{h}$ of $H$ (resp. $H'$) such that the elements of $\mathbf{h}$ contained in $H \bigcap N_w$ (resp. $H' \bigcap N_w$) form an ordered basis of $H \bigcap N_w$ (resp. $H' \bigcap N_w$).

We deduce $(H \bigcap N_w)^{p^k} = H^{p^k} \bigcap N_w$, $(H' \bigcap N_w)^{p^{k'}} = H'^{p^{k'}} \bigcap N_w$. Since the highest $p$-elementary divisor of $H'^{p^{k'}} \bigcap N_w$ in $H^{p^k} \bigcap N_w$ clearly is smaller than or equal to the highest $p$-elementary divisor of $H'^{p^{k'}}$ in $H^{p^k}$ it follows that the pair $((H \bigcap N_w, r), (H' \bigcap N_w, r'))$ satisfies (+) as well (using the same numbers $k, k'$). - Analogously one proves that $(H \bigcap P_w, r) \le (H' \bigcap P_w, r')$. □

If the reductive group $\mathbb{G}$ is split over $L$ then the root spaces $\mathfrak{g}_a$ are one-dimensional, and the group $N_w$ is "as good as $L$-uniform", in the following sense: Write $\Phi_w = \{a_1, \ldots, a_d\}$. Let $\mathfrak{x}_j \in \mathfrak{g}_{a_j}$ be a generator of the $o_L$-lattice $\mathrm{Log}\big(U_{a_j}\big)$. Let $(v_1, \ldots, v_n)$ be a $\mathbb{Z}_p$-basis of $o_L$, $h_{ij} := \mathrm{Exp}(v_i\, \mathfrak{x}_j)$. Then $(v_1\, \mathfrak{x}_j, \ldots, v_n\, \mathfrak{x}_j)$ is a $\mathbb{Z}_p$-basis of $\mathrm{Log}\big(U_{a_j}\big)$, and as $N_w$ is directly spanned by its subgroups $U_{a_j}$ every element $g \in N_w$ has a unique representation

$$\begin{aligned} g &= g_1 \quad \ldots \quad g_d && \big(g_j \in U_{a_j}\big) \\ &= h_{11}^{t_{11}} \ldots h_{n1}^{t_{n1}} \ \ldots \ h_{1d}^{t_{1d}} \ldots h_{nd}^{t_{nd}} && (t_{ij} \in \mathbb{Z}_p) \quad . \end{aligned}$$

**Lemma 11.5.** Assume $\mathbb{G}$ to be split over $L$ and $L \,|\, \mathbb{Q}_p$ unramified, and define $(h_{ij})$ as above. Let $H \subset N_w$ be an $L$-uniform normal open subgroup such that $\mathrm{Log}(H)$ is spanned by root space vectors. There are integers $\alpha(1, 1), \ldots, \alpha(n, d) \in \mathbb{N}$ such that the elements $h_{11}^{p^{\alpha(1,1)}}, \ldots, h_{nd}^{p^{\alpha(n,d)}}$ form an ordered basis of $H$. A system of representatives of $N_w / H$ is given by the elements $h_{11}^{t_{11}} \ldots h_{nd}^{t_{nd}}$ $(0 \le t_{ij} < p^{\alpha(i,j)})$.

*Proof.* Let $(\mathfrak{x}'_1, \ldots, \mathfrak{x}'_d)$ be an $o_L$-basis of $\mathrm{Log}(H)$ with $\mathfrak{x}'_j \in \mathfrak{g}_{a_j}$ $(1 \le j \le d)$. As $\mathrm{Log}(H) \subset \mathrm{Log}(N_w)$ is an $o_L$-sublattice we have $\mathfrak{x}'_j = c_j\, \mathfrak{x}_j$ for some $c_j \in o_L$, and since $L$ is unramified over $\mathbb{Q}_p$ we may choose the $\mathfrak{x}'_j$ such that $\mathfrak{x}'_j = p^{\alpha(j)}\, \mathfrak{x}_j$ for some $\alpha(j) \in \mathbb{N}$. The elements $\mathrm{Exp}\big(v_i\, \mathfrak{x}'_j\big) = \mathrm{Exp}(v_i\, \mathfrak{x}_j)^{p^{\alpha(j)}} = h_{ij}^{p^{\alpha(j)}}$ $(1 \le i \le n,\ 1 \le j \le d)$ constitute an ordered basis of $H$. This proves the first assertion.

Since $H \subset N_w$ is normal the family $\big(h_{11}^{t_{11}} \ldots h_{nd}^{t_{nd}}\big)_{0 \le t_{ij} < p^{\alpha(i,j)}}$ contains a system of representatives (cf. the proof of Lemma 5.7). Suppose $h_{11}^{t_{11}} \ldots h_{nd}^{t_{nd}} = h_{11}^{u_{11}} \ldots h_{nd}^{u_{nd}}\, h$ with $0 \le t_{ij}, u_{ij} < p^{\alpha(i,j)}$ and $h \in H$. We have to show that $t_{ij} = u_{ij}$.

Put $h^{(j)} := h_{1j}^{u_{1j}} \ldots h_{nj}^{u_{nj}}$ $(1 \le j \le d)$. Then $h_{11}^{u_{11}} \ldots h_{nd}^{u_{nd}}\, h = h^{(1)} \ \ldots \ h^{(d)}\ h_1 \ldots h_d$ with uniquely determinded elements $h_j \in H \bigcap U_{a_j}$, and we claim that further

(+) $$h^{(1)} \dots h^{(d)}\, h_1 \dots h_d = h^{(1)}\, h'_1 \dots h^{(d)}\, h'_d$$

with suitable $h'_j \in H \cap U_{a_j}$.

To see this, consider arbitrary elements $n_1, \dots, n_r \in N_w$ and $g \in H$, each contained in some root subgroup. If $g \neq 1$ and $U_a$ is the root subgroup containing $g$ then let $k(g)$ denote the *height* of the root $-a$, i.e. the unique natural number $k$ such that $-a = b_1 + \dots + b_k$ with simple roots $b_1, \dots, b_k$ (with respect to $\Delta_w$). If $g = 1$ we put $k(g) := \infty$. Since $H$ is normal in $N_w$ we have

(++) $$n_1 \dots n_r\, g = g\, n_1\, g_1 \dots n_r\, g_r$$

with commutators $g_j = (n_j, g) \in H$ $(1 \leq j \leq r)$; moreover, each $g_j$ is a finite product of root subgroup elements $g_{i\,j} \in H$ satisfying $k(g_{i\,j}) \geq k(g) + 1$ ([3], Proposition 14.5, in particular the assertion marked $(*)$). Similarly,

(+++) $$g\, n_1 \dots n_r = g'_1\, n_1 \dots g'_r\, n_r\, g$$

with commutators $g'_j = (g, n_j) \in H$ $(1 \leq j \leq r)$, each one a finite product of root subgroup elements $g'_{i\,j} \in H$ satisfying $k\left(g'_{i\,j}\right) \geq k(g) + 1$.

As there is a bound $k_0 \in \mathbb{N}$ such that $k(g) > k_0$ implies $g = 1$, a finite number of recursive applications of transformations of the form (++) or (+++) to the left hand side of (+) gives the right hand side of (+). Indeed, in a first step we use (++) to move each $h_j$ according to the root subgroup it belongs to and obtain

$$h^{(1)}\, h^{(2)} \dots h^{(d)}\, h_1\, h_2 \dots h_d = h^{(1)}\, h_1\, h^{(2)}\, h_2\, h'_2 \dots h^{(d)}\, h_d\, h'_d$$

where each appearing $h'_j$ is a finite product of root subgroup elements $h'_{j\,k}$ with $k\left(h'_{j\,k}\right) \geq k(h_j) + 1$. In a second step, each of these $h'_{j\,k}$ (unless equal to 1) is moved according to the root subgroup it belongs to, using either (++) or (+++). Along the way commutators appear, being products of root subgroup elements $h''_{j\,k\,l}$ with $k\left(h''_{j\,k\,l}\right) \geq k(h_j) + 2$. In a third step the $h''_{j\,k\,l}$ (unless equal to 1) are moved, and so on. After finitely many steps all commutators are equal to 1, and we are done.

From (+) we deduce $h_{1\,j}^{t_{1\,j}} \dots h_{nj}^{t_{nj}} \left(h_{1j}^{u_{1j}} \dots h_{nj}^{u_{nj}}\right)^{-1} \in H$ for all $j$ and thereby, since the groups $U_{a_j}$ are commutative, $h_{ij}^{t_{ij}}\, h_{ij}^{-u_{ij}} \in H$ and hence $u_{i\,j} = t_{i\,j}$ for all $i, j$. □

We finish this section by stating a result describing the $D(G_0, K)$- modules

$$M^w_{(H,r)} = D_{(H,r)}(G_0, K) \underset{D(G_0,K)}{\otimes} M^w = D_{(H,r)}(G_0, K) \underset{D(P_w,K)}{\otimes} \chi'_w$$

that were considered in Section 9 for the pairs $(H, r) \in \Xi(G)$.

**Proposition 11.6 (Frommer-Orlik-Strauch).** There is a pair $(H_0, r_0) \in \Xi(G_0)$ such that the natural $K$-linear map

$$(*) \qquad D_{(H \cap N_w, r)}(N_w, K) \underset{K}{\otimes} \chi'_w \to D_{(H,r)}(G_0, K) \underset{D(P_w,K)}{\otimes} \chi'_w$$

is a topological isomorphism for all $(H, r) \geq (H_0, r_0)$.

*Proof.* For fixed $H$ this is Proposition 3.4.2 of [23]; the case $L = \mathbb{Q}_p$ is already contained in [16], Proposition 7. We adapt the proof to our directed set $\Xi = \Xi(G_0)$.

Let $(H, r) \in \Xi$. Let us assume that the $D(P_w, K)$-action on $\chi'_w$ extends to a *continuous* $D_{(H \cap P_w, r)}(P_w, K)$-action. Then

$$D_{(H,r)}(G_0, K) \underset{D(P_w,K)}{\otimes} \chi'_w = D_{(H,r)}(G_0, K) \underset{D_{(H \cap P_w, r)}(P_w,K)}{\otimes} \chi'_w.$$

Moreover, the decomposition $G_0 = N_w P_w$ induces an isomorphism of $K$-Banach spaces

$$D_{(H,r)}(G_0, K) = D_{(H \cap N_w, r)}(N_w, K) \underset{K}{\hat{\otimes}} D_{(H \cap P_w, r)}(P_w, K)$$

(cf. [23], Proposition 3.3.4; note that by Proposition 11.3 our groups $H$, $H \cap N_w$, $H \cap P_w$ satisfy the same conditions as the corresponding groups $P_0$, $U^-_{w,0}$, $P^+_{w,0}$ constructed in Section 3.3.3 of that paper). Hence the map $(*)$ is bijective. The canonical topology on the finitely generated $D_{(H,r)}(G_0, K)$-module $D_{(H,r)}(G_0, K) \otimes_{D(P_w,K)} \chi'_w$ is equal to the quotient topology with respect to the obvious surjection $D_{(H,r)}(G_0, K) \to D_{(H,r)}(G_0, K) \otimes_{D(P_w,K)} \chi'_w$. Hence the map $(*)$ is a continuous bijection and consequently, by the open mapping theorem, a topological isomorphism.

It remains to show that there exists a pair $(H_0, r_0) \in \Xi$ such that for any $(H, r) \geq (H_0, r_0)$ the $D(P_w, K)$-action on $\chi'_w$ extends to a continuous $D_{(H \cap P_w, r)}(P_w, K)$-action. Using Proposition 3.4.2 (i) of [23] we find a pair $(H_0, r_0)$ and an extension to a continuous $D_{(H_0 \cap P_w, r_0)}(P_w, K)$-action. Let $(H, r) \geq (H_0, r_0)$. We then know (Corollary 11.4 and Lemma 9.2) that there is a continuous $K$-algebra homomorphism

$$D_{(H \cap P_w, r)}(P_w, K) \to D_{(H_0 \cap P_w, r_0)}(P_w, K)$$

inducing the desired continuous extension of the $D(P_w, K)$-action on $\chi'_w$. $\square$

## 12. The principal series of $\mathrm{SL}_2(\mathbb{Q}_p)$

We specialize the setting of Section 11 to the following situation: Let $G$ be the $\mathbb{Q}_p$- analytic group $\mathrm{SL}_2(\mathbb{Q}_p)$, $S \subset G$ the standard torus, $P \subset G$ the lower Borel subgroup containing $S$, and let $N \subset G$ be the subgroup of upper strictly triangular matrices, with Lie algebra $\mathfrak{n}$. Let $\chi : \mathbb{Q}_p^\times \to K^\times$ be a locally analytic character, extended to $P$ via the map $\begin{pmatrix} a & \\ * & a^{-1} \end{pmatrix} \mapsto a^{-1}$. We want to calculate the $\mathfrak{n}$-character of the locally analytic $G$- representation

$$V_\chi = \mathrm{Ind}_P^G(\chi)$$

on the subset $S_0 \cap G^{\mathrm{reg}}$ of regular elements.

As representatives for the two elements of the Weyl group $W$ of $(G, S)$ we choose $w_+ = \begin{pmatrix} 1 & \\ & 1 \end{pmatrix}$, $w_- = \begin{pmatrix} & -1 \\ 1 & \end{pmatrix}$. Adopting the previous notations we get the following subgroups of $G$: $G_0 = \begin{pmatrix} \mathbb{Z}_p^\times & p\mathbb{Z}_p \\ \mathbb{Z}_p & \mathbb{Z}_p^\times \end{pmatrix} \cap G$, $S_0 = \left\{ \begin{pmatrix} a & \\ & a^{-1} \end{pmatrix} ; a \in \mathbb{Z}_p^\times \right\}$, $P_{w_+} = \begin{pmatrix} \mathbb{Z}_p^\times & \\ \mathbb{Z}_p & \mathbb{Z}_p^\times \end{pmatrix} \cap G$, $P_{w_-} = \begin{pmatrix} \mathbb{Z}_p^\times & p\mathbb{Z}_p \\ & \mathbb{Z}_p^\times \end{pmatrix} \cap G$, $N_{w_+} = \begin{pmatrix} 1 & p\mathbb{Z}_p \\ & 1 \end{pmatrix}$, $N_{w_-} = \begin{pmatrix} 1 & \\ \mathbb{Z}_p & 1 \end{pmatrix}$. Identifying the two latter groups in the obvious way with $\mathbb{Z}_p$ the underlying locally analytic vector spaces of the components $(V_\chi)_{w_+}$, $(V_\chi)_{w_-}$ are copies of $C^{\mathrm{an}}(\mathbb{Z}_p, K)$. Put

$$\mathfrak{x}_+ := \begin{pmatrix} 0 & p \\ 0 & 0 \end{pmatrix} \in \mathfrak{n}_{w_+} = \mathfrak{n},$$
$$\mathfrak{x}_- := \begin{pmatrix} 0 & 0 \\ 1 & 0 \end{pmatrix} \in \mathfrak{n}_{w_-}.$$

For a function $f \in C^{\mathrm{an}}(U, K)$ on an open subset $U \subset \mathbb{Q}_p$ we denote by $\frac{d}{dx} f : U \to K$ the derivative of $f$. E.g., the derivative of $\chi$ satisfies

$$\frac{d}{dx} \chi(z) = c(\chi)\, z^{-1}\, \chi(z)$$

where we let $c(\chi) := \frac{d}{dx} \chi(1)$.

**Lemma 12.1.** (i) The action of $G_0$ on $(V_\chi)_{w_+}$ is given by

$$\left( \begin{pmatrix} 1 & pb \\ & 1 \end{pmatrix} f \right)(z) = f(z - b),$$

$$\left(\begin{pmatrix} 1 & \\ c & 1 \end{pmatrix} f\right)(z) = \chi(1 - p\,c\,z)^{-1}\, f\!\left(\tfrac{z}{1-p\,c\,z}\right),$$

$$\left(\begin{pmatrix} a & \\ & a^{-1} \end{pmatrix} f\right)(z) = \chi(a)^{-1}\, f\!\left(a^{-2}\, z\right).$$

The action of $G_0$ on $(V_\chi)_{w_-}$ is given by

$$\left(\begin{pmatrix} 1 & p\,b \\ & 1 \end{pmatrix} f\right)(z) = \chi\,(1 - p\,b\,z)^{-1}\, f\!\left(\tfrac{z}{1-p\,b\,z}\right),$$

$$\left(\begin{pmatrix} 1 & \\ c & 1 \end{pmatrix} f\right)(z) = f(z - c),$$

$$\left(\begin{pmatrix} a & \\ & a^{-1} \end{pmatrix} f\right)(z) = \chi(a)\, f\!\left(a^{2}\, z\right).$$

(ii) The Lie algebra action on $(V_\chi)_{w_+}$ satisfies

$$(\mathfrak{x}_+\, f)\,(z) = -\tfrac{d}{d\,x}\, f\,(z)$$

$$(\mathfrak{x}_-\, f)\,(z) = c(\chi)\, p\, z\, f(z) + p\, z^2\, \tfrac{d}{d\,x}\, f(z).$$

The Lie algebra action on $(V_\chi)_{w_-}$ satisfies

$$(\mathfrak{x}_-\, f)\,(z) = -\tfrac{d}{d\,x}\, f\,(z),$$

$$(\mathfrak{x}_+\, f)\,(z) = c(\chi)\, p\, z\, f(z) + p\, z^2\, \tfrac{d}{d\,x}\, f\,(z).$$

*Proof:* In our case the Bruhat-Tits decomposition $G_1 = G_0\, W\, P$ is realised by

$$\begin{pmatrix} a & b \\ c & d \end{pmatrix} = \begin{cases} \begin{pmatrix} 1 & b/d \\ 0 & 1 \end{pmatrix} w_1 \begin{pmatrix} d^{-1} & \\ c & d \end{pmatrix}; & |b| < |d|, \\ \begin{pmatrix} 1 & \\ d/b & 1 \end{pmatrix} w_2 \begin{pmatrix} -b^{-1} & \\ -a & -b \end{pmatrix}; & |b| \geq |d|, \end{cases}$$

and for the projections $\pi_w^\pm$ we have $\pi_{w_+}^-(g) = \begin{pmatrix} 1 & b/d \\ & 1 \end{pmatrix}$, $\pi_{w_+}^+(g) = \begin{pmatrix} 1/d & \\ c & d \end{pmatrix}$, $\pi_{w_-}^-(g) = \begin{pmatrix} 1 & \\ c/a & 1 \end{pmatrix}$, $\pi_{w_-}^+(g) = \begin{pmatrix} a & b \\ & 1/a \end{pmatrix}$ for $g = \begin{pmatrix} a & b \\ c & d \end{pmatrix} \in G_0$. Now apply Lemma 11.1. □

In Section 11 we saw that dually to the $G_0$-stable decomposition $V_\chi = (V_\chi)_{w_+} \oplus (V_\chi)_{w_-}$ we have the decomposition of $D(G_0, K)$- modules

$$M := (V_\chi)'_b = M^+ \oplus M^-$$

where $M^\pm = D(G_0, K) \otimes_{D(P_{w_\pm},K)} \chi'_{w_\pm} \simeq D(N_{w_\pm}, K) \otimes_K \chi'_{w_\pm}$, and that there is a pair $(H_0, r_0) \in \Xi$ such that for all $(H, r) \geq (H_0, r_0)$,

$$M_{(H,r)} = M^+_{(H,r)} \oplus M^-_{(H,r)}$$

where $M^\pm_{(H,r)} = D_{(H \cap N_{w_\pm}, r)}(N_{w_\pm}, K) \otimes_K \chi'_{w_\pm}$ (Proposition 11.6). We fix elements

$$h_+ := 1 \in \mathbb{Z}_p \quad \text{corresponding to} \quad \begin{pmatrix} 1 & p \\ 0 & 1 \end{pmatrix} \in N_{w_+},$$
$$h_- := 1 \in \mathbb{Z}_p \quad \text{corresponding to} \quad \begin{pmatrix} 1 & 0 \\ 1 & 1 \end{pmatrix} \in N_{w_-}$$

under the above identification, and as usual we put $b_\pm := \delta_{h_\pm} - 1 \in M^\pm$. Then $\mathfrak{x}_\pm = \log(1 + b_\pm)$.

**Lemma 12.2.** Define operators $\Delta, \Pi_z : M^\pm \to M^\pm$ $(z \in \mathbb{Z}_p)$ by

$$\Delta(\textstyle\sum a_n b_\pm^n) := (1 + b_\pm) \sum n\, a_n\, b_\pm^{n-1},$$
$$\Pi_z(\textstyle\sum a_n b_\pm^n) := \sum a_n \big((1 + b_\pm)^z - 1\big)^n.$$

Then

$$\lambda(x \mapsto x\, f(x)) = \Delta\, \lambda\,(f),$$
$$\lambda(x \mapsto f(z\,x)) = \Pi_z\, \lambda\,(f)$$

for all $\lambda \in M^\pm$ and $f \in C^{\mathrm{an}}(\mathbb{Z}_p, K)$.

*Proof.* This is contained in [29], Lemma 4.3. □

Fix a pair $(H, r) \geq (H_0, r_0)$ in $\Xi$. Recall (Remark 7.3) that $\epsilon(r)$ is a $p$-power with the property that for all $z \in N_{w_\pm}^{\epsilon(r)}$ we have $\delta_z = \exp(\operatorname{Log} z)$ (convergent series in $M^\pm_{(H,r)}$). Moreover, the group $H \cap N_{w_\pm}$ is of the form $m^\pm\, \mathbb{Z}_p$ with a $p$-power $m^\pm$. We choose $\mathbb{N}_{<\epsilon(r)\, m^\pm} = \{0, \ldots, \epsilon(r)\, m^\pm - 1\}$ as a system of representatives for $N_{w_\pm} \big/ (H \cap N_{w_\pm})^{\epsilon(r)}$. For $n \in \mathbb{N}$ and $z \in \mathbb{Z}_p$ we put

$$T^\pm_{n,z} := \mathfrak{x}_\pm^n\, \delta_z = \log(1 + b_\pm)^n\, (1 + b_\pm)^z \in M^\pm.$$

According to Corollary 7.4 the family $\big(T^\pm_{n,z}\big)_{n \in \mathbb{N}, i \in \mathbb{N}_{<\epsilon(r)\, m^\pm}}$ is a basis of the Banach space $M^\pm_{(H,r)}$. We further introduce the elements

$$\tilde{T}^{\pm}_{n,z} := \sum_{j=0}^{n} \binom{n}{j} c_{n,j,z}\, T^{\pm}_{j,z} \in M^{\pm} \quad (z \neq 0),$$

$$\tilde{T}^{\pm}_{n,0} := T^{\pm}_{n,0} \in M^{\pm},$$

where for $j \leq n$ we put

$$c_{n,j,z} := z^{j-n}\, (c(\chi) + j)\, (c(\chi) + j + 1) \dots (c\,(\chi) + n - 1).$$

From now on, we suppose that $|c(\chi)| \leq 1$. Moreover, by eventually increasing $(H, r)$ we attain $||x_{\pm}||_{(H,r)} \geq 1$ (Remark 7.3 (2) a)). These assumptions are needed in order to prove the following result:

**Lemma 12.3.** The family $(\tilde{T}^{\pm}_{n,i})_{n\in\mathbb{N}, i\in\mathbb{N}_{<\epsilon(r)\, m^{\pm}}}$ is a basis of the Banach space $M^{\pm}_{(H,r)}$.

*Proof.* For every $n \in \mathbb{N}$ and every $z \in \mathbb{Z}_p$ the matrix $\left(\binom{k}{j} c_{k,j,z}\right)_{0\leq j,k\leq n}$ possesses the inverse $\left((-1)^{j+k} \binom{k}{j} c_{k,j,z}\right)_{0\leq j,k\leq n}$: indeed, $\sum_{l=0}^{n} \binom{l}{j} c_{l,j,z}\, (-1)^{l+k} \binom{k}{l} c_{k,l,z} = c_{k,j,z} \sum_{l=0}^{n} (-1)^{l+k} \binom{k}{l} \binom{l}{j} = \delta_{j,k}$ (Kronecker's symbol) and $\sum_{l=0}^{n} (-1)^{j+l} \binom{l}{j} c_{l,j,z} \binom{k}{l} c_{k,l,z} = c_{k,j,z} \sum_{l=0}^{n} (-1)^{j+l} \binom{k}{l} \binom{l}{j} = \delta_{j,k}$, since $c_{l,j,z}\, c_{k,l,z} = c_{k,j,z}$ for all $j \leq l \leq k$, $c_{j,j,z} = 1$ for all $j$, and $\sum_{l=0}^{n} (-1)^{l+k} \binom{k}{l} \binom{l}{j} = \delta_{j,k}$ for all $j$, $k \leq n$.

It follows that every $\lambda \in M^{\pm}_{(H,r)}$ may be written as

$$\lambda = \sum_{n\in\mathbb{N}, i\in\mathbb{N}_{<\epsilon(r)\, m^{\pm}}} b_{n,i}\, T^{\pm}_{n,i} = \sum_{i\in\mathbb{N}_{<\epsilon(r)\, m^{\pm}}} \sum_{n\in\mathbb{N}} b_{n,i} \sum_{j\in\mathbb{N}} (-1)^{j+n} \binom{n}{j} c_{n,j,i}\, \tilde{T}^{\pm}_{j,i}$$

with uniquely determined elements $b_{n,i} \in K$. We have to show that for each $i$ the last two sums on the right hand side commute, i.e. letting $\lambda_{n,j} := b_{n,i}\, (-1)^{j+n} \binom{n}{j} c_{n,j,i}\, \tilde{T}^{\pm}_{j,i}$ we have to show that for every $\varepsilon > 0$ there is exists an $n_0 \in \mathbb{N}$ such that $||\lambda_{n,j}||_{(H,r)} < \varepsilon$ whenever $n \geq n_0$ and $j \geq 0$ (cf. [25], Section 23). But by our assumption on $c(\chi)$ we have $|c_{j,k,i}| \leq 1$ for all $j$, $k$, and therefore

$$||\lambda_{n,j}||_{(H,r)} \leq \left\| b_{n,i} \binom{n}{j} \tilde{T}^{\pm}_{j,i} \right\|_{(H,r)} = \left\| b_{n,i} \binom{n}{j} \sum_{k=0}^{j} \binom{j}{k} c_{j,k,i}\, T^{\pm}_{k,i} \right\|_{(H,r)} \leq |b_{n,i}| \cdot \max_{0\leq k\leq n} \left\| T^{\pm}_{k,i} \right\|_{(H,r)}.$$

Moreover, the assumption $\|\mathfrak{x}_\pm\|_{(H,r)} \ge 1$ implies that $\max_{0\le k\le n} \left\|T^\pm_{k,i}\right\|_{(H,r)} = \left\|T^\pm_{n,i}\right\|_{(H,r)}$. Since the sum $\sum_{n\in\mathbb{N}} b_{n,i}\, T^\pm_{n,i}$ converges in $M^\pm_{(H,r)}$ we indeed have $\|\lambda_{n,j}\|_{(H,r)} \le |b_{n,i}| \left\|T^\pm_{n,i}\right\|_{(H,r)} < \varepsilon$ for large $n$, independently of $j$. □

**Lemma 12.4.** (i) The action of a torus element $s = \begin{pmatrix} a & \\ & a^{-1} \end{pmatrix} \in S_0$ on $M^+$ satisfies

$$s \cdot T^+_{n,z} = \chi(a)\, a^{2n}\, T^+_{n,a^2 z}.$$

The Lie algebra action on $M^+$ satisfies

$$\mathfrak{x}_+ \cdot T^+_{n,z} = T^+_{n+1,z}.$$

(ii) The action of $s = \begin{pmatrix} a & \\ & a^{-1} \end{pmatrix} \in S_0$ on $M^-$ satisfies

$$s \cdot \tilde{T}^-_{n,z} = \chi\left(a^{-1}\right) a^{-2n}\, \tilde{T}^-_{n,a^{-2} z}.$$

The Lie algebra action on $M^-$ satisfies

$$\mathfrak{x}_+ \cdot \tilde{T}^-_{n,z} = -p\, z^2\, \tilde{T}^-_{n+1,z}\ (z \ne 0),$$
$$\mathfrak{x}_+ \cdot \tilde{T}^-_{n,0} = -p\, n(c(\chi) + n - 1)\, \tilde{T}^-_{n-1,0}.$$

*Proof.* (i) For $y \in \mathbb{Z}_p$ and $f \in C^{\mathrm{an}}(\mathbb{Z}_p, K)$ define $f_y \in C^{\mathrm{an}}(\mathbb{Z}_p, K)$ by $f_y(x) = f(y\,x)$. According to Lemma 12.1 (i) and Lemma 12.2, $\left(s\, T^+_{n,z}\right)(f) = T^+_{n,z}\left(s^{-1} f\right) = \chi(a)\ T^+_{n,z}\left(f_{a^2}\right) = \chi(a)\ \Pi_{a^2}\left(\log(1+b_+)^n\, (1+b_+)^z\right)(f) = \chi(a)\ a^{2n}\ \log(1+b_+)^n\ (1+b_+)^{a^2 z}\ (f) = \chi(a)\ a^{2n}\ T^+_{n,a^2 z}\ (f)$, whence the first assertion of (i). Since the action of $D(N_{w_+}, K)$ on $M^+$ is just left multiplication the second assertion is also clear.

(ii) Similarly as in part (i) we calculate $\left(s\, T^-_{n,z}\right)(f) = \chi\left(a^{-1}\right) T^-_{n,z}\left(f_{a^{-2}}\right) = \chi\left(a^{-1}\right) a^{-2n}\, T^-_{n,a^{-2} z}(f)$ which already proves the first assertion for $z = 0$. For $z \ne 0$ we conclude

$$\begin{aligned}
s\, \tilde{T}^-_{n,z} &= \textstyle\sum_{j=0}^n \binom{n}{j} z^{j-n} (c(\chi)+j) \ldots (c(\chi)+n-1)\, s\, T^-_{j,z} \\
&= \chi\left(a^{-1}\right) \textstyle\sum_{j=0}^n \binom{n}{j} a^{-2j} z^{j-n}\, (c(\chi)+j) \ldots (c(\chi)+n-1)\, T^-_{j,a^{-2} z} \\
&= \chi\left(a^{-1}\right) a^{-2n} \textstyle\sum_{j=0}^n \binom{n}{j} \left(a^{-2} z\right)^{j-n} (c(\chi)+j) \ldots (c(\chi)+n-1)\, T^-_{j,a^{-2} z} \\
&= \chi\left(a^{-1}\right) a^{-2n}\, \tilde{T}^-_{j,a^{-2} z}\,.
\end{aligned}$$

Finally, by Lemma 12.1 (ii) the Lie algebra action on $M^-$ satisfies $(\mathfrak{x}_-\, \lambda)(f) = \lambda(-\mathfrak{x}_-\, f) =$

$\lambda\big(\frac{d}{dx}\, f\big)$ and $\big(-p^{-1}\,\mathfrak{x}_+\,\lambda\big)(f) = \lambda\big(p^{-1}\,\mathfrak{x}_+\, f\big) = c(\chi)\,\Delta\,\lambda(f) + \Delta^2\,\lambda\big(\frac{d}{dx}\, f\big) = c(\chi)\,\Delta\,\lambda(f) + \mathfrak{x}_-\,\Delta^2\,\lambda(f)$ in the notation of Lemma 12.2. Turning to our elements $T^-_{n,z}$, $\tilde{T}^-_{n,z}$ we obtain

$$
\begin{aligned}
\Delta\, T^-_{n,z} &= (1+b_-)\left(\frac{n\log(1+b_-)^{n-1}}{1+b}\,(1+b_-)^z + \log(1+b)^n\, z\,(1+b)^{z-1}\right) \\
&= n\, T^-_{n-1,z} + z\, T^-_{n,z}\,, \\
\Delta^2\, T^-_{n,z} &= n\,(n-1)\, T^-_{n-2,z} + 2\,n\,z\, T^-_{n-1,z} + z^2\, T^-_{n,z}\,, \\
-p^{-1}\,\mathfrak{x}_+\, T^-_{n,z} &= c(\chi)\, n\, T^-_{n-1,z} + c(\chi)\, z\, T^-_{n,z} \\
&+ n\,(n-1)\, T^-_{n-1,z} + 2\,n\,z\, T^-_{n,z} + z^2\, T^-_{n+1,z} \\
&= n\,(c\,(\chi) + n - 1)\, T^-_{n-1,z} + z\,(c\,(\chi) + 2\,n)\, T^-_{n,z} + z^2\, T^-_{n+1,z}\,,
\end{aligned}
$$

thereby proving the second assertion of (ii) in case $z = 0$; for $z \neq 0$ we calculate

$$
\begin{aligned}
-p^{-1}\,\mathfrak{x}_+\,\tilde{T}^-_{n,z} &= \textstyle\sum_{i=0}^{n}\binom{n}{i}\, z^{i-n}(c(\chi)+i)\ldots(c(\chi)+n-1)\left(-p^{-1}\,\mathfrak{x}_+\, T^-_{i,z}\right) \\
&= \textstyle\sum_{i=-1}^{n-1}\binom{n}{i+1}\, z^{i+1-n}(c\,(\chi)+i)\ldots(c\,(\chi)+n-1)\,(i+1)\, T^-_{i,z} \\
&+ \textstyle\sum_{i=0}^{n}\binom{n}{i}\, z^{i-n}(c\,(\chi)+i)\ldots(c\,(\chi)+n-1)\, z\,(c(\chi)+2\,i)\, T^-_{i,z} \\
&+ \textstyle\sum_{i=1}^{n+1}\binom{n}{i-1}\, z^{i-1-n}(c\,(\chi)+i-1)\ldots(c\,(\chi)+n-1)\, z^2\, T^-_{i,z} \\
&= \textstyle\sum_{i=0}^{n+1}\binom{n+1}{i}\, z^{i+1-n}\,(c\,(\chi)+i)\ldots(c\,(\chi)+n)\, T^-_{i,z} \\
&= z^2\,\tilde{T}^-_{n+1,z}\,,
\end{aligned}
$$

because $\binom{n}{i+1}(i+1) + \binom{n}{i}(c(\chi)+2\,i) + \binom{n}{i-1}(c(\chi)+i-1) = \binom{n-1}{i}n + 2\binom{n-1}{i-1}n + \binom{n-1}{i-2}n + \binom{n}{i}c(\chi) + \binom{n}{i-1}c(\chi) = \binom{n+1}{i}(c(\chi)+n)$ by the recursion formula. This proves the second assertion of (ii). □

Now fix a number $k \in \mathbb{N}$; in case $-c(\chi) \in \mathbb{N}$ assume $k > -c(\chi)$. Note that since the decomposition $M = M^+ \oplus M^-$ is $\mathfrak{n}$-stable we have $M_{(H,r)} \big/ \mathfrak{n}^k\, M_{(H,r)} = \big(M^+_{(H,r)} \big/ \mathfrak{n}^k\, M^+_{(H,r)}\big) \oplus \big(M^-_{(H,r)} \big/ \mathfrak{n}^k\, M^-_{(H,r)}\big)$.

**Corollary 12.5.** (i) Put $I^+_k := \mathbb{N}_{<k} \times \mathbb{N}_{<\epsilon(r)\, m^+}$. The elements

$$T^+_{n,i} + \mathfrak{n}^k\, M^+_{(H,r)} \quad ((n,\, i) \in I^+_k)$$

constitute a basis of the $K$-vector space $M^+_{(H,r)} \big/ \mathfrak{n}^k\, M^+_{(H,r)}$.

(ii) Put

$$I_k^- := \mathbb{N}_{<k} \times \mathbb{N}_{<\epsilon(r)\,m^-} - \begin{cases} \{(n,\,0);\, n > -c(\chi)\} & \text{if } -c(\chi) \in \mathbb{N}, \\ \{(n,\,0);\, n \geq 0\} & \text{if } -c(\chi) \notin \mathbb{N}. \end{cases}$$

The elements

$$\tilde{T}^-_{n,i} + \mathfrak{n}^k\, M^-_{(H,r)} \quad ((n,\,i) \in I_k^-)$$

constitute a basis of the $K$-vector space $M^-_{(H,r)} \big/ \mathfrak{n}^k\, M^-_{(H,r)}$.

*Proof.* This follows from the description of the action of $\mathfrak{n}$ on $M^\pm$ given in Lemma 12.4. For (ii), in the case $-c(\chi) \in \mathbb{N}$ observe that $\mathfrak{x}_+\, \tilde{T}^-_{n,0} = 0$ if and only if $n = 0$ or $n = 1 - c(\chi)$, and that $\mathfrak{x}_+\, \tilde{T}^-_{n,0}$ is a nonzero scalar multiple of $\tilde{T}^-_{n-1,0}$ in all other cases; hence $\mathfrak{x}_+^k \cdot \mathrm{Span}\big(\tilde{T}^-_{n,0};\, n \geq 0\big) = \mathrm{Span}\big(\tilde{T}^-_{n,0};\, n \geq 1 - c(\chi)\big)$ provided $k \geq 1 - c(\chi)$. □

For any $z \in \mathbb{Z}_p$ we define elements $Q^\pm(z) \in \epsilon(r)\, m^\pm\, \mathbb{Z}_p$ and $R^\pm(z) \in \mathbb{N}_{<\epsilon(r)\,m^\pm}$ by $z = Q^\pm(z) + R^\pm(z)$.

**Lemma 12.6.** Let $s = \begin{pmatrix} a & \\ & a^{-1} \end{pmatrix} \in S_0$.

(i) The coefficients $b^+_{(n,i),(n',i')} \in K$ defined by the equations

$$s\, T^+_{n',i'} = \textstyle\sum_{(n,i)\in I_k^+} b^+_{(n,i),(n',i')}\, T^+_{n,i} \quad ((n',\,i') \in I_k^+)$$

in $M^+_{(H,r)} \big/ \mathfrak{n}^k\, M^+_{(H,r)}$ satisfy

$$b^+_{(n,i),(n',i')} = 0 \;\text{ if } n < n' \text{ or } i \neq R^+\big(a^2\, i'\big),$$

$$b^+_{(n,i),(n,i)} = \chi(a)\, a^{2\,n} \;\text{ if } i = R^+\big(a^2\, i\big).$$

(ii) The coefficients $b^-_{(n,i),(n',i')} \in K$ defined by the equations

$$s\, \tilde{T}^-_{n',i'} = \textstyle\sum_{(n,i)\in I_k^-} b^-_{(n,i),(n',i')}\, \tilde{T}^-_{n,i} \quad ((n',\,i') \in I_k^-)$$

in $M^-_{(H,r)} \big/ \mathfrak{n}^k\, M^-_{(H,r)}$ satisfy

$$b^-_{(n,i),(n',i')} = 0 \;\text{ if } n > n' \text{ or } i \neq R^-\big(a^{-2}\, i'\big),$$

$$b^-_{(n,i),(n,i)} = \chi^{-1}(a) \sum_{j,m\in\mathbb{N}} d_{n,j,m}\, a^{-2\,j}\, \big(a^{-2} - 1\big)^m \;\text{ if } i = R^-\big(a^{-2}\, i\big),$$

where

$$d_{n,j,m} := \frac{1}{m!}\, c_{j+m,j,1} (-1)^{j+m+n} \binom{n}{j} \binom{j+m}{n}.$$

*Proof.* (i) First, $s\, T^{+}_{n',i'} = \chi(a)\, a^{2\,n'}\, T^{+}_{n',a^2\,i'}$ by Lemma 12.4 (i). But $T^{+}_{n',a^2\,i'} = \varkappa_+^{n'}\, \delta_1^{a^2\,i'} = \varkappa_+^{n'}\, \delta_1^{R^+\left(a^2\,i'\right)}\, \delta_1^{Q^+\left(a^2\,i'\right)}$, and $\delta_1^{Q^+\left(a^2\,i'\right)} = \exp\!\left(Q^+\!\left(a^2\,i'\right)\log\delta_1\right) = \exp\!\left(Q^+\!\left(a^2\,i'\right)\varkappa_+\right) = \sum_{m\geq 0} \frac{1}{m!}\, Q^+\!\left(a^2\,i'\right)^m \varkappa_+^m$. Hence

$$s\, T^{+}_{n',i'} = \chi(a)\, a^{2\,n'} \sum_{m\geq 0} \frac{1}{m!}\, Q^+\!\left(a^2\,i'\right)^m\, T^{+}_{n'+m,R^+\left(a^2\,i'\right)}.$$

This proves our assertion, since $n' + m = n'$ implies $\frac{1}{m!}\, Q^+\!\left(a^2\,i\right)^m = 1$.

(ii) Similarly as in part (i), this time using Lemma 12.4 (ii), we calculate

$$s\, T^{-}_{n',i'} = \chi\!\left(a^{-1}\right) a^{-2\,n'} \sum_{m\geq 0} \frac{1}{m!}\, Q^-\!\left(a^{-2}\,i'\right)^m\, T^{-}_{n'+m,R^-\left(a^{-2}\,i'\right)}.$$

Hence

$$\begin{aligned}
s\, \tilde{T}^{-}_{n',i'} &= \sum_{j\geq 0} \binom{n'}{j} c_{n',j,i'}\, s\, T^{-}_{j,i'} \\
&= \sum_{j,m\geq 0} \binom{n'}{j} c_{n',j,i'}\, \chi\!\left(a^{-1}\right) a^{-2\,j}\, \frac{1}{m!}\, Q^-\!\left(a^{-2}\,i'\right)^m\, T^{-}_{j+m,R^-\left(a^{-2}\,i'\right)} \\
&= \sum_{j,m,n\geq 0} \binom{n'}{j} c_{n',j,i'}\, \chi\left(a^{-1}\right) a^{-2\,j}\, \frac{1}{m!} \\
&\qquad Q^-\!\left(a^{-2}\,i'\right)^m (-1)^{j+m+n} \binom{j+m}{n} c_{j+m,n,R\left(a^{-2}\,i'\right)}\, \tilde{T}^{-}_{n,R^-\left(a^{-2}\,i'\right)}.
\end{aligned}$$

It already follows that $b_{(n',i'),(n,i)} = 0$ if $n > n'$ or $i \neq R^-\!\left(a^{-2}\,i'\right)$. If $i = R^-\!\left(a^{-2}\,i\right)$ then $Q^-\!\left(a^{-2}\,i\right) = a^{-2}\,i - i = i\left(a^{-2} - 1\right)$, and therefore

$$\begin{aligned}
b_{(n,i),(n,i)} &= \sum_{j,m} \binom{n}{j} c_{n,j,i}\, \chi\left(a^{-1}\right) a^{-2\,j}\, \frac{1}{m!}\, i^m\!\left(a^{-2} - 1\right)^m (-1)^{j+m+n} \binom{j+m}{n} c_{j+m,n,i} \\
&= \sum_{j,m} \binom{n}{j} \chi\left(a^{-1}\right) a^{-2\,j}\, \frac{1}{m!} \left(a^{-2} - 1\right)^m (-1)^{j+m+n} \binom{j+m}{n} c_{j+m,j,1}
\end{aligned}$$

because $c_{n,j,i}\, c_{j+m,n,i} = c_{j+m,j,i} = i^{-m}\, c_{j+m,j,1}$. □

**Theorem 12.7.** Assume $|c(\chi)| \leq 1$. The $G_0$-representation $V_\chi$ possesses an $\mathfrak{n}$-character on $S_0 \cap G^{\mathrm{reg}}$ given by

$$\theta_{V_\chi}\left(\begin{pmatrix} a & \\ & a^{-1}\end{pmatrix}\right) = \frac{\chi(a)^{-1}}{\left|1-a^{-2}\right|_p\left(1-a^{-2}\right)} + \frac{\chi(a)}{\left|1-a^{2}\right|_p\left(1-a^{2}\right)} - \frac{\chi(a)\,a^{2\,m_0}}{1-a^2}$$

where

$$m_0 := \begin{cases} 1-c(\chi) & \text{if } -c(\chi)\in\mathbb{N},\\ 0 & \text{otherwise.}\end{cases}$$

*Proof.* Using Proposition 9.5 we will treat $(V_\chi)_{w_+}$ and $(V_\chi)_{w_-}$ separately.

**Step 1:** We fix a pair $(H, r)$ and a number $k\in\mathbb{N}$ (both large enough), and we determine the formal character of the representation $M^{\pm}_{(H,r)}\big/\mathfrak{n}^k\,M^{\pm}_{(H,r)}$ of $S_0$.

First of all this representation possesses a formal character $\Theta^{\pm}_k\in\mathbb{Z}[X(S_0)]$ (over the field $C$) whose value in any $s=\begin{pmatrix} a & \\ & a^{-1}\end{pmatrix}\in S_0$ is equal to the trace of the operator on $M^{\pm}_{(H,r)}\big/\mathfrak{n}^k\,M^{\pm}_{(H,r)}$ defined by $s$ (Remark 2.2 (ii)). Hence, in the notations of Lemma 12.6 (i),

$$\mathrm{ev}_{S_0}\left(\Theta^{+}_k\right)(s) = \sum_{(n,i)\in I^{+}_k} b^{+}_{(n,i),(n,i)} = \chi^{-1}(s)\,l(s)\sum_{n=0}^{k-1} a^{2\,n}$$

where $l(s)\in\mathbb{N}$ is the number of fixpoints of the permutation $i\mapsto R^{+}\left(a^2\,i\right)$ of the set $\mathbb{N}_{<\epsilon(r)\,m^{+}}$. Let $\varepsilon: S_0\to K^{\times}$ denote the character $\begin{pmatrix} a & \\ & a^{-1}\end{pmatrix}\mapsto a^{-1}$, and let $\Psi^{+}\in\mathbb{Z}[X(S_0)]$ denote the formal character of the matrix representation of $S_0$

$$\begin{pmatrix} a & \\ & a^{-1}\end{pmatrix}\mapsto\left(\delta_{i,R^{+}\left(a^2\,i'\right)}\right)_{(i,i')\in\mathbb{N}_{<\epsilon(r)\,m^{+}}\times\mathbb{N}_{<\epsilon(r)\,m^{+}}}.$$

Then $l(s)=\mathrm{ev}_{S_0}(\Psi^{+})\,(s)$, and hence

$$\mathrm{ev}_{S_0}\left(\Theta^{+}_k\right)(s) = \mathrm{ev}_{S_0}\left(e\left(\chi^{-1}\right)\cdot\Psi^{+}\cdot\sum_{n=0}^{k-1} e(\varepsilon)^{-2\,n}\right)(s)$$

by multiplicativity of the evaluation map. But since this map is also injective (Remark 1.2 (i)) we deduce

$$\Theta^{+}_k = e\left(\chi^{-1}\right)\cdot\Psi^{+}\cdot\sum_{n=0}^{k-1} e(\varepsilon)^{-2\,n}.$$

Similarly, in the notations of Lemma 12.6 (ii),

$$\mathrm{ev}_{S_0}(\Theta_k^-)(s) = \sum_{(n,i)\in I_k^-} b^-_{(n,i),(n,i)} = \sum_{n=0}^{k-1}\sum_{i=1}^{\epsilon(r)\,m^- -1} b^-_{(n,i),(n,i)} + \sum_{n=0}^{-c(\chi)} b^-_{(n,0),(n,0)};$$

here we understand the sum $\sum_{n=0}^{-c(\chi)}(\,...)$ to be zero if $-c(\chi)\notin\mathbb{N}$. Letting $l^-(s)\in\mathbb{N}$ denote the number of fixpoints of the permutation $i\mapsto R^-\left(a^{-2}\,i\right)$ of the set $\mathbb{N}^*_{<\epsilon(r)\,m^-} = \mathbb{N}_{<\epsilon(r)\,m^-} - \{0\}$ we have

$$\sum_{n=0}^{k-1}\sum_{i=1}^{\epsilon(r)\,m^- -1} b^-_{(n,i),(n,i)} = \chi(s)\,l^-(s)\sum_{n=0}^{k-1}\sum_{j,m\in\mathbb{N}} d_{n,j,m}\,a^{-2\,j}\left(a^{-2}-1\right)^m.$$

Moreover, if $-c(\chi)\in\mathbb{N}$ then we calculate

$$\begin{aligned}\sum_{n=0}^{-c(\chi)} b^-_{(n,0),(n,0)} &= \chi(s)\sum_{n=0}^{-c(\chi)}\sum_{j=0}^{n}\sum_{m=0}^{\infty}\frac{c_{j+m,j,1}}{m!}\,a^{-2\,j}\left(a^{-2}-1\right)^m(-1)^{j+m+n}\binom{n}{j}\binom{j+m}{n}\\ &= \chi(s)\sum_{j=0}^{-c(\chi)}\sum_{m=0}^{-c(\chi)-j}\sum_{n=0}^{j+m}\frac{c_{j+m,j,1}}{m!}\,a^{-2\,j}\left(a^{-2}-1\right)^m(-1)^{j+m+n}\binom{n}{j}\binom{j+m}{n}\\ &= \chi(s)\sum_{j=0}^{-c(\chi)} a^{-2\,j};\end{aligned}$$

the second equality holds since $c_{j+m,j,1}=0$ whenever $j\le -c(\chi) < j+m$, the third equality follows from $\sum_{n\in\mathbb{N}}(-1)^{j+m+n}\binom{n}{j}\binom{j+m}{n}=\delta_{m,0}$. Hence

$$\begin{aligned}\mathrm{ev}_{S_0}(\Theta_k^-)(s) &= \chi(s)\,l^-(s)\sum_{n=0}^{k-1}\sum_{j,m\in\mathbb{N}} d_{n,j,m}\,\varepsilon(s)^{2\,j}\left(\varepsilon(s)^2-1\right)^m\\ &+ \chi(s)\sum_{j=0}^{-c(\chi)}\varepsilon(s)^{2\,j}.\end{aligned}$$

Again by the injectivity of $\mathrm{ev}_{S_0}$ we deduce

$$\Theta_k^- = e(\chi)\,e(\Psi^-)\sum_{n=0}^{k-1}\sum_{j,m\in\mathbb{N}} d_{n,j,m}\,e(\varepsilon)^{2\,j}\left(e(\varepsilon)^2-1\right)^m + e(\chi)\sum_{j=0}^{-c(\chi)} e(\varepsilon)^{2\,j}$$

where $\Psi^-\in\mathbb{Z}[\![X(S_0)]\!]$ denotes the formal character of the matrix representation of $S_0$

$$\begin{pmatrix} a & \\ & a^{-1}\end{pmatrix}\mapsto\left(\delta_{i,R^-\left(a^{-2}\,i'\right)}\right)_{(i,i')\in\mathbb{N}^*_{<\epsilon(r)\,m^-}\times\mathbb{N}^*_{<\epsilon(r)\,m^-}}.$$

**Step 2:** We keep the pair $(H,r)$ from the first step and form the limit $k\to\infty$ in $\mathbb{Z}[\![X(S_0)]\!]$ of the formal characters just obtained:

$$\Theta^+_{(H,r)} := \lim_{k\to\infty} \Theta^+_k \;=\; e\big(\chi^{-1}\big)\,\Psi^+ \sum_{j\in\mathbb{N}} e(\varepsilon)^{-2j}\,,$$

$$\Theta^-_{(H,r)} := \lim_{k\to\infty} \Theta^-_k \;=\; e(\chi)\,e(\Psi^-) \sum_{n,j,m\in\mathbb{N}} d_{n,j,m}\, e(\varepsilon)^{2j} \left(e(\varepsilon)^2-1\right)^m \;+\; e(\chi) \sum_{j=0}^{-c(\chi)} e(\varepsilon)^{2j}$$

$$= e(\chi)\,\Psi^- \sum_{j\in\mathbb{N}} e(\varepsilon)^{2j} \;+\; e(\chi)\sum_{j=0}^{-c(\chi)} e(\varepsilon)^{2j}$$

since

$$\begin{aligned}\textstyle\sum_{n,j,m} d_{n,j,m}\, e(\varepsilon)^{2j}\left(e(\varepsilon)^2-1\right)^m &= \textstyle\sum_{j,m} \frac{1}{m!}\, c_{j+m,j,1}\, e(\varepsilon)^{2j}\left(e(\varepsilon)^2-1\right)^m \sum_n (-1)^{j+m+n}\binom{n}{j}\binom{j+m}{n} \\ &= \textstyle\sum_{j,m} \frac{1}{m!}\, c_{j+m,j,1}\, e(\varepsilon)^{2j}\left(e(\varepsilon)^2-1\right)^m \delta_{m,0} \\ &= \textstyle\sum_j e(\varepsilon)^{2j}\,.\end{aligned}$$

**Step 3:** Evaluation of $\Theta^\pm_{(H,r)}$ on the subset $S_0 \cap G^{\text{reg}} = \{\begin{pmatrix} a & \\ & a^{-1}\end{pmatrix};\ a \neq \pm 1\}$ of regular elements: In $\mathbb{Z}[\![X(S_0)]\!]$ we have the equalities $\left(1-e(\varepsilon)^{\mp 2}\right)\sum_{j\in\mathbb{N}} e(\varepsilon)^{\mp 2j} = 1$, whence $\Theta^\pm_{(H,r)}$ is contained in the submodule $\mathbb{Z}[\![X(S_0)]\!]_{S_0\cap G^{\text{reg}}}$ (cf. Definition 1.1), and

$$\text{ev}_{S_0\cap G^{\text{reg}}}\Big(\textstyle\sum_{j\in\mathbb{N}} e(\varepsilon)^{\mp 2j}\Big)\left(\begin{pmatrix} a & \\ & a^{-1}\end{pmatrix}\right) = \left(1-a^{\pm 2}\right)^{-1}.$$

We claim that if $\epsilon(r)\,m^\pm = |\epsilon(r)\,m^\pm|_p^{-1} \geq \left|a^{\pm 2}-1\right|_p^{-1}$ then

$$\text{ev}_{S_0}(\Psi^+)\left(\begin{pmatrix} a & \\ & a^{-1}\end{pmatrix}\right) = \left|1-a^2\right|_p^{-1},$$

$$\text{ev}_{S_0}(\Psi^-)\left(\begin{pmatrix} a & \\ & a^{-1}\end{pmatrix}\right) = \left|1-a^{-2}\right|_p^{-1} - 1.$$

Indeed, $\text{ev}_{S_0}(\Psi^+)(s)$ is the number of fixpoints of the permutation $i \mapsto R^+\big(a^2 i\big)$ of the set $\mathbb{N}_{<\epsilon(r)\,m^+}$. For each $i$ we have $i = R^+\big(a^2 i\big) \iff v_p\big(a^2 i - i\big) \geq v_p(\epsilon(r)\,m^+) \iff v_p(i) \geq v_p(\epsilon(r)\,m^+) - v_p\big(a^2-1\big)$; since $v_p(\epsilon(r)\,m^+) \geq v_p\big(a^2-1\big)$ by assumption, there are precisely $p^{v_p\left(a^2-1\right)} = \left|1-a^2\right|_p^{-1}$ elements $i$ in $\mathbb{N}_{<\epsilon(r)\,m^+}$ satisfying that condition (Lemma 3.2). Similarly, there are $p^{v_p\left(a^{-2}-1\right)} = \left|1-a^{-2}\right|_p^{-1}$ elements $i$ in $\mathbb{N}_{<\epsilon(r)\,m^-}$ satisfying $i = R^-\big(a^{-2} i\big)$. Since $\text{ev}_{S_0}(\Psi^-)(s)$ is the number of fixpoints of the permutation $i \mapsto R^-\big(a^{-2} i\big)$ of the set $\mathbb{N}^*_{<\epsilon(r)\,m^-}$ the second formula follows from the additional observation that 0 always is a fixpoint of that permutation. This proves our claim.

By multiplicativity of the evaluation map we obtain

$$
\begin{aligned}
\mathrm{ev}_{S_0\cap G^{\mathrm{reg}}}\left(\Theta^+_{(H,r)}\right)(s) &= \chi(a)\left|1-a^2\right|_p^{-1}\left(1-a^2\right)^{-1},\\
\mathrm{ev}_{S_0\cap G^{\mathrm{reg}}}\left(\Theta^-_{(H,r)}\right)(s) &= \chi\left(a^{-1}\right)\left(\left|1-a^{-2}\right|_p^{-1}-1\right)\left(1-a^{-2}\right)^{-1}\\
&+ \chi\left(a^{-1}\right)\sum_{j=0}^{-c(\chi)} a^{-2j}\\
&= \chi\left(a^{-1}\right)\left|1-a^{-2}\right|_p^{-1}\left(1-a^{-2}\right)^{-1}\\
&- \begin{cases} \chi\left(a^{-1}\right)\left(1-a^{-2}\right)^{-1} & \text{if } -c(\chi)\notin\mathbb{N}\\ \chi\left(a^{-1}\right)\left(1-a^{-2}\right)^{-1}a^{-2(1-c(\chi))} & \text{if } -c(\chi)\in\mathbb{N}\end{cases}
\end{aligned}
$$

provided $\epsilon(r)\,m^{\pm} \geq \left|a^{\pm 2}-1\right|_p^{-1}$.

**Step 4:** For any given $s$, there certainly exists a pair $(H, r) = (H_1, r_1) \in \Xi$ satisfying $\epsilon(r)\,m^{\pm} \geq \left|a^{\pm 2}-1\right|_p^{-1}$. But then this condition is also satisfied for each pair $(H, r) \geq (H_1, r_1)$: this follows from Lemma 9.2 (ii) and the fact that $\epsilon(r)\,m^{\pm}$ is equal to the order of $N_{w_\pm}\big/(N_{w_\pm}\cap H)^{\epsilon(r)}$. Furthermore, the expression $\mathrm{ev}_{S_0\cap G^{\mathrm{reg}}}\left(\Theta^{\pm}_{(H,r)}\right)(s)$ does not depend on the particular choice of $(H, r)$. Hence the family $\left(\mathrm{ev}_{S_0\cap G^{\mathrm{reg}}}\left(\Theta^{\pm}_{(H,r)}\right)\right)_{(H,r)}$ converges weakly. Composing the limit function with the map $s \mapsto s^{-1}$ then gives the asserted formula. □

**Remark 12.8.** Here is an illustration of how our definition is compatible with known results. Suppose that $\chi$ is a *smooth* character, i.e. $c(\chi) = 0$. We can then consider the smooth principal series representation $V_\chi^\infty = \mathrm{Ind}_P^G(\chi)^\infty$ of $G$ induced from $\chi$, and with the above method one can show that this representation possesses an $\mathfrak{n}$-character on $S_0 \cap G^{\mathrm{reg}}$, given by

$$
\theta_{V_\chi^\infty}\left(\begin{pmatrix} a & \\ & a^{-1}\end{pmatrix}\right) = \frac{\chi(a)^{-1}}{\left|1-a^{-2}\right|_p} + \frac{\chi(a)}{\left|1-a^{2}\right|_p}.
$$

(However, since the action of $\mathfrak{n}$ is trivial the computations become much simpler here; the two components $\left(V_\chi^\infty\right)_{w_\pm}$ identify with $C^\infty(\mathbb{Z}_p, K)$.) As predicted by Theorem 10.5 this formula coincides with the known character formula of the smooth principal series of $\mathrm{SL}_2(\mathbb{Q}_p)$ (cf. [17], p. 200; note that they use *normalized* induction).

Letting $\varepsilon$ denote the character $\begin{pmatrix} a & \\ & a^{-1}\end{pmatrix} \mapsto a^{-1}$ we obtain the exact sequence of representations of $G$

$$0 \longrightarrow V_\chi^\infty \overset{\subset}{\longrightarrow} V_\chi \overset{\frac{d}{dx}}{\longrightarrow} V_{\varepsilon^2 \chi} \longrightarrow 0\,;$$

here $\frac{d}{dx}$ denotes the derivation on each component $(V_\chi)_{w_\pm} = C^{\mathrm{an}}(\mathbb{Z}_p, K)$. Theorem 12.7 yields

$$\theta_{V_\chi}\left(\begin{pmatrix} a & \\ & a^{-1}\end{pmatrix}\right) = \frac{\chi(a)^{-1}}{\left|1-a^{-2}\right|_p \left(1-a^{-2}\right)} + \frac{\chi(a)}{\left|1-a^{2}\right|_p \left(1-a^{2}\right)} - \frac{\chi(a)\, a^2}{1-a^2},$$

$$\theta_{V_{\varepsilon^2 \chi}}\left(\begin{pmatrix} a & \\ & a^{-1}\end{pmatrix}\right) = \frac{\chi(a)^{-1}\, a^{-2}}{\left|1-a^{-2}\right|_p \left(1-a^{-2}\right)} + \frac{\chi(a)\, a^2}{\left|1-a^{2}\right|_p \left(1-a^{2}\right)} - \frac{\chi(a)\, a^2}{1-a^2}$$

(note that $c\left(\varepsilon^2 \chi\right) = 2 \notin -\mathbb{N}$). It is easily verified that the alternating sum of the three characters is zero.

More generally, assume $-c(\chi) \in \mathbb{N}$. The character $\chi$ is then the product of the algebraic character $\varepsilon^{c(\chi)}$ and the smooth character $\varepsilon^{-c(\chi)}\,\chi$. As an application of the theory developped in [29] and [30] Schneider and Teitelbaum prove that in this case (and only in this case) $V_\chi$ is reducible. More exactly they construct an exact sequence of $G$- representations

$$0 \longrightarrow V^{\mathrm{alg}}_{\varepsilon^{c(\chi)}} \otimes V^\infty_{\varepsilon^{-c(\chi)}\,\chi} \longrightarrow V_\chi \longrightarrow V_{\varepsilon^{2-2\,c(\chi)}\,\chi} \longrightarrow 0$$

where the algebraic induction $V^{\mathrm{alg}}_{\varepsilon^{c(\chi)}} = \mathrm{ind}_P^G(\varepsilon^{c(\chi)})$ is at the same time the irreducible $\mathbb{Q}_p$- rational $G$-representation of highest weight $\varepsilon^{c(\chi)}$ w.r.t. the upper Borel subgroup (cf. [30], Section 4). We have $c\left(\varepsilon^{2-2\,c(\chi)}\,\chi\right) = 2 - c(\chi) \notin -\mathbb{N}$, hence Theorem 12.7 gives the characters

$$\theta_{V_\chi}\left(\begin{pmatrix} a & \\ & a^{-1}\end{pmatrix}\right) = \frac{\chi(a)^{-1}}{\left|1-a^{-2}\right|_p \left(1-a^{-2}\right)} + \frac{\chi(a)}{\left|1-a^{2}\right|_p \left(1-a^{2}\right)} - \frac{\chi(a)\, a^{2-2\,c(\chi)}}{1-a^2},$$

$$\theta_{V_{\varepsilon^{2-2c(\chi)}\chi}}\left(\begin{pmatrix} a & \\ & a^{-1}\end{pmatrix}\right) = \frac{\chi(a)^{-1}\, a^{-2+2\,c(\chi)}}{\left|1-a^{-2}\right|_p \left(1-a^{-2}\right)} + \frac{\chi(a)\, a^{2-2\,c(\chi)}}{\left|1-a^{2}\right|_p \left(1-a^{2}\right)} - \frac{\chi(a)\, a^{2-2\,c(\chi)}}{1-a^2}.$$

The character of the smooth $G$-representation $V^\infty_{\varepsilon^{-c(\chi)}\,\chi}$ is given by

$$\theta_{V^\infty_{\varepsilon^{-c(\chi)}\chi}}\left(\begin{pmatrix} a & \\ & a^{-1}\end{pmatrix}\right) = \frac{\chi(a)^{-1}\, a^{c(\chi)}}{\left|1-a^{-2}\right|_p} + \frac{\chi(a)\, a^{-c(\chi)}}{\left|1-a^{2}\right|_p}.$$

Finally, for the character of the rational $G$-representation $V^{\mathrm{alg}}_{\varepsilon^{c(\chi)}}$ the Weyl character formula yields

$$\left(a - a^{-1}\right)\theta_{V^{\mathrm{alg}}_{\varepsilon^{c(\chi)}}}\left(\begin{pmatrix} a & \\ & a^{-1}\end{pmatrix}\right) = a^{c(\chi)-1} - a^{1-c(\chi)}.$$

(cf. [18], p. 298, Theorem 7.1.1). Again, it is easy to check that the alternating sum of the characters is zero:

$$\theta_{V_{\varepsilon^{2-2c(\chi)}\chi}}\left(\begin{pmatrix} a & \\ & a^{-1}\end{pmatrix}\right) - \theta_{V_\chi}\left(\begin{pmatrix} a & \\ & a^{-1}\end{pmatrix}\right) = \frac{\chi(a)^{-1}\left(a^{-2+2\,c(\chi)}-1\right)}{\left|1-a^{-2}\right|_p\left(1-a^{-2}\right)} + \frac{\chi(a)\left(a^{2-2\,c(\chi)}-1\right)}{\left|1-a^{2}\right|_p\left(1-a^{2}\right)}$$

$$= \frac{\chi(a)^{-1}\,a^{c(\chi)}}{\left|1-a^{-2}\right|_p}\cdot\frac{a^{-1+c(\chi)}-a^{1-c(\chi)}}{a-a^{-1}} + \frac{\chi(a)\,a^{-c(\chi)}}{\left|1-a^{2}\right|_p}\cdot\frac{a^{1-c(\chi)}-a^{-1+c(\chi)}}{a^{-1}-a}$$

$$= \theta_{V^{\infty}_{\varepsilon^{-c(\chi)}\chi}}\left(\begin{pmatrix} a & \\ & a^{-1}\end{pmatrix}\right)\cdot\theta_{V^{\mathrm{alg}}_{\varepsilon^{c(\chi)}}}\left(\begin{pmatrix} a & \\ & a^{-1}\end{pmatrix}\right).$$

In determining the $\mathfrak{n}$-character of $V_\chi$ it became visible that the calculations on the component $(V_\chi)_{w_+}$ are much easier to do than on $(V_\chi)_{w_-}$, the reason being that $\mathfrak{n} = \mathfrak{n}_{w_+}$ acts on the dual $\left((V_\chi)_{w_+}\right)'_b \simeq D(N_w, K)$ simply by left multiplication. This fact will be exploited further in the next section.

## 13. The principal series of the Iwahori subgroup

We resume the assumptions and notations of Section 11 and fix a Weyl group element $w \in W$; thus $G_0$ is the Iwahori subgroup of the reductive group $G$ of the same type as $P$, and we have the decomposition $G_0 = N_w \cdot P_w$ (direct span). Let

$$\mathfrak{n}_w = \mathfrak{g}_{a_1} \oplus \ldots \oplus \mathfrak{g}_{a_d}$$

be the root space decomposition of the Lie algebra $\mathfrak{n}_w$ of $N_w$. Let $G^{\mathrm{reg}} \subset G$ denote the set of regular elements, and put $S' := S_0 \bigcap G^{\mathrm{reg}}$. By the end of this section we will be able to prove the following result:

**Theorem 13.1.** Assume $\mathbb{G}$ to be split over $L$ and $L$ unramified over $\mathbb{Q}_p$. The $G_0$- representation $\mathrm{Ind}_{P_w}^{G_0}(\chi_w)$ possesses an $\mathfrak{n}_w$-character $\theta$ on $S'$, given by

$$\theta(s) = \chi_w(s) \prod_{1 \le j \le d} \frac{1}{\left(1 - a_j^{-1}(s)\right) \left|1 - a_j^{-1}(s)\right|_L}.$$

The dual $D(G_0, K)$-module of $\mathrm{Ind}_{P_w}^{G_0}(\chi_w)$ is given by $M = \left(\mathrm{Ind}_{P_w}^{G_0}(\chi_w)\right)'_b = D(G_0, K) \otimes_{D(P_w, K)} \chi'_w$ $\simeq D(N_w, K) \otimes_K \chi'_w$. According to Proposition 11.6 there is a pair $(H_0, r_0)$ such that for all $(H, r)$ $\ge (H_0, r_0)$ we have

$$M_{(H,r)} \simeq D_{(H \cap N_w, r)}(N_w, K) \underset{K}{\otimes} \chi'_w.$$

With these identifications the subalgebra $D(N_w, K) \subset D(G_0, K)$ simply acts by left multiplication on $M$ and on $M_{(H,r)}$.

Let $(\mathfrak{x}_1, \ldots, \mathfrak{x}_d)$ be a basis of the $o_L$-lattice $\mathrm{Log}(H \cap N_w) \subset \mathfrak{n}_w$ consisting of root space vectors $\mathfrak{x}_i \in \mathfrak{g}_{a_i}$; such a basis exists by Proposition 11.3. Let $(v_1, \ldots, v_n)$ be a $\mathbb{Z}_p$-basis of $o_L$, and put

$$h_{ij} := \mathrm{Exp}(v_i\, \mathfrak{x}_j) \quad (1 \le i \le n,\ 1 \le j \le d).$$

We apply Lemma 11.5 to obtain a sequence of $p$-powers $m = (m_{11}, \ldots, m_{nd})$ such that $(h_{11}^{m_{11}}, \ldots,$ $h_{nd}^{m_{nd}})$ is an ordered basis of $H \cap N_w$ and

$$\left(h_{11}^{\alpha_{11}} \ldots h_{nd}^{\alpha_{nd}}\right)_{0 \le \alpha < m}$$

is a system of representatives of $N_w / (H \cap N_w)$ (where by "$0 \le \alpha < m$" we mean "$0 \le \alpha_{ij} < m_{ij}$ for all $i$, $j$"). In particular,

$$\left(h_{11}^{\alpha_{11}} \ldots h_{nd}^{\alpha_{nd}}\right)_{0 \le \alpha < \epsilon(r) \cdot m}$$

is a system of representatives of $N_w / (H^{\epsilon(r)} \cap N_w)$. Then Corollary 7.4 tells us that the family

$$\left(\mathfrak{x}_1^{\beta_1} \ldots \mathfrak{x}_d^{\beta_d}\, \delta_{h_{11}}^{\alpha_{11}} \ldots \delta_{h_{nd}}^{\alpha_{nd}}\right)_{\alpha < \epsilon(r)\, m,\ \beta \ge 0}$$

is a basis of the $K$-Banach space $M_{(H,r)}$. For $k \in \mathbb{N}$ let $M(k) \subset M_{(H,r)}$ denote the closed $K$- vector subspace generated by the elements

$$\lambda_{\alpha,\beta} := \mathfrak{x}_1^{\beta_1}\, \delta_{h_{11}}^{\alpha_{11}} \ldots \delta_{h_{n1}}^{\alpha_{n1}} \quad \ldots \quad \mathfrak{x}_d^{\beta_d}\, \delta_{h_{1d}}^{\alpha_{1d}} \ldots \delta_{h_{nd}}^{\alpha_{nd}}$$

($\alpha \in \mathbb{N}^{nd}$, $\alpha < \epsilon(r)\, m$, $\beta \in \mathbb{N}^d$, $\beta \not< k$). From the formula $[\mathfrak{g}_a, \mathfrak{g}_b] = \mathfrak{g}_{a+b}$ for all roots $a$, $b$ (where we understand $\mathfrak{g}_{a+b} = 0$ if $a + b$ is not a root) it follows that $[\mathfrak{n}_w, \mathfrak{n}_w]$ is generated by the one-dimensional spaces $\mathfrak{g}_{a+b}$ ($a, b \in \{a_1, \ldots, a_d\}$) and, by iteration, that the $\nu$-th lower central series member $C^\nu(\mathfrak{n}_w)$ is generated by the spaces $\mathfrak{g}_{b_1+\ldots+b_\nu}$ $(b_1, \ldots, b_\nu \in \{a_1, \ldots, a_d\})$. Hence from Lemma 8.5 (ii) and Lemma 8.6, applied to the family

$$\mathbf{g} = \left(h_{1\,j}{}^{\alpha_{1\,j}} \ldots h_{n\,j}{}^{\alpha_{n\,j}}\right)_{1\le j\le d, \alpha<\epsilon(r)\,m},$$

we obtain that the filtrations $(M(k))_{k\in\mathbb{N}}$ and $\left(\mathfrak{n}_w^k\, M_{(H,r)}\right)_{k\in\mathbb{N}}$ are nested into each other. Once we know that the spaces $M(k)$ are stable under $S_0$ we may apply Lemma 2.7 to calculate the limit of the sequence $\left(\mathrm{Ch}\left(M_{(H,r)} \big/ \mathfrak{n}_w^k\, M_{(H,r)}\right)\right)_{k\in\mathbb{N}}$ in $\mathbb{Z}[\![X(S_0)]\!]$ by means of the sequence $(\mathrm{Ch}(M_{(H,r)} / M(k)))_{k\in\mathbb{N}}$.

The action of an element $s \in S_0$ on $M$ is given by

$$s\,\lambda = \chi_w^{-1}(s)\,\lambda^{(s)}$$

where $\lambda^{(s)}$ denotes the distribution $f \mapsto \lambda(f_s)$ and $f_s$ denotes the function $\left(x \mapsto f\left(s\,x\,s^{-1}\right)\right)$ (Lemma 11.1 (i)). In the next two lemmata we calculate $\lambda^{(s)}$ for the elements $\lambda = \lambda_{\alpha,\beta}$ defined above.

For a tuple $\beta \in \mathbb{N}^d$ and $s \in S_0$ write

$$a(s)^\beta := a_1(s)^{\beta_1} \ldots a_d(s)^{\beta_d} \in L^\times;$$

thus we obtain a character $a^\beta = \left(s \mapsto a(s)^\beta\right) \in X(S_0)$.

**Lemma 13.2.** Let $s \in S_0$, $\alpha \in \mathbb{N}^{n\,d}$, $\beta \in \mathbb{N}^d$. Then

$$\lambda_{\alpha,\beta}^{(s)} = a(s)^\beta \left(\prod_{1\le i\le n} \delta_{s\,h_{i1}^{\alpha_{i1}}\,s^{-1}}\right) \mathfrak{x}_1^{\beta_1} \ldots \left(\prod_{1\le i\le n} \delta_{s\,h_{id}^{\alpha_{id}}\,s^{-1}}\right) \mathfrak{x}_d^{\beta_d}.$$

*Proof.* First of all,

$$\delta_x(f_s) = \delta_{s\,x\,s^{-1}}(f)$$

for all $x \in N_w$ and $f \in C^{\mathrm{an}}(N_w, K)$. Next, the linear form $\mathfrak{x}_j$ on $C^{\mathrm{an}}(N_w, K)$ satisfies

$$\mathfrak{x}_j(f) = \lim_{t\to 0} \frac{f\left(\mathrm{Exp}\left(t\,\mathfrak{x}_j\right)\right) - f(1)}{t}$$

(cf. [30], Section 2). But by [7], III.4.4 Corollary 3 of Proposition 8, and since $\mathfrak{x}_j$ is a root vector of weight $a_j$, we calculate $f_s(\mathrm{Exp}(t\,\mathfrak{x}_j) = f\left(s\,\mathrm{Exp}(t\,\mathfrak{x}_j)\,s^{-1}\right) = f(\mathrm{Exp}(\mathrm{Ad}(s)\,(t\,\mathfrak{x}_j))) = f(\mathrm{Exp}(a_j(s)\,t\,\mathfrak{x}_j))$. Hence

$$\mathfrak{x}_j(f_s) = \lim_{t\to 0} \frac{f\big(\mathrm{Exp}\big(a_j(s)\, t\, \mathfrak{x}\big)\big) - f(1)}{t} = a_j(s)\, \mathfrak{x}_j(f).$$

Thus we have shown $\delta_{h_{ij}}^{(s)} = \delta_{s\,x\,s^{-1}}$ and $\mathfrak{x}_j^{(s)} = a_j(s)\, \mathfrak{x}_j$ for all $i$, $j$. The lemma follows now from the convolution formula: for any two distributions $\lambda$, $\mu$ we have

$$\begin{aligned}(\lambda\, \mu)\,(f_s) &= \mu(x \mapsto \lambda(y \mapsto f_s(y\,x))) \\ &= \mu\Big(\big(x \mapsto \lambda\big((y \mapsto f(y\,x))_s\big)\big)_s\Big) \\ &= (\lambda^{(s)}\, \mu^{(s)})\,(f). \ \square\end{aligned}$$

Notations: For any $c \in \mathbb{Z}_p$ we define elements $Q_{ij}(c) \in \epsilon(r)\, m_{ij}\, \mathbb{Z}_p$ and $R_{ij}(c) \in \mathbb{N}_{<\epsilon(r)\, m_{ij}}$ by $c = Q_{ij}(c) + R_{ij}(c)$. For any $a \in o_L$ and any $1 \le i' \le n$ we write $a\, v_{i'} = \sum_{1\le i\le n} b(a)_{i,i'}\, v_i$ with coefficients $b(a)_{i,i'} \in \mathbb{Z}_p$.

**Lemma 13.3.** Let $s \in S_0$, $\alpha \in \mathbb{N}^{n\,d}$, $\beta \in \mathbb{N}^d$ with $\alpha < \epsilon(r)\, m$. Then

$$\lambda_{\alpha,\beta}^{(s)} = a(s)^\beta \sum_{\gamma\in\mathbb{N}^{n\,d}} \frac{1}{\gamma!}\, Q(s,\, \alpha)^\gamma\, \lambda_{R(s,\alpha),\beta+|\gamma|}$$

where $\gamma! := \prod_{i,j} \gamma_{ij}!$, $|\gamma| := (\gamma_{1j} + \ldots + \gamma_{nj})_{1\le j\le d}$, $R(s,\, \alpha) := \big(R_{ij}\big(\sum_{1\le i'\le n} \alpha_{i'\,j}\, b(a_j(s))_{i\,i'}\big)\big)_{1\le i\le n, 1\le j\le d}$, $Q(s,\, \alpha) := \big(Q_{ij}\big(\sum_{1\le i'\le n} \alpha_{i'\,j}\, b(a_j(s))_{i\,i'}\big)\, v_i\big)_{1\le i\le n, 1\le j\le d}$.

*Proof.* We work on with the result of the preceding lemma: Fix $0 \le j \le d$. Then

$$\begin{aligned}\prod_{1\le i'\le n} \delta_{s\, h_{i'j}^{\alpha_{i'j}}\, s^{-1}} &= \prod_{1\le i'\le n} \delta_{s\, \mathrm{Exp}\left(\alpha_{i'\,j}\, v_{i'}\, \mathfrak{x}_j\right) s^{-1}} = \prod_{1\le i'\le n} \delta_{\mathrm{Exp}\left(\mathrm{Ad}(s)\, \left(\alpha_{i'\,j}\, v_{i'}\, \mathfrak{x}_j\right)\right)} \\ = \prod_{1\le i'\le n} \delta_{\mathrm{Exp}\left(a_j(s)\, \alpha_{i'\,j}\, v_{i'}\, \mathfrak{x}_j\right)} &= \prod_{1\le i'\le n} \delta_{\mathrm{Exp}\left(\sum_{1\le i\le n} \alpha_{i'\,j}\, b(a_j(s))_{i,i'}\, v_i\, \mathfrak{x}_j\right)} = \prod_{1\le i\le n} \delta_{\mathrm{Exp}\left(\sum_{1\le i'\le n} \alpha_{i'\,j}\, b(a_j(s))_{i,i'}\, v_i\, \mathfrak{x}_j\right)} \\ &= \prod_{1\le i\le n} \delta_{\mathrm{Exp}\left(R_{ij}\, v_i\, \mathfrak{x}_j\right)}\, \delta_{\mathrm{Exp}\left(Q_{ij}\, v_i\, \mathfrak{x}_j\right)} = \prod_{1\le i\le n} \delta_{h_{ij}}^{R_{ij}}\, \delta_{\mathrm{Exp}\left(Q_{ij}\, v_i\, \mathfrak{x}_j\right)}\end{aligned}$$

where we put $R_{ij} = R_{ij}\big(\sum_{1\le i'\le n} \alpha_{i'\,j}\, b(a_j(s))_{i\,i'}\big)$ and $Q_{ij} = Q_{ij}\big(\sum_{1\le i'\le n} \alpha_{i'\,j}\, b(a_j(s))_{i\,i'}\big)$. Now $\mathrm{Exp}(Q_{ij}\, v_i\, \mathfrak{x}_j) = h_{ij}^{\;Q_{ij}}$ is contained in $H^{\epsilon(r)}$, so according to Remark 7.3 (2) d), $\delta_{\mathrm{Exp}\left(Q_{ij}\, v_i\, \mathfrak{x}_j\right)} = \exp(Q_{ij}\, v_i\, \mathfrak{x}_j) = \sum_{\gamma\in\mathbb{N}} \frac{1}{\gamma!}\, (Q_{ij}\, v_i)^\gamma\, \mathfrak{x}_j^\gamma$. Hence

$$\left(\prod_{1\le i\le n} \delta_{s\, h_{ij}^{\alpha_{ij}}\, s^{-1}}\right) \mathfrak{x}_j^{\beta_j} = \left(\prod_{1\le i\le n} \left(\delta_{h_{ij}}^{R_{ij}} \sum_{\gamma\in\mathbb{N}} \frac{1}{\gamma!}\, (Q_{ij}\, v_i)^\gamma\, \mathfrak{x}_j^\gamma\right)\right) \mathfrak{x}_j^{\beta_j}$$

$$= \sum_{\gamma \in \mathbb{N}^n} \frac{1}{\gamma_1! \ldots \gamma_n!} \left( \prod_{1 \le i \le n} (Q_{ij}\, v_i)^{\gamma_i} \right) \left( \prod_{1 \le i \le n} \delta_{h_{ij}}^{R_{ij}} \right) \mathfrak{x}_j^{\gamma_1 + \ldots + \gamma_n + \beta_j}$$

Combined with the preceding lemma this gives

$$\lambda_{\alpha,\beta}^{(s)} = a(s)^{\beta} \left( \prod_{1 \le i \le n} \delta_{s\, h_{i1}^{\alpha_{i1}}\, s^{-1}} \right) \mathfrak{x}_1^{\beta_1} \ldots \left( \prod_{1 \le i \le n} \delta_{s\, h_{id}^{\alpha_{id}}\, s^{-1}} \right) \mathfrak{x}_d^{\beta_d}$$

$$= a(s)^{\beta} \sum_{\gamma \in \mathbb{N}^{n\,d}} \frac{1}{\gamma!} \left( \prod_{\substack{1 \le i \le n \\ 1 \le j \le d}} (Q_{ij}\, v_i)^{\gamma_{ij}} \right) \left( \prod_{1 \le i \le n} \delta_{h_{i1}}^{R_{i1}} \right) \mathfrak{x}_1^{\gamma_{11} + \ldots + \gamma_{n1} + \beta_1} \ldots \left( \prod_{1 \le i \le n} \delta_{h_{id}}^{R_{id}} \right) \mathfrak{x}_d^{\gamma_{1d} + \ldots + \gamma_{nd} + \beta_d}$$

which is precisely the desired formula. □

**Corollary 13.4.** The space $M(k)$ is stable under the action of $S_0$. □

**Corollary 13.5.** The $S_0$-representation $M_{(H,r)} / M(k)$ is isomorphic to a block matrix representation of $S_0$ of the form

$$\begin{pmatrix} \chi_w^{-1}\, a^{(0,\ldots,0)}\, A & & & 0 \\ & \chi_w^{-1}\, a^{(1,0,\ldots,0)}\, A & & \\ & & \ddots & \\ * & & & \chi_w^{-1}\, a^{(k-1,\ldots,k-1)}\, A \end{pmatrix}$$

where $A$ is the matrix representation of $S_0$ defined by $A(s)_{(\alpha,\alpha')} = \delta_{\alpha',R(s,\alpha)}$ ($\alpha, \alpha' \in \mathbb{N}_{<\epsilon(r)\,m}$; here $\delta$ denotes Kronecker's symbol).

*Proof.* A $K$-basis of $M_{(H,r)} / M(k)$ is given by the elements $\lambda_{\alpha,\beta}$ ($\alpha < \epsilon(r)\,m$, $\beta < k$). Since, in the notations of Lemma 13.3, "$\beta + |\gamma| = \beta$" implies "$\frac{1}{\gamma!}\, Q(s, \alpha)^{\gamma} = 1$" the corollary follows immediately from that lemma. □

We continue to denote by $A$ the $S_0$-representation defined in the preceding corollary. Recall that $S' = S_0 \cap G^{\text{reg}}$.

**Lemma 13.6.** Let $s \in S'$.

(i) $a_j(s) \neq 1$ for all $1 \le j \le d$.

(ii) If $\epsilon(r)\, m_{ij} \ge |a_j(s) - 1|_L^{-1}$ for all $1 \le i \le n$, $1 \le j \le d$ then

$$\mathrm{tr}(A(s)) = |a_1(s) - 1|_L^{-1} \dots |a_d(s) - 1|_L^{-1}.$$

*Proof.* (i) follows from the description of $G^{\mathrm{reg}}$ given in [34], p. 197; or from [3], Lemma 12.2.

(ii) The trace of $A(s)$ is equal to the number of fixpoints of the permutation $\alpha \mapsto R(s, \alpha)$ of the set $\mathbb{N}^{n\,d}_{<\epsilon(r)\,m}$. We observe

$$\begin{aligned} & \alpha = R(s, \alpha) \\ \Leftrightarrow \quad & \alpha_{i\,j} = R\left(\textstyle\sum_{1\le i'\le n} b(a_j(s))_{i\,i'}\, \alpha_{i'\,j}\right) \text{ for all } i,\ j \\ \Leftrightarrow \quad & v_p\left(\textstyle\sum_{1\le i'\le n} b(a_j(s))_{i\,i'}\, \alpha_{i'\,j} - \alpha_{i\,j}\right) \ge v_p(\epsilon(r)\, m_{i\,j}) \text{ for all } i,\ j. \end{aligned}$$

But

$$v_p\Big(\det\Big(\big(b(a_j(s))_{i,i'}\big)_{1\le i,i'\le n} - 1\Big)\Big) = v_p\big(N_{L|\mathbb{Q}_p}(a_j(s) - 1)\big) = v_L(a_j(s) - 1)$$

(cf. [22], I. (2.2) and II. (4.8); recall that $L\,|\,\mathbb{Q}_p$ is unramified). Since we assume $v_p(\epsilon(r)\,m_{i\,j}) \ge v_L(a_j(s) - 1)$ there are, for each $j$, precisely $p^{v_L(a_j(s)-1)}$ tuples $(\alpha_{1\,j}, \dots, \alpha_{d\,j})$ in $\mathbb{N}_{<\epsilon(r)\,m_{1\,j}} \times \dots \times \mathbb{N}_{<\epsilon(r)\,m_{d\,j}}$ satisfying the above condition (independently of the particular values of $m_{i\,j}$ and $\epsilon(r)$). Consequently the number of fixpoints equals $\prod_j p^{v_L(a_j(s)-1)} = \prod_j |a_j(s) - 1|_L^{-1}$.

**Lemma 13.7.** (i) The representations $A$ and $M_{(H,r)} \big/ \mathfrak{n}^k\, M_{(H,r)}$ of $S_0$ $(k \in \mathbb{N})$ are finitely trigonalisable over $C$. The sequence of formal characters $\big(\mathrm{Ch}\big(M_{(H,r)} \big/ \mathfrak{n}^k\, M_{(H,r)}\big)\big)_{k\in\mathbb{N}}$ converges in $\mathbb{Z}[\![X(S_0)]\!]$ to the element

$$\Phi_{(H,r)} := e\big(\chi_w^{-1}\big) \cdot \mathrm{Ch}(A) \cdot \sum_{\beta\in\mathbb{N}^d} e\big(a^\beta\big).$$

(ii) $\Phi_{(H,r)}$ is evaluable on $S'$.

(iii) Let $s \in S'$. If $\epsilon(r)\, m_{i\,j} \ge |a_j(s) - 1|_L^{-1}$ for all $1 \le i \le n$, $1 \le j \le d$ then

$$\mathrm{ev}_{S'}(\Phi_{(H,r)})\,(s) = \chi_w(s)^{-1} \prod_{1\le j\le d} \frac{1}{(1 - a_j(s))\, |1 - a_j(s)|_L}.$$

*Proof.* (i) We have seen that the decreasing sequences $M(0) \supset M(1) \supset M(2) \supset \dots$ and $\mathfrak{n}^0\, M_{(H,r)} \supset \mathfrak{n}^1\, M_{(H,r)} \supset \mathfrak{n}^2\, M_{(H,r)} \supset \dots$ are nested into each other. It is clear that the finite-dimensional $S_0$-representations $A$, $M_{(H,r)}\,/\,M(k)$, $M_{(H,r)} \big/ \mathfrak{n}^k\, M_{(H,r)}$ are finitely trigonalisable over $C$. According to Corollary 13.5, and using Corollary 2.6, we have

$$\mathrm{Ch}(M_{(H,r)}/M(k)) = e\left(\chi_w^{-1}\right)\cdot \mathrm{Ch}(A)\cdot \sum_{\beta<k} e\left(a^{\beta}\right),$$

hence

$$\lim_{k\to\infty} \mathrm{Ch}(M_{(H,r)}/M(k)) = e\left(\chi_w^{-1}\right)\cdot \mathrm{Ch}(A)\cdot \sum_{\beta\in\mathbb{N}^d} e\left(a^{\beta}\right).$$

By Lemma 2.7 this is also the limit of the sequence $\left(\mathrm{Ch}\left(M_{(H,r)}\big/\mathfrak{n}^k M_{(H,r)}\right)\right)_{k\in\mathbb{N}}$.

(ii) In $\mathbb{Z}[\![X(S_0)]\!]$ we have the equality

$$\begin{aligned} e(1) &= (e(1)-e(a_1))\left(\textstyle\sum_{\beta_1\in\mathbb{N}} e(a_1)^{\beta_1}\right)\ldots(e(1)-e(a_d))\left(\textstyle\sum_{\beta_d\in\mathbb{N}} e(a_d)^{\beta_n}\right) \\ &= (e(1)-e(a_1))\ldots(e(1)-e(a_d))\textstyle\sum_{\beta\in\mathbb{N}^d} e\left(a^{\beta}\right). \end{aligned}$$

Hence $(e(1)-e(a_1))\ \ldots\ (e(1)-e(a_d))\ \Phi_{(H,r)} = e\left(\chi_w^{-1}\right)\mathrm{Ch}(A)$, where $e\left(\chi_w^{-1}\right)\mathrm{Ch}(A)$ is contained in $\mathbb{Z}[X(S_0)]$ and, regarding Lemma 13.6 (i), $(e(1)-e(a_1))\ \ldots\ (e(1)-e(a_d))$ is contained in the set $\mathcal{S}_{S'}$ (cf. Section 1). This means that $\Phi_{(H,r)}$ is evaluable on $S'$.

(iii) By Remark 1.2 (iii) we may calculate the respective values of $e\left(\chi_w^{-1}\right)$, $\mathrm{Ch}(A)$, $\sum_{\beta\in\mathbb{N}^d} e\left(a^{\beta}\right)$ in $s\in S'$ separately. In case of $e\left(\chi_w^{-1}\right)$ and $\mathrm{Ch}(A)$ these are simply the usual character values in $s$ (Remark 2.2 (ii)), which means in the latter case: the trace of the finite dimensional operator $A(s)$ as calculated in Lemma 13.6. Finally, the proof of part (ii) above shows that $\mathrm{ev}\left(\sum_{\beta\in\mathbb{N}^d} e\left(a^{\beta}\right)\right)(s) = (1-a_1(s))^{-1}\ \ldots\ (1-a_d(s))^{-1}$. Multiplication of these values gives the asserted formula. □

**Corollary 13.8.** The sequence of functions $(\mathrm{ev}_{S'}(\Phi_{(H,r)}))_{(H,r)\in\Xi}$ converges weakly to the function

$$S'\to K,\ s\mapsto \chi_w(s)^{-1}\prod_{1\le j\le d}\frac{1}{\left(1-a_j(s)\right)|1-a_j(s)|_L}.$$

*Proof.* For every $s\in S'$ there certainly exists a pair $(H_1, r_1)\ge(H_0, r_0)$ satisfying the hypothesis of part (iii) of Lemma 13.7. If $(H, r)\ge(H_1, r_1)$ then $H^{\epsilon(r)}\subset {H_1}^{\epsilon(r_1)}$ (Lemma 9.2), and so the pair $(H, r)$ satisfies that hypothesis as well. It remains to observe that the formula in Lemma 13.7 (iii) does not depend on the particular pair $(H, r)$ provided that hypothesis is satisfied. □

This completes the proof of Theorem 13.1.

## References

[1] Amice, Y.: Interpolation $p$-adique, *Bull. Soc. math. France 92*, 1964, pp. 117-180

[2] Bernstein, J.: *Representations of p-adic groups*, lecture notes written by K. E. Rumelhart, 1992, available at `http://www.math.uchicago.edu/~mitya/langlands.html`

[3] Borel, A.: *Linear Algebraic Groups*, 2nd ed., Springer, 1991

[4] Bourbaki, N.: *Algebra*, Chap. 1-3. Springer, 1989

[5] Bourbaki, N.: *Algèbre*, Chap. 4-7. Masson, Paris, 1981

[6] Bourbaki, N.: *General topology,* Chap. 1-4., Springer, 1988

[7] Bourbaki, N.: *Lie groups and Lie algebras*, Chap. 1-3. Springer, 1989

[8] Bourbaki, N.: *Groupes et algèbres de Lie*, Chap. 4-6. Masson, Paris, 1981

[9] Bourbaki, N.: *Varietes differentielles et analytiques. Fascicule de resultats*, Hermann, Paris, 1967

[10] Cartier, P.: Representations of $\mathfrak{p}$-adic groups: A survey, *Automorphic forms, representations and L-functions, Proc. symp. pure math. 33,* AMS, 1977, pp. 111-165

[11] Demazure, M., Gabriel, P.: *Groupes Algébriques*, Masson, Paris / North-Holland, Amsterdam, 1970

[12] Dixon, J.D., Du Sautoy, M.P.F., Mann, A., Segal, D.: *Analytic pro p-groups* , 2nd Edition, Cambridge University Press 1999

[13] Dixmier, J.: *Enveloping algebras,* AMS 1996

[14] Emerton, M.: Locally analytic vectors in representations of locally $p$-adic analytic groups, *preprint*, 2004, available at `www.math.northwestern.edu/~emerton/preprints.html`

[15] Féaux de Lacroix, C.T.: Einige Resultate über die topologischen Darstellungen $p$- adischer Liegruppen auf unendlich dimensionalen Vektorräumen über einem $p$- adischen Körper, *Schriftenreihe Math. Inst, Univ. Münster,* 3. Serie, Heft 23, 1999

[16] Frommer, H.: The locally analytic principal series of split reductive groups, *Preprintreihe SFB 478 Münster*, Heft 265, 2003

[17] Gel'fand, M., Graev, M. I., Pyatetskii-Shapiro, I. I.: *Representation theory and automorphic forms,* Academic Press, 1990

[18] Goodman, R., Wallach, N.: *Representations and invariants of the classical groups,* Cambridge University Press, 1998

[19] Harish-Chandra, A submersion principle and its applications, *Papers dedicated to the memory of V. K. Patodi*, Indian Academy of Sciences, Bangalore, and the Tata Institute of Fundamental Research, Bombay, 1980, pp. 95-102; also available in: Harish-Chandra, *Collected Papers Volume IV,* Springer, 1984, pp. 439-446

[20] Kohlhaase, J.: Invariant distributions on *p*-adic analytic groups, *Duke Math. J. 137*, 2007, pp 19-62.

[21] Lazard, M.: Groupes analytiques *p*-adiques. *Publ. math. I.H.E.S. 26*, 1965, pp. 5-219

[22] Neukirch, J.: *Algebraische Zahlentheorie.* Springer, 1992.

[23] Orlik, S., Strauch, M.: A. J.: On the irreducibility of locally analytic principal series representations, *Preprint-Reihe des Mathematischen Instituts der Universität Leipzig*, Heft 06/2006

[24] van Rooij, A.C.M.: *Non-archimedean functional analysis.* Marcel Dekker, New York, 1978

[25] Schikhof, W.H.: *Ultrametric calculus*, Cambridge University Press, 1984

[26] Schmidt, T.: Auslander regularity of *p*-adic distribution algebras, *Preprint-Reihe des SFB 478 Münster*, Heft 451, 2006

[27] Schmidt, T.: *Analytic vectors in continuous p-adic representations*, preprint, 2008.

[28] Schneider, P.: *Nonarchimedean functional analysis*, Springer, 2002.

[29] Schneider, P., Teitelbaum, J.: Locally analytic distributions and *p*- adic representation theory, with applications to $\mathrm{GL}_2$, *Journ. AMS 15*, 2002, pp. 443-468

[30] Schneider, P., Teitelbaum, J.: $U(\mathfrak{g})$-finite locally analytic representations, *Representation theory 5*, 2001, pp. 111-128

[31] Schneider, P., Teitelbaum, J.: *p*-adic Fourier theory, *Documenta Math. 6*, 2001, pp. 447- 481

[32] Schneider, P., Teitelbaum, J.: Algebras of *p*- adic distributions and admissible representations, *Invent. math. 153*, 2003, pp. 145-196

[33] Schneider, P., Teitelbaum, J.: Duality for admissible locally analytic representations, *Representation theory 9*, 2005, pp. 297-326

[34] Silberger, A. J.: *Introduction to harmonic analysis on reductive p-adic groups*, Mathematical notes, Princeton University Press, 1979

[35] Tits, J.: Reductive groups over local fields, *Automorphic forms, representations and L-functions, Proc. symp. pure math. 33*, AMS, 1977, pp. 29-69

Ralf Diepholz

Westfälische Wilhelms-Universit ät Münster, Einsteinstrasse 62, 48161 Münster, Germany
ninteman@uni- muenster.de